\documentclass{amsart}
\usepackage{amssymb,latexsym,graphics}
\usepackage[mathscr]{eucal}

\newtheorem{thm}{Theorem}[section]

\newtheorem{lem}[thm]{Lemma}
\newtheorem{prop}[thm]{Proposition}

\theoremstyle{remark}
\newtheorem{rem}[thm]{Remark}

\theoremstyle{definition}

\newtheorem{notation}[thm]{Notation}
\theoremstyle{definition}

\newcommand{\lp}[2]{\Vert \, #1 \, \Vert_{#2}}
\newcommand{\td}[1]{\widetilde{#1}}

\newcommand{\bL}{\underline{L}}

\newcommand{\dlambda}{{\dot{\lambda}}}
\def\pa{\partial}

\newcommand{\ret}{\vspace{.3cm}}

\begin{document}

\title[Singularities in the Critical $O(3)$ $\sigma$-Model]
{On the Formation of Singularities in the Critical $O(3)$ $\sigma$-Model}
\author{Igor Rodnianski}
\address{Department of Mathematics, Princeton University}
\author{Jacob Sterbenz}
\address{Department of Mathematics, University of California, San Diego}
\thanks{This work was conducted when I.R. was visiting the Department of Mathematics, MIT. He was
also partially supported by the NSF grant DMS-0406627. J.S. was supported by an NSF postdoctoral
fellowship.}
\date{}
\commby{}


\begin{abstract}
We study the phenomena of  energy concentration for the critical  $O(3)$ sigma model, also known
as the wave map flow from $\mathbb{R}^{2+1}$ Minkowski space into the sphere
$\mathbb{S}^2$. We establish rigorously and constructively existence of a set of smooth
initial data resulting in a dynamic finite time formation of singularities. The construction and
analysis is done in the context of the
k-equivariant  symmetry reduction, and we restrict to  maps with
homotopy class  $k\geqslant 4$. The concentration mechanism we uncover is
essentially  due to a
resonant self-focusing (shrinking) of a corresponding harmonic map. We show that the phenomenon
is generic (e.g. in
certain Sobolev spaces) in that it persists under small perturbations of initial data,
while the resulting blowup is bounded by a
log-modified self-similar asymptotic.
\end{abstract}

\maketitle

\section{Introduction}

One of the simplest non-trivial models of Quantum Field Theory is
based on the $(2+1)$ dimensional Lorentz invariant $O(3)$ classical
$\sigma$-model. It is a nonlinear scalar field Lagrangian theory for
a map $\Phi:\,\mathbb{R}^{2+1}\to \mathbb{S}^2\subset \mathbb{R}^3$
with the Lagrangian density:
\begin{equation}
        \mathcal{L}[\Phi] \ = \frac 12 \partial_\alpha \Phi\cdot
        \partial_\beta \Phi \ m^{\alpha\beta} \ , \label{WM_L}
\end{equation}
where $m_{\alpha\beta}$ is the Minkowski metric on $\mathbb{R}^{2+1}$.
Evolution of the nonlinear scalar field
$\Phi$ is described by the Euler-Lagrange equations:
\begin{equation}
    \Box\Phi=-\Phi \left (\pa^\alpha\Phi\cdot \pa_\alpha
    \Phi\right) \ . \label{eq:EL}
\end{equation} \\

The equation \eqref{eq:EL} belongs to the more general class of
``wave-map'' problems, in which $\Phi$ is a map from Minkowski space
$\mathbb{R}^{2+1}$ to a Riemannian manifold $(\mathcal {M}, g)$. The
map $\Phi$ is a solution of the Euler-Lagrange equations:
\begin{equation}
    D^\alpha\pa_\alpha\Phi=0 \ , \label{basic_WM}
\end{equation}
corresponding to the Lagrangian density:
\begin{equation}
        \mathcal{L}[\Phi] \ = \ \frac{1}{2}
    \ g_{ij}\  \partial_\alpha \Phi^i
        \partial_\beta \Phi^j \ m^{\alpha\beta} \ . \label{WM_action}
\end{equation}
Here $\{\Phi^i\}$ denote local coordinates on
$\mathcal{M}$, which in turn (under the map) depend on the
Minkowski variables $\{x^\alpha\}_{\alpha=0,1,2}$.
$D$ is the pullback of the Levi-Civita
connection to the (trivial) bundle $\Phi^*( T\mathcal{M})$.
In terms of
the local coordinates $\{\Phi^i\}$ this pull-back connection acting
on sections of $\Phi^*( T\mathcal{M})$ reads:
\begin{align}
        D_\alpha \ &= \ \partial_\alpha + \overline{\Gamma}_{\alpha
        j}^k \ , &\overline{\Gamma}_{\alpha j}^k \ &= \
        \Gamma_{ij}^k(\Phi) \partial_\alpha \Phi^i \ , \label{frame_chris}
\end{align}
where $\Gamma_{ij}^k$ is the Christoffel symbol in the
coordinates $\{\Phi^i\}$. The wave-map equation
\eqref{basic_WM} then has the intrinsic form:
\begin{equation}
        \partial^\alpha \partial_\alpha \Phi^k \ = \ -
        \Gamma_{ij}^k(\Phi)\, \partial^\alpha \Phi^i \partial_\alpha \Phi^j \ .
        \label{basic_WM2}
\end{equation}\ret

The goal of this paper is to establish and rigorously analyze a catastrophic
instability in the $(2+1)$ dimensional $O(3)$ $\sigma$-model represented by the
equation \eqref{eq:EL}. We will exhibit a spontaneous and monotonic
self-focusing mechanism responsible for a dynamic formation of singularities
for a rather large  and stable set of initial data. This will be done through an entirely
explicit and  constructive\footnote{This is in contrast
to some of the examples of the focussing nonlinear Schr\"odinger and wave equations,
where a finite time  blow-up can be shown by non-constructive
arguments (see \cite{G_NLS} and \cite{L_NLW}). See however the work of Martel/Merle
\cite{MM} on the critical KdV problem for an example of a constructive finite time blow up mechanism.}
description of this phenomena. Our basic result is as follows:\\

\begin{thm}
For every $0 < \epsilon\ll 1$ and $4\leqslant k$ there exists  a set of smooth initial data
$(\Phi_0,\dot{\Phi}_0)\in (\mathbb{S}^2, T\mathbb{S}^2)$ with energy $E=4\pi k + \epsilon^2$
(i.e. the Dirichlet energy defined below), and a finite time $T^{**}=T^{**}(\Phi_0,\dot{\Phi}_0)$,
such that the corresponding solution $\Phi(t,x)$ of
problem \eqref{eq:EL}  remains smooth on the interval $[0,T^{**})$ and
develops a singularity at $T^{**}$. More specifically, there exists a (smooth) decomposition
$\Phi=\underline{\Phi} + \mathcal{R}$, such that as $t\to T^{**}$
we have that for any large $0<M$ an $L^\infty$
bound of the form\footnote{More precise asymptotic
behavior in terms of the energy concentration will be given below, including
both upper and lower bounds.}:
\begin{equation}
    \frac{M}{(T^{**}-t)} \ \leqslant \
    \sup_{x\in\mathbb{R}^2} \, |\nabla_x \underline{\Phi}|
    \leqslant   \frac{\sqrt {|\ln (T^{**}-t)|}}{T^{**}-t}\  \ , \label{rough_lower_blow}
\end{equation}
as well as a uniform bound on the energy  of the remainder:
\begin{equation}
    E[\mathcal{R}] \ \lesssim \ \epsilon^2 \ . \label{remainder_bound}
\end{equation}
Furthermore, sufficiently small equivariant perturbations of $(\Phi_0,\dot{\Phi}_0)$
also lead to blowup with the bounds \eqref{rough_lower_blow}--\eqref{remainder_bound}.
\end{thm}\ret\ret

The problem of a finite time breakdown of solutions of the problem \eqref{eq:EL}
has been a subject of intense study. From a purely analytical perspective,
the context is the global regularity theory for the general wave-map equations
\eqref{basic_WM}, where it is suspected that the formulation of singularities is
ultimately tied to certain convexity properties of the target manifold $\mathcal{M}$.\\

From a more physical or gauge theoretic perspective, and in a
specific context of the $O(3)$ model, the issue of possible
singularity development is thought to be connected to the
\emph{incompleteness} of a certain moduli space which characterizes
the associated static solutions, and provides an approximation for
the dynamical evolution through the so called {\it geodesic hypothesis}. Due to its analytic and physical
interest, the equation \eqref{eq:EL} has
also been a popular subject of numerical and heuristic studies which all universally pointed
in the direction of singularity formation.\\

Before proceeding with a more detailed description of our main
result, we believe it is useful to give a more thorough description
of these various points of view. This begins with a discussion of
the static solutions of \eqref{eq:EL}, that is harmonic maps into
the sphere. Historically, one of the primary motivating factors of
interest in the $O(3)$ $\sigma$--model was due to the richness of
the set of its static solutions. An ingenious procedure of Belavin
and Polyakov \cite{BP_sigma} allows one to find these solutions in a
given homotopy class characterized by the topological degree:
\begin{equation}
    k\ =\ \frac 1{4\pi}\,
    \int_{\mathbb {R}^{2}} \Phi_* \left (dA_{\mathbb {S}^2}\right)\ , \notag
\end{equation}
as solutions of the \emph{first order} ``Bogomol'nyi equations'' (see \cite{B_stab}).
To realize this, one factors the energy functional:
\begin{equation}
    V[\Phi]\ =\ \frac{1}{2} \int_{\mathbb{R}^2} \nabla_x\Phi
    \cdot\nabla_x\Phi\ dx \ , \label{var_energy}
\end{equation}
which is the potential part of the Dirichlet type energy:
\begin{equation}
    E[\Phi](t)\ =\ \frac{1}{2}
    \, \int_{\mathbb{R}^2} \left(\pa_t\Phi\cdot\pa_t \Phi+ \nabla_x\Phi\cdot\nabla_x\Phi\right) \ dx
    \ = \  T[\Phi](t) + V[\Phi](t) \ , \label{eq:Energy}
\end{equation}
associated with the action of the Lagrangian \eqref{WM_L}:
\begin{equation}
    \int_{\mathbb{R}^{2+1}}\ \mathcal{L}[\Phi] \ dx\, dt \ = \ -\,
    \int_{\mathbb{R}}\ \Big( T[\Phi] -  V[\Phi]\Big)(t)\
    dt \ . \label{the_action}
\end{equation}
Using the notation $\epsilon_{ij}$ for the antisymmetric tensor on
two indices, this factorization reads:
\begin{align}
  \begin{split}
    V[\Phi]\ &= \ \frac{1}{4}\,
    \int_{\mathbb{R}^2}\ \left[(\pa_i\Phi\pm \epsilon_{i}^{\ \, j}\Phi\times\pa_j\Phi)
    \cdot (\pa^i\Phi\pm \epsilon^{ij}\Phi\times\pa_j\Phi)\right]\  dx\ \\
    &\ \ \ \ \ \ \ \ \ \ \ \ \ \ \ \ \ \ \ \ \ \ \ \ \ \ \ \ \ \
    \ \ \ \ \ \  \ \ \ \pm \
    \frac{1}{2}\, \int_{\mathbb{R}^2} \epsilon^{ij} \Phi\cdot (\pa_i \Phi\times
    \pa_j\Phi)\ dx \ ,
  \end{split}\label{general_lines}\\
    &=\ \frac{1}{4}\, \int_{\mathbb{R}^2}\
     \left[(\pa_i\Phi\pm \epsilon_{i}^{\ \, j}\Phi\times\pa_j\Phi)
    \cdot (\pa^i\Phi\pm \epsilon^{ij}\Phi\times\pa_j\Phi)\right]\  dx
     \ \pm \ 4\pi k \ . \notag
\end{align}
from which it is more or less immediate that an absolute minimum
of the energy functional $V[\Phi]$ in a given topological sector $k$ must
be a solution of the equation:
\begin{equation}
    \pa_i \Phi \pm \epsilon_{i}^{\ \, j} \Phi\times
    \pa_j \Phi \ = \ 0 \ . \label{eq:Bog}
\end{equation}
In terms of complex coordinates on $\mathbb{R}^2$ and
$\mathbb{S}^2$, the identities \eqref{eq:Bog} are seen to be nothing
other than the Cauchy-Riemann and conjugate Cauchy-Riemann equations
(this is a general phenomena, see \cite{H_int}). Therefore, the
moduli space $M_k$ of static energy minimizing solutions of \eqref{eq:EL}
in a homotopy sector $k$ can be identified with the rational
maps (in $z$ or $\bar{z}$ using complex variables)
$I:\mathbb{C}\to\mathbb{C}$ with degree $k$. Of particular
importance to us will be the \emph{$k$-equivariant} static solutions
(of positive polarity) which are defined via the relation
$I(e^{i\theta}z) = e^{i\theta k} I(z)$. We will label such solutions
by $I^k$ and refer to them as \emph{solitons}.\\

Having uncovered the structure of the space of static minimizing solutions, one is then led to
the ``geodesic'' ansatz alluded to above
for the approximate dynamics of time-dependent solutions  (see also \cite{Manton}).
To understand this, the first thing to notice is that  minimizers of the variational
problem \eqref{var_energy} are highly degenerate. Specifically, they
are invariant under the full conformal
group of linear fractional transformations acting on $\mathbb{C}$. If one restricts
to $k$-equivariant solutions, then most of this symmetry is lost, and the only remaining
degree of freedom which fixes the energy \eqref{var_energy} is the scaling transformations:
\begin{equation}
    I^k(t,x)\ \to \ I^k_\lambda(t,x)\ =\ I^k(\lambda t, \lambda x) \ . \label{sc_trans}
\end{equation}
Based on this, one would expect that the path of least resistance according to the action
\eqref{the_action} would be for (symmetric) solutions close to the
family of static solutions to ``slide'' along the moduli space $M_k$
via the transformations  \eqref{sc_trans}. That is, for a fully dynamic
solution $\Phi$ sufficiently close to some $I^k$, i.e. $E[\Phi]=4\pi k + \epsilon^2$,
there should be a splitting\footnote{This is very similar to what is
done in  the modulational theory of dispersive solitons, and we will expound on
this in much more detail in the sequel.} as follows:
\begin{equation}
    \Phi(t,x) \ = \ I^k\big(\lambda(t) x\big) + \{\hbox{small error}\} \ ,
    \label{intro_decomp}
\end{equation}
and the goal is to understand the \emph{lower dimensional} dynamics
of the parameter $\lambda(t)$. Plugging the ansatz
\eqref{intro_decomp} in the action \eqref{the_action} yields the
following effective Lagrangian\footnote{It must be kept in mind that
this calculation is purely heuristic as, the original Lagrangian
\eqref{WM_L} itself is only a \emph{formal} way to derive the
equations \eqref{eq:EL}.} for $\lambda(t)$:
\begin{equation}
    \underline{\mathcal{L}}[\lambda](t) \ = \ \mathcal{C}_k \frac{\dlambda^2}{\lambda^4}(t)
      + 4\pi k + \{\hbox{small error}\} \ , \label{eff_lang}
\end{equation}
where $\dlambda=\frac{d\lambda}{dt}$ and
where the normalization constant $\mathcal{C}_k$ is given by the explicit integral
(note that this is only finite for $2\leqslant k$):
\begin{equation}
    \mathcal{C}_k \ = \ -\, \frac{1}{2}\,
    \int_{\mathbb{R}^2}\ \lp{r\partial_r I^k}{}^2\ dx  \ . \notag
\end{equation}
Here $\lp{\cdot}{}^2$ is the norm on $\mathbb{R}^3$. The effective dynamics generated by
\eqref{eff_lang} are now given by the formula:
\begin{equation}
    -\, \frac{d^2}{dt^2}\big(\lambda^{-1}\big) \
    = \ \frac{d}{dt}\left(\frac{\dlambda}{\lambda^2}\right) \ = \
    \{\hbox{small error}\} \ . \label{eff_dynamics}
\end{equation}
If one were to ignore the contribution on the right hand side of
this last equation, then the evolution generated by \eqref{eff_lang} would
imply that dynamically the soliton radius collapses as a linear
function of time, or equivalently that $\lambda\sim(T^{**}-t)^{-1}$
for some fixed $T^{**}$. In this sense,
the moduli space $M_k$ is said to be incomplete.\\

While the above scenario is appealing for its simple geometric and physical
motivation, it has been rigorously known for some time that it cannot be quite correct.
This is due to the fundamental regularity results of Shatah and Tahvildar-Zadeh
\cite{ST_WM1}--\cite{ST_WM2} (see also \cite{CS_WM1}--\cite{CS_WM2}) which rules
out the existence of purely (i.e. linear) self-similar collapse:\\

\begin{thm}[Regularity theory for symmetric wave-maps]\label{C_th}
Let $\Phi$ be an equivariant solution to the equation \eqref{eq:EL} with smooth
Cauchy data. There exists an $\epsilon > 0$ with the following property:
Let $T^*$ be any time such that this solution is
$C^\infty$ for all times $0\leqslant t < T^*$ and such that the
following condition holds:\\

\begin{itemize}
    \item For any $0\leqslant t < T^*$ the energy content
    $E_{\mathcal{B}(t)}[\Phi(t)]$  inside the ball  $\mathcal{B}(t)$
    centered at the origin $r=0$,
    is such that\ $\lim_{t\to T^*}E_{\mathcal{B}(t)}[\Phi(t)]\leqslant \epsilon$
    whenever $|\mathcal{B}(t)| = (T^*-t)\cdot o(1)$.
\end{itemize}\ret

Then the wave-map $\Phi$ extends past $T^*$ as a $C^\infty$
solution. That is, if the energy $E[\Phi(t)]$ concentrates at most
at a (linear) self-similar rate up to time $T^*$, then the solution
cannot break down at time $T^*$.
\end{thm}\ret

This theorem shows that the error terms on the right hand side of
\eqref{eff_dynamics} \emph{cannot} be ignored, and that any complete
theory of how $\lambda(t)$ should evolve must take them into
account. In fact, the above result  leaves the question of breakdown
in finite time for the equation \eqref{eq:EL} open to a much wider
range of possibilities because while it gives a necessary lower
bound on any possible blowup rate for $\lambda(t)$, it does not
give \emph{any} upper bound in case collapse might occur.\\

At this point we should further mention  that the general wave-map
equations \eqref{basic_WM} have also been  studied intensely
from an analytic perspective. For the static case of
\eqref{basic_WM} we point out the references \cite{H_moving} and
Chapter 8 of \cite{Jost}, and the references therein. In the case of
dynamic solutions and the Cauchy problem, the only general
understanding of the equations \eqref{basic_WM} that is yet
available is for the local theory (see \cite{KS}) and the small data scale
invariant (global) results of \cite{Tat_WM1} (in the Besov case)
and \cite{Tao_WM1}, \cite{J_WM1}, and \cite{Tat_WM2} for the case of small energy.
Large data global regularity has been conjectured in the case of a hyperbolic space $\mathbb{H}^2$
target, while  singularity formation has
been expected for the $O(3)$ $\sigma$-model for some time. We explain this more in a moment.
We point out to the reader that this is in stark
contrast to what is known for the parabolic analog of
\eqref{basic_WM} (i.e. the harmonic map heat flow), where the global
regularity theory at all energy levels is much better understood
(see \cite{ES_original}, \cite{CDY_blowup},  \cite{S_heat}).\\

However, in the case where the dynamic solutions of \eqref{basic_WM}
possess a large amount of symmetry, there has been considerable
progress toward our understanding of the general Cauchy problem in the case of
arbitrarily large initial data. The global regularity question was
first handled in the work of  Shatah/Tahvildar-Zadeh (see \cite{ST_WM1}--\cite{ST_WM2})
and Christodoulou/Tahvildar-Zadeh (see \cite{CS_WM1}--\cite{CS_WM2}),
where the context is spherical
symmetry or more generally k-equivariance. There is also the
important and closely related work of Struwe (see
\cite{S_WM1}--\cite{S_WM3}), where breakdown is studied in the
general (symmetric) case including maps into $\mathbb{S}^2$. From
all these works, it is known that if the target
manifold $\mathcal{M}$ is ``geodesically convex'', then symmetric
solutions to \eqref{basic_WM} cannot break down in finite time.
Furthermore, this behavior has been shown to be stable under small rough
perturbations in  the recent work of \cite{J_WM2}. Finally, it is known in general (i.e. without
geodesic convexity) that if a symmetric solution to \eqref{basic_WM}
does break down in finite time, then the singularity formation must
be tied to the existence of a static solution to \eqref{basic_WM},
and in fact will rescale to a non-trivial harmonic map in the
limit.\\

While the works mentioned above furnish a great deal of
understanding, they also leave completely open the issue of whether
or not singularities do in fact form in the specific case of dynamic
solutions to the equation \eqref{eq:EL}. The most convincing
evidence to date that breakdown does occur in finite time is the
analytic work of R. C\^ote \cite{C_inst} on strong asymptotic
instability in the energy space, and  the many numerical studies
that have been performed (see for example \cite{B_WM}, \cite{IL},
\cite{LPZ_stab}, \cite{LS}, and \cite{PZ_shrink}). We mention here
that the work \cite{LS} suggests a universal $\log$-modified
self-similar behavior similar to
\eqref{rough_lower_blow}.\\

In this work we show that singularities will form in finite time for
the critical $O(3)$ $\sigma$-model in such a way that the stable dynamics
is bounded by a $\log$-modified self-similar collapse
that is not so far\footnote{In terms of power law type behavior.}
from what is predicted by \eqref{eff_dynamics}.  One of the major points of this
paper is to uncover the precise analytic mechanism which is
responsible for this upper bound.\\

Before closing this subsection, let us make several remarks.
The
first is that  the $O(3)$ $\sigma$-model also enjoys many analogies
with other more complicated field theories such as the $(4+1)$
Yang-Mills and the $(3+1)$ Yang-Mills-Higgs equations. For this
reason, the it has  been an important testing
ground for ideas concerning the structural behavior of these more
complicated models. We would like mention here the work of Bizon,
Ovchinnikov, and Sigal for the case of Yang-Mills instantons
\cite{BS_YM}, which proposes a collapse scenario for the critical
$(4+1)$-dimensional Yang-Mills equations similar to what we deal
with here. The reader will see that some of our methods are
inspired by certain calculations performed in that paper.\\

Secondly, existence of finite time blow-up solutions had been
had been known for some time in the case of a super-critical higher
dimensional wave map problem with Minkowski space $\mathbb{R}^{n+1}$
with $n>2$ as a base and a rotationally symmetric Riemannian manifold
${\mathcal M}$ as a target manifold. The construction of blow-up solutions
is based on existence of $k$-equivariant {\it self-similar} solutions of finite
energy for the higher dimensional wave map problem. Such solutions have been
exhibited in the work of Shatah (see \cite{Sh}) for the
$\mathbb{R}^{3+1}\to\mathbb {S}^3$ problem. This was later
extended to other target manifolds in \cite{ST_WM2} and higher dimensions
$n\ge 4$ in \cite{CST}. In the latter work it was also shown that for $n\ge 7$
self-similar blow up
can occur even in the case when the target manifold is negatively curved.\\

Thirdly, an interesting issue that we would like to draw the
readers attention to
here is that the k-equivariant heat flow corresponding
to the instance of \eqref{eq:EL} we study here is known to be
\emph{globally regular}  (it
is expected that the corresponding Schr\"odinger flow is also
globally regular). That is, for the equivariant maps into the
sphere $\mathbb{S}^2$ with the homotopy index $k>1$, the
harmonic map heat  flow \emph{does not} break down
in finite time \cite{SG}. The reason why finite time breakdown can occur in the wave flow
analog of this problem is essentially due to the second order
nature of the equations. See Remark \ref{res_rem} below for more
thorough discussion.\\\

Lastly, we point out that our work is essentially independent of
previous techniques used for wave-maps, although we find it
extremely useful to keep in mind that the self-similar blow-up is
\emph{a priori} ruled out by Theorem \ref{C_th}. However, we do
refer to that result for the statement  of ``small energy implies
regularity'', which underlies much of what we do in the sequel. We
again stress that the fundamental structure we rely on in this paper
is the ``quasi-integrable'' and  ``super-symmetric" aspects of the
\emph{static} (elliptic) case of the equation \eqref{eq:EL}.
Specifically, the fact that such solutions may be constructed by
solving the first order Bogomol'nyi equations, as opposed to the
full second order Euler-Lagrange equations. These aspects enter
prominently into our analysis of the \emph{time-dependent}
problem.\\

In the remainder of this section we will give a detailed discussion
of the symmetry reduction we use in this work, as well as two
separate statements of our main theorem. \\

\noindent {\it Acknowledgments:} \,\,\,The authors would like to thank M.
Grillakis and J. Shatah for valuable discussions of the symmetric wave-map problem,
and J. Krieger and W. Schlag for valuable discussions regarding
their recent work \cite{KS_SW} and related stability problems of wave
equation solitons. We are also deeply indebted to Pierre
Raphael for pointing out to us a subtle miscalculation of the blowup rate in an earlier
version of the paper and for other valuable discussions. 
The authors would like to thank the MIT and UCSD
mathematics departments for their hospitality while this work was
being conducted.

\ret

\subsection{Symmetric reduction of the problem, and the
statement of the main theorem}

As we have already mentioned,  we will restrict our study of the
system \eqref{eq:EL} by enforcing some fairly rigid\footnote{It
would be extremely interesting to remove these in some way. For
example by either studying large initial deviations from an
equivariant soliton, or by studying small non-equivariant
perturbations of an equivariant soliton.} symmetry \emph{and}
``size'' assumptions. As seen in the introduction, the class of
solutions one has access to under these restrictions already
exhibits some interesting and striking phenomena, and in many ways
is still quite far from being understood. We now give an alternative
derivation of the symmetry assumption we use here.  This corresponds
to solutions behaving rigidly with respect to rotations on the base
manifold $\mathbb{R}^{2+1}$. That is, we require that along some
fixed time-line $(t,0)\in\mathbb{R}^{2+1}$ a rotation of $2\pi$
corresponds to a rotation of $2k\pi$ on the sphere $\mathbb{S}^2$
about some fixed axis. This type of symmetry dictates a more or less
canonical set of coordinates on the target, which is simply polar
coordinates centered about the axis of rotation. We write this in
the usual way in terms of two angles:
\begin{align}
        \{\Phi^1,\Phi^2\} \ &= \ (\phi,\theta) \ ,
        &g \ &= \ ds^2 \ = \ d\phi^2 + \sin^2(\phi) d\theta^2 \ ,
        \label{local_sphere}
\end{align}
where we restrict $0\leqslant \phi \leqslant \pi$ and $0\leqslant
\theta < 2\pi$, with $\phi=\pi,0$ the respective north and south
poles of the rotation axis.\\

With this choice of coordinates, our symmetry assumption boils down
to the simple relation $\theta \equiv k\Theta $ where $(t,r,\Theta)$
are polar coordinate on the base $\mathbb{R}^{2+1}$ and we mod with
respect to $2\pi$. In this case, the only remaining degree of
freedom is given by the quantity $\phi$ which can only depend on the
variables $(t,r)$. Because all of the Christoffel symbols
$\Gamma^\phi_{ij}$ vanish except in the case $i=j=\theta$, where we
have by a simple calculation $\Gamma_{\theta\theta}^\phi =
-\frac{1}{2} \frac{d}{d\phi}(\sin^2\phi)$, the general system
\eqref{basic_WM2} reduces to the single equation:
\begin{align}
        -\partial_t^2\phi + (\partial_r^2 +
        \frac{1}{r}\partial_r)\phi \ &=\ k^2\, \frac{\sin(2\phi)}{2r^2}
        \ , &k\ \ \in \ \mathbb{N}^+ \ , \label{EL_eqs_red}
\end{align}
where we implicitly enforce the boundary conditions $\phi(0)=0$ and
$\phi(\infty)=\pi$.\\

Before we continue, it is useful for us to record here the formula
for the Lagrangian density \eqref{WM_action} under this
k-equivariant symmetry reduction and in terms of the local
coordinates on the sphere $(\phi,\theta)$:
\begin{equation}
    \mathcal{L}[\phi] \ = \ \frac{1}{2}\ \left[
    -(\partial_t\phi)^2 + (\partial_r\phi)^2 + \frac{k^2}{2r^2}(1-\cos(2\phi))
    \right]\ . \label{red_action}
\end{equation}
In this notation the conserved energy \eqref{eq:Energy} becomes:
\begin{equation}
    E[\phi] \ = \ \pi\,
    \int_{\mathbb{R}^+}\ \left[
    (\partial_t\phi)^2 +
    (\partial_r\phi)^2 + \frac{k^2}{r^2}\sin^2(\phi)
    \right]\ rdr \ . \label{red_energy}
\end{equation}
The statement of our main theorem is now the following:\\

\begin{thm}[Finite time energy concentration for wave-maps]\label{main_thm}
Consider the full wave-map equation \eqref{basic_WM2} with $\mathbb{S}^2$ target
under the
equivariant restriction to equation \eqref{EL_eqs_red}. Then for any
integer $4\leqslant k$, and for any sufficiently small constant $0 < c_0\ll 1$
with the property that for any $\epsilon\leqslant c_0^2$,
we can find a set of smooth (in the sense of the full
map on $\mathbb{R}^{2+1}$) Cauchy data:
\begin{align}
    \phi(0) \ &= \ \phi_0^\epsilon \ ,
    &\partial_t \phi\, (0) \ &= \ \dot{\phi}_0^\epsilon \ , \notag
\end{align}
with energy size $E[\phi^\epsilon] = 4\pi k + \epsilon^2$ such that
this solution collapses at a finite time $T^{**}$. More
specifically, this solution collapses at a rate bounded by a ``$\log$-modified
self-similar'' dynamic in the sense that there exists a universal
time independent profile $\underline{\phi}^k$, and a real parameter
$0 < \lambda(t)$, such that:
\begin{equation}
    E\Big[\phi(t,r) - \underline{\phi}^k(\lambda r)\Big]
    \ \lesssim \ \epsilon^2
\end{equation}
and such that for any $0<M$ and times sufficiently close to $T^{**}$ one has the bound:
\begin{equation}
   \frac{M}{(T^{**}-t)} \ \leqslant \  \lambda(t) \ \leqslant \
   c_0^\frac{1}{4}
    \frac{\sqrt{|\ln(T^{**}-t)|}}{(T^{**}-t)}
    \ . \label{b_rates}
\end{equation}
Finally, this type of blowup is  stable within the class of initial data
in the sense that there exists a weighted Sobolev space $H^{s,m}$
(see \eqref{Hspace_def} for a definition), such that the
$c_0\epsilon$ ball about
$(\phi_0^\epsilon,\dot{\phi}_0^\epsilon)$ in $H^{s,m}$ also leads to
collapse with the same universal profile $\underline{\phi}^k$ and
the same bound \eqref{b_rates}.
\end{thm}\ret

\ret

\subsection{The family of static solutions and a modulational
version of Theorem \ref{main_thm}}\label{mod_state_sub}

As we have mentioned previously, it is well known from work of
Struwe (see again \cite{S_WM1}) that any blowup of the form
described in Theorem \ref{main_thm} must in fact be a ``bubbling
off'' of a static solution to the equation \eqref{EL_eqs_red}. That
is, after rescaling the solution $\phi(t)$ as described in Theorem
\ref{main_thm}, the resulting profile should be a solution to
\eqref{EL_eqs_red}. In the sequel, we will actually take the
converse approach and give an explicit construction of such bubbling
off solutions. This will be done in a way which is generally
consistent with the decomposition \eqref{intro_decomp} of
the introduction.
Our method also naturally shows that this process  is
reached from a generic (in the symmetric sense) set of initial data,
and that it enjoys a certain universality which is embodied by the
blowup rate \eqref{b_rates}.\\

To get things started, we derive the formula for the energy
minimizer of the (full) action \eqref{red_action}. This is
just a recalculation of lines \eqref{general_lines} in the current
notation. Completing the square in the spatial terms in the
energy \eqref{red_energy} we can write it as:
\begin{align}
    E[\phi] \ &= \ \pi\,
    \int_{\mathbb{R}^+}\ \left[
    (\partial_t\phi)^2 + \big(\partial_r\phi-\frac{k}{r}\sin(\phi) \big)^2
    \right]\ rdr \ + \ 2\pi\,
    \int_0^\infty\ k\sin(\phi)\partial_r\phi \ dr
    , \notag \\
    &= \ \pi\,
    \int_{\mathbb{R}^+}\ \left[
    (\partial_t\phi)^2 + \big(\partial_r\phi-\frac{k}{r}\sin(\phi) \big)^2
    \right]\ rdr \ + \ 4k\pi \ . \label{difference_energy}
\end{align}
Thus, one has the universal lower bound $4k\pi \ \leqslant
E[\phi]$, which can be reached if we can find a function $I^k$ with
the property that $I^k(0)=0$ and $I^k(\infty)=\pi$, and which
satisfies the following equations:
\begin{align}
    \partial_t I^k \ &= \ 0 \ , \notag\\
    r\partial_r I^k \ &= \ k\sin(I^k) \ . \label{B_eq}
\end{align}
We shall refer to the solution $I^k$ as the
harmonic map {\it soliton}.
A direct calculation reveals that
the function $I^k$ is given by the explicit formula:
\begin{equation}
    I^k(r) \ = \ 2\tan^{-1}(r^k) \ . \label{instanton}
\end{equation}
We also denote:
\begin{align}
    I(r)\ &:=\ I^k(r) \ ,
    &J(r)\ &:=\ r\partial_r I(r) \ . \label{JI_defs}
\end{align}
Note that since the equations \eqref{B_eq} are homogeneous, the general solution
$I_\lambda$ is only defined up to a rescaling, as we have already mentioned on line
\eqref{sc_trans} above. Now, in terms of these objects we can state the
following more technical and precise version  Theorem
\ref{main_thm}, which is what we shall actually prove in the
sequel:\\

\begin{thm}[Modulational version of the main theorem]\label{mod_thm}
Consider the reduced wave-map equation \eqref{EL_eqs_red} with $4\leqslant k$.
Suppose we are given a pair of sufficiently small positive constants
$\epsilon,c_0$ with $\epsilon\leqslant c_0^2$, and an initial data set of the
form\footnote{Note that our choice of initial data already requires $2\leqslant k$ as for
$k=1$ we have that $\|J\|_{L^2(rdr)}=\infty$.}:
\begin{align}
    \phi(0) \ &= \ I + u_0 \ ,
    &\partial_t\phi\, (0) \ &= \ \frac{\epsilon}{\pi}\, \lp{J}{L^2(rdr)}^{-2}\cdot J +
    g_0 \ , \label{special_C_data}
\end{align}
where:
\begin{equation}
    \int_{\mathbb{R}^+}\ u_0\cdot J\ rdr \ = \ 0
    \ , \notag
\end{equation}
and obeys the smallness condition:
\begin{equation}
    \lp{(u_0,g_0)}{H^{2,1}}^2 \ \leqslant \ c_0^2 \epsilon^2 \ ,
    \label{special_smallness}
\end{equation}
where we have set:
\begin{equation}
    \lp{(u_0,g_0)}{H^{2,1}}^2 \ = \
    \sum_{i=0}^1 \
    \int_{\mathbb{R}^+}\ (1+r^2)^{1-i}
    \big[ (\partial^i_r g_0)^2 +
    \frac{(g_0)^2}{r^2}
    + (\partial_r^{i+1} u_0)^2
    + \frac{(\partial_r^i u_0)^2}{r^2} \big]\ rdr \ . \label{Hspace_def}
\end{equation}
Then we have that the following is true: There exists a continuous
time dependent parameter $\lambda(t)$ with $\lambda(0)=1$, and such
that the solution $\phi$ to \eqref{EL_eqs_red} with initial data
\eqref{special_C_data} splits into the sum:
\begin{equation}
     \phi(t,r) \ = \ I(\lambda(t)r) + u(t,r) \ , \label{first_I_decomp}
\end{equation}
where the ``remainder'' term $u$ satisfies the bounds:
\begin{equation}
    \int_{\mathbb{R}^+}\ \big[ (\partial_t u)^2 + (\partial_r u)^2
    + \frac{u^2}{r^2} \big]\ rdr \ \lesssim \  \epsilon^2 \ ,
     \label{orbit_smallness}
\end{equation}
for all times the solution exists. Furthermore, there exists a
finite time $T^{**}$ such that $\lim_{t\to T^{**}}
\lambda(t)\to\infty$. Finally, this parameter obeys the
following bounds for times $t$ sufficiently close to
$T^{**}$:
\begin{equation}
    \frac{M}{(T^{**}-t)} \ \leqslant \  \lambda(t) \ \leqslant \
   c_0^\frac{1}{4}
    \frac{\sqrt{|\ln(T^{**}-t)|}}{(T^{**}-t)}
    \ . \label{blowup_rate}
\end{equation}
\end{thm}\ret

\begin{rem}\label{lower_bound_rem}
The lower bound in the blowup rate \eqref{blowup_rate} follows easily from
the orbital stability bound \eqref{orbit_smallness} and Theorem \ref{C_th}.
Therefore, in the sequel we shall concentrate on establishing blowup
with the \emph{upper} bound on line \eqref{blowup_rate}. The reader should
note however that the presence of the extra small constant $c_0^\frac{1}{4}$,
which may go to zero with $\epsilon$, indicates that the true blowup rate
is even closer to self similar than the $\sqrt{|\ln(T^{**}-t)|}$ correction.
We will return to this delicate issue in a later work.
\end{rem}\ret

\begin{rem}
The extra decay provided by the bounds
\eqref{special_smallness} is not essential to what we do here and the result also holds
in the space $H^{2,0}$.
It is assumed here as a convenience that will simplify the exposition.
However, the extra regularity afforded to us in the norm
\eqref{Hspace_def} will be used in a crucial way. We also remark
that the norm \eqref{Hspace_def} is consistent with smoothness of
the derivative of the initial data when considered as  map from
$\mathbb{R}^2$ into the pullback bundle $\cup_{x\in\mathbb{R}^{2+1}}
\Phi_x^*( T\mathbb{S}^2)$. This is a consequence of some simple
calculations involving the frame Christoffel symbols
\eqref{frame_chris}.
\end{rem}\ret

\begin{rem}\label{res_rem}
As we have already mentioned, the blowup mechanism we exhibit here is ignited
by a spectral phenomenon. The choice of
initial data \eqref{special_C_data} guarantees
that the time derivative of the wave-map $\dot{\phi}_0$ has a ``strong"
projection onto the ``ground state'' $J(r)$ of the equation
\eqref{EL_eqs_red} linearized around the soliton
$I(r)$. This is precisely the coefficient in a Riccati equation
for the scaling parameter $\lambda(t)$ (see formula \eqref{basic_ODE} below).
The Riccati equation generates the first self-similar epoch of collapse which lasts
on the time interval of size $\sim \epsilon^{-1}$ and, as the projection
of the time derivative $\pa_t\phi$ on the ``ground state" $J_\lambda$ of the
modulated soliton $I_\lambda$ goes to zero, is eventually replaced by a
more violent accelerated regime leading to the blow-up.
For this initial phenomenon to take place it is crucial
that the linearized  ground state $J$ is an $L^2(rdr)$ function and
that the projection of the time derivative of $\phi$  on the ground state $J$
is initially non-trivial.
That is, one of the main
things which makes our analysis possible is that
the first order\footnote{That is, in this notation the second order
wave equation \eqref{EL_eqs_red} can be written as a first order
system.} field quantities $(\phi-I_\lambda,\partial_t\phi)$ can not be \emph{both}
orthogonal to the eigenfunction of  \eqref{EL_eqs_red} linearized
around $I(\lambda r)$ (unless one restricts
the initial data to a co-dimension one submanifold).
\\

In this regard there are some interesting open questions connected
with the value of the homotopy index.
For $k=2$ the linearized ground state is still in $L^2$, so it is
likely that an adaptation of our methods is possible. This is
important because it is this case which is most closely related to
Yang-Mills (see the next remark). For the unit homotopy class,
$k=1$, the situation appears to be more complicated. In this case
the linearized ground state just misses $L^2$ by a $\log$. The major open
problem here seems to be whether there is complete instability of
the kind stated in Theorem \ref{main_thm}, or if small enough
perturbations (in some space) are asymptotically stable, with blowup
occurring as some kind of ``critical phenomena'' depending
delicately on the size of the perturbation. Another interesting
thing is that in the case of $k=1$, there are some numerical
simulations which seem to indicate that the blowup, while taking
place, occurs at an algebraically \emph{different} rate from
\eqref{blowup_rate} (see again \cite{B_WM}). On the other hand,
there are further numerical and heuristic results (see \cite{LS})
which suggest the validity of the $\log$-modified behavior even in
this case $(k=1)$. We believe that both the $k=1,2$ cases of (a
possible analog of) Theorem \ref{main_thm} deserve further serious
investigation in terms of numerics, heuristics, and theory.
\end{rem}\ret

\begin{rem}
Another important issue we call the readers attention to
is that in the case of the critical Yang-Mills, heuristic
arguments as well as numerical evidence point to blowup with a
modified self-similar asymptotic of the same form as
\eqref{blowup_rate} (see again \cite{BS_YM}).
In fact, the spherically symmetric reduction of the $(4+1)$-dimensional Yang-Mills
equations is very closely connected with the homotopy $k=2$ case for the $O(3)$
$\sigma$-model.  This strongly suggests that the methods we
develop here will transfer to the case of the Yang-Mills model as well, and this
will be the subject of a forthcoming work of the authors.
\end{rem}\ret


\subsection{A few more calculations}

Before continuing on, we list here  some simple formulas
involving the unit solitons $I^k$ which will be of particular
importance to us in the sequel:
\begin{align}
    r\partial_r I^k \ = \ k\sin(I^k) \ &= \ k\frac{2r^k}{1+r^{2k}} \ ,
    &\cos(I^k) \ &= \ \frac{1-r^{2k}}{1+r^{2k}} \ , \label{I_form1}\\
    \sin(2I^k) \ &= \ 4\frac{r^k - r^{3k}}{(1+r^{2k})^2} \ ,
    &\cos(2I^k) \ &= \ \frac{1-6r^{2k} + r^{4k}}{(1+r^{2k})^2}
    \ . \label{I_form2}
\end{align}
Also, in the sequel we will refer to any specific instance of $I^k$
as simply $I$, and we remind the reader that we are assuming
$4\leqslant k$.\\

\ret

\subsection{Vanishing of the wave-map}

We end this section by recording and proving a simple geometric
lemma which will be of central importance to us throughout the
sequel. We will show that a
wave-map $\Phi$ together with its derivative vanish
at the origin $r=0$ when computed in the pair of
local coordinates \eqref{local_sphere} and $(r,\Theta)$.\\

\begin{lem}[Admissibility condition for the wave-map $\Phi$]\label{admiss_lem}
Let $\Phi$ be a smooth k-equivariant function from the plane
$\mathbb{R}^2$ into the sphere $\mathbb{S}^2$. Then if $2\leqslant
k$ one has that:
\begin{align}
        |\partial_r \phi| \ &\leqslant \ C_\phi\, r \ ,
        &0\ \leqslant \ r \ &\leqslant \ 1 \ . \label{admiss_bound}
\end{align}
\end{lem}\ret

\begin{rem}
Note that the condition $2\leqslant k$ for estimate
\eqref{admiss_bound} is crucial, as the formula \eqref{instanton}
shows for unit homotopy class instanton $I^1$.
\end{rem}

\begin{proof}[Proof of the estimate \eqref{admiss_bound}]
Our first step is to establish that $\partial_r\phi$ is
continuous and vanishes at $r=0$.\\

First of all, notice that along any fixed radial line
$\Theta=const$, the vector-field $\partial_r$ is a
continuous section of $T\mathbb{R}^2$. The same is true of the
field $\frac{1}{r}\partial_\Theta$. Furthermore, one has that:
\begin{equation}
        \lim_{\substack{\Theta = 0 \\ r\to 0}}\ \partial_r \ = \
        -\, \lim_{\substack{\Theta = \frac{\pi}{2} \\ r\to 0}}\
        \frac{1}{r}\partial_\Theta \ . \notag
\end{equation}
Therefore, by continuity we must have that:
\begin{equation}
        \lim_{\substack{\Theta = 0 \\ r\to 0}}\ \lp{\partial
        \Phi(\partial_r)}{}^2 \ = \
        \lim_{\substack{\Theta = \frac{\pi}{2} \\ r\to 0}}\
        \lp{\partial\Phi(\frac{1}{r}\partial_\Theta)}{}^2 \ , \notag
\end{equation}
where $\lp{\cdot}{}^2$ refers to the metric \eqref{local_sphere}.
Computing both sides of this last equation, we see that not only is
$\partial_r\phi$ continuous (and hence bounded) on the interval
$[0,1]$, but that we also have:
\begin{equation}
         \lim_{ r\to 0}\ |\partial_r\phi| \ = \
        \lim_{ r\to 0}\
        \frac{k|\sin(\phi)|}{r} \ . \label{r0_iden}
\end{equation}
Using now the fact that $\phi(0)=0$ to write $\phi(r)=\int_0^r
\partial_y\phi(y)\, dy$, upon substitution of this integral into the right hand
side of \eqref{r0_iden} , we see from the fundamental theorem of
calculus and the condition $2\leqslant k$, that we must in fact have
$\partial_r\phi(0)=0$.\\

It remains to show that $\partial_r\phi$ vanishes uniformly (with
non-uniform constant) in $r$. To do this, we compute the Dirichlet
energy:
\begin{equation}
        e(\phi) \ = \ g_{ij}
    \Big[\partial_r\Phi^i \partial_r\Phi^j + \frac{1}{r^2 }\partial_\Theta \Phi^i
    \partial_\Theta \Phi^j\Big]
     \ = \
        |\partial_r\phi|^2 + \frac{k^2}{r^2}|\sin(\phi)|^2 \ . \notag
\end{equation}
This is a $C^\infty$ function on $\mathbb{R}^2$ which depends on the
radial variable only. Furthermore, we have that $e(\phi)(0)=0$. Thus
it is necessary that:
\begin{align}
        |e(\phi)| \ &\leqslant \ C \, r^2 \ ,
        &0\ \leqslant \ r \ &\leqslant \ 1\ . \notag
\end{align}
for some constant that depends on $\phi$. In particular, we have the
bound \eqref{admiss_bound}.
\end{proof}

\ret\ret

\section{Some notational conventions and an overview}

In this section, we will first list some standard notational
conventions that will be useful throughout the sequel. We then give
a quick technical overview of the main result.

\ret

\subsection{Some notation}\label{not_sect}

Throughout this paper, we shall employ the standard notation
$A\lesssim B$ to mean $A\leqslant CB$ for two quantities $A$ and
$B$, where $C$ is a fixed constant. There is no uniformity in this
notation for \emph{separate} instances. That is, separate
occurrences of $\lesssim$ on the same page will not necessarily imply
$C$ is the same for each. Another, less standard, notation which
will be of great use is the following:\\

\begin{notation}\label{abs_F_not}
For any pair of non-negative integers  $0\leqslant m,n$ we will
denote by $F^{m,n}$ any $C^\infty(\mathbb {R}^+)$ function which satisfies the
following bounds:
\begin{equation}
    \big|(r\partial_r)^i F^{m,n}\big| \ \lesssim \ C_i\,
    \frac{r^m}{(1+r)^{m+n}} \ . \notag
\end{equation}
We will also use a shorthand notation for the case $m=0$. Here we shall set
$F^n=F^{0,n}$, so that we have:
\begin{equation}
    \big|(r\partial_r)^i F^n\big| \ \lesssim \
    \frac{C_i}{(1+r)^n} \ . \notag
\end{equation}
We also denote the $\lambda$ rescaling of these functions by
$F^{m,n}_\lambda(r) = F^{m,n}(\lambda r)$, and similarly for
$F^n_\lambda$. Finally, let us remark that different instances of
$F^{m,n}_\lambda$ on any line, or between lines, can mean separate
functions.
\end{notation}\ret

This notation will occur so frequently in the sequel that is will be
useful for us to record here several instances which involve either
time or space derivatives, or multiplication by powers of $r$.
Collectively these are the following, where we assume
$j\in\mathbb{Z}$ is such that $-m\leqslant j$ in the first identity
and $1\leqslant m$ in the third:
\begin{align}
    r^j F^{m,n}_\lambda \ &= \ \lambda^{-j}
    F^{m+j,n-j}_\lambda \ ,
    &\partial_t F^{m,n}_\lambda \ &= \
    \dlambda\lambda^{-1} F^{m,n}_\lambda \ ,
    &\partial_r F^{m,n}_\lambda \ &= \
    \lambda F^{m-1,n+1}_\lambda \ . \label{F_rules}
\end{align}
All of these are immediate from Definition \ref{abs_F_not} above.\\

Also, in the sequel we will in general use the $\lambda$-subscript
notation to denote the $\lambda$ rescaling of a given function. For
example $I_\lambda(r) = I(\lambda r)$.

\ret

\subsection{An overview}\label{o_sect}

As
is perhaps already clear at this point, our method for establishing
Theorem \ref{mod_thm} is to control directly a certain modulational
equation for the time dependent scaling parameter $\lambda(t)$, and to show
that this evolves according to a blow-up ODE. Thus, in this sense our
work is  closely related in spirit to the modulational stability approach
originally pioneered by M. Weinstein (see \cite{W_mod}) and later
sharpened by Buslaev and Perelman (see \cite{BP_NLS}) to study
solitons dynamics of the focussing non-linear Schr\"odinger equation.
The major difference however is
that we are actually trying to show that there is an extremely
strong asymptotic \emph{instability}. This of course introduces a
serious problem when one tries to control the non-linear equation
\eqref{EL_eqs_red} linearized around the modulated soliton
$I_\lambda$. Nonetheless we begin by using the decomposition $\phi = I_\lambda + u$
from line \eqref{first_I_decomp},
and then linearizing \eqref{EL_eqs_red} around $I_\lambda$:
\begin{equation}
    \partial_t^2 u + H_\lambda\, u \ = \ -\ddot{I}_\lambda \ + \
    \mathcal{N}(u) \ , \label{eq:linear}
\end{equation}
where the Hamiltonian is given by $H_\lambda=-\pa_r^2 - r^{-1} \pa_r +Q_\lambda(r)$, and
the nonlinear term $\mathcal {N}(u)$ containing quadratic and higher order
terms in $u$ (also containing a factor of $r^{-2}$). \\

Our first task, dealt with in Section \ref{orbit_sect},
is to prove orbital stability of the
modulated soliton $I_\lambda$ under the condition that the ``radiation" part of the
solution (i.e. $u$) is orthogonal to the function $J_\lambda=r\pa_r I_\lambda$, which is
the unique eigenfunction of the Hamiltonian $H_\lambda$. The latter
is a consequence of the fact that $I_\lambda$ realizes an absolute minimum of the
energy \eqref{red_energy} associated with the full nonlinear problem \eqref{EL_eqs_red}.
The orthogonality condition  provides us with an ODE for the
scaling parameter $\lambda(t)$, which is coupled to the radiation term $u$:
\begin{equation}
    \dlambda \Big(\langle J,J \rangle -
        \langle u(\lambda^{-1}r) , r\partial_r J\rangle
        \Big)  \ = \  \Big (\lambda\langle\phi_t ,
        J_{\lambda}\rangle\Big )\lambda^2 \ , \label{exp_ODE}
\end{equation}
while the orbital stability statement will give us a very
weak control of the remainder $u$:
\begin{equation}
    \int_{\mathbb{R}^+}\ \left[ (\partial_t\phi)^2
        + (\partial_r u)^2 + \frac{k^2}{r^2}u^2\right]\
        rdr \ \lesssim \ \epsilon^2 \ . \label{exp_orbit}
\end{equation}
Due to the coupling between the scaling parameter $\lambda(t)$ and
the radiation $u$, to control $\lambda(t)$ to the extent that we can
show $\lambda(t)\to\infty$ in finite time requires much better
control on the radiation term $u$. The usual procedure for dealing
with this is to scale out the modulational parameter $\lambda(t)$ at
each fixed time, so \emph{time independent} spectral methods can be
used to control the linearized equation. This procedure works well
if one can prove that there is a \emph{slow} limit of the parameter
$\lambda(t)$, but it obviously causes a catastrophe if $\lambda(t)$
grows rapidly. In
this case, a truly non-linear approach is needed.\\

The reason why standard  non-linear estimates, for example the kind used to
prove orbital (Lyapunov) stability (i.e. \eqref{exp_orbit} above
), are not sufficient to reach the blow-up time
$\sim \epsilon^{-1}$  is basically due to the fact that their
application to the ODE \eqref{exp_ODE} is not truly scale invariant.
That is, the use of \emph{fixed time} estimates which result from
orbital stability analysis causes a loss  (of scaling) when one
integrates over time. Such integrations seem unavoidable when
analyzing \eqref{exp_ODE}. To overcome this problem requires
uncovering a non-linear dispersion phenomenon\footnote{A non-linear
dispersion phenomenon (of a different nature) has been observed and
used by Merle and Raphael in their work \cite{MR} on the blow-up
analysis for the critical focussing non-linear Schr\"odinger
equation.} in the equation \eqref{eq:linear} for the radiation term
$u$. That such a dispersive process indeed takes place is in some
sense the miracle of the equation \eqref{eq:linear}. More
specifically, as the soliton $I_\lambda$ collapses it actually
repels the excess radiation away from the origin. This evacuation
process only causes the soliton to collapse at a faster rate, and it
is what is ultimately responsible for the acceleration of self-similar
behavior governed by the LHS of the
blow-up rate \eqref{blowup_rate}.\\

One of the most interesting issues in this
paper is the mechanism by which this ``repulsive'' behavior of the linearized
equations manifests itself mathematically. This is where the
``quasi-integrable" system aspect of the static version of \eqref{EL_eqs_red}
comes in. As we have already discussed in the introduction, static
solutions to \eqref{EL_eqs_red} are generated by the {\it first order}
Bogomol'nyi equation \eqref{B_eq}. When one linearizes \eqref{EL_eqs_red} around
these static solutions, the corresponding Hamiltonian $H_\lambda$ splits as
a product of two first order operators which are adjoints of each
other. That is:
\begin{equation}
    H_\lambda\  = \ A^*_\lambda A_\lambda \ , \notag
\end{equation}
where $A_\lambda$ is the linearization
of \eqref{B_eq}. The reason why this splitting is so useful is that the
Hamiltonian $H_\lambda$ also possesses its super-symmetric companion:
\begin{equation}
    \tilde H_\lambda\ = \ A_\lambda A^*_\lambda \ , \notag
\end{equation}
with $A^*_\lambda$ and $A_\lambda$ being the analogs of the creation and annihilation
operators, and $H_\lambda$ and $\tilde H_\lambda$ related to each other according to
the remarkable intertwining relation:
\begin{equation}
    A_\lambda H_\lambda \ = \ \tilde H_\lambda A_\lambda \ . \label{eq:intertwin}
\end{equation}
Such a splitting elucidates the non-negativity of $H_\lambda$, and identifies the function
$J_\lambda$, which is the kernel of $A_\lambda$, as the ground (vacuum) state of $H_\lambda$.
In addition, the intertwining property \eqref{eq:intertwin} allows us to simply conjugate the
problem \eqref{eq:linear} to one whose linear part involves the more manageable
Hamiltonian $\tilde H_\lambda$:
\begin{equation}
    \pa_t^2 (A_\lambda u) + \tilde H_\lambda (A_\lambda u) \ = \
    - A_\lambda (\ddot{I}_\lambda) +
    A_\lambda \mathcal {N}(u) + [\pa_t^2,A_\lambda]  u \ . \label{eq:conj}
\end{equation}
The Hamiltonian $\tilde H_\lambda$, which is obtained from
$H_\lambda$ by the process of ``removing" its ground state, is of
the explicit form $\tilde H_\lambda=-\pa_r^2 -r^{-1} \pa_r +
V_\lambda(r)$ and involves a \emph{space-time repulsive}
time-dependent potential. This means that  for the problem
\eqref{eq:conj} one may proceed via purely physical space methods,
and there is no difficulty in handling extremely violent growth of the
scaling parameter $\lambda(t)$. What we can do is to establish quite
strong (i.e. scale invariant) integrated and fixed time energy
estimates  (i.e. so called Morawetz type estimates), while keeping
precise track of the influence of the source  terms on the right
hand side of \eqref{eq:conj} involving the scaling parameter
$\lambda$. To undo the conjugation procedure embodied in
\eqref{eq:conj}, we only need to use the fact that $u$ is orthogonal
to the kernel of $A_\lambda$, because through a little elementary
functional analysis this allows one to turn our Morawetz estimates
into ones involving only the term $u$ (as opposed to $A_\lambda u$).
Once these estimates
are established it is possible to show, through a somewhat lengthy
calculation, that after a long self-similar epoch where
$C_0\dlambda\sim \epsilon_0 \lambda^2$,  the modulation ODE for
$\lambda(t)$  enters another monotonic  regime where it takes the
final form:
\begin{equation}
    C_0\dot\lambda(t) \ \sim  \
    \epsilon_0 \lambda^2(t) - \lambda^2(t)
    \int_0^t O\Big( \frac{\dlambda^4}{\lambda^7}(s)\Big)ds \ . \notag
\end{equation}
It is this  ODE which leads to the blow-up,
and also gives the tight \emph{upper} bounds in \eqref{blowup_rate}. The reader should
compare this last formula to the blowup ODE for modulated Yang-Mills
instantons derived in \cite{BS_YM} through heuristic arguments.\\

\noindent We now turn to the details of all of this. As is common in this
type of work, many of our assumptions will be bootstrapped. We shall follow the
outline:\\

\begin{itemize}
    \item In Section \ref{spectral_sect}
    we discuss the Hamiltonians $H_\lambda$ and $\tilde H_\lambda$.\\
    \item In Section \ref{orbit_sect} we derive the basic
    ODE for the scaling parameter $\lambda(t)$ and
    prove orbital stability statement.\\
    \item In Section \ref{ODE_sect}, assuming certain estimates on
    the ``radiation" part of the solution, we obtain a much
    refined \emph{closed form} (i.e.
    without explicit dependence on the $u$) of the modulation ODE, establish its
    monotonic and algebraic properties,
    and prove blow-up along with an explicit rate bound. \\
    \item In Section \ref{mor_sect}, assuming the monotonic properties
    of the scaling parameter $\lambda(t)$,
    we prove the integrated space-time and fixed time bounds for the radiation by
    making use of the conjugated Hamiltonian $\tilde H_\lambda$.\\
    \item In Appendix \ref{C_comp_app} give some further explicit computations
    needed in the analysis of the blowup ODE.\\
    \item In Appendix \ref{coercive_app}
    we establish some general coercive properties for the class of first order
    operators related to $A_\lambda$.
\end{itemize}

\ret\ret

\section{The Linearized Equations and a Basic Spectral Calculation}\label{spectral_sect}

Our purpose here is to derive and record certain calculations
involving the linearization of the equation \eqref{EL_eqs_red}
around a time dependent modulation of the soliton $I(\lambda(t)r)$.
That is, we decompose the full solution as:
\begin{equation}
    \phi(t,r) \ = \ I(\lambda(t)r) + u(t,r) \ . \label{I_u_decomp}
\end{equation}
This yields the following set of formulas for equation \eqref{EL_eqs_red}
linearized around $I_\lambda$:
\begin{equation}
    \partial_t^2 u + H_\lambda\, u \ = \ -\ddot{I_\lambda} \ +\ \mathcal{N}(u) \ , \label{lin_eq}\\
\end{equation}
where we have set:
\begin{equation}
    \mathcal {N}(u) \ = \ \frac{k^2\sin(2I_\lambda)}{2r^2}
    \cdot(1-\cos(2u)) \ +\  \frac{k^2\cos(2I_\lambda)}{r^2}
    \cdot(u-\frac{1}{2}\sin(2u))
    \ . \notag
\end{equation}
Here $H_\lambda$ is the linearized Hamiltonian:
\begin{align}
    H_\lambda \ &= \ -\partial_r^2 -\frac{1}{r}\partial_r
    +\frac{k^2}{r^2}\cos(2I_\lambda) \ , \notag\\
    &= \ A_\lambda^*A_\lambda \ , \label{the_fact}
\end{align}
where the first order operators $A,A^*$ are given by:
\begin{align}
    A_\lambda \ &= \ -\partial_r + \frac{k}{r}\cos(I_\lambda) \ ,
    &A^*_\lambda \ &= \ \partial_r + \frac{1}{r} +
    \frac{k}{r}\cos(I_\lambda)
    \ . \label{A_ops}
\end{align}
The spectrum of $H_\lambda$, defined as a self-adjoint operator on $L^2(rdr)$ and obtained by
taking the closure of $C_0^\infty(\mathbb{R}^+)$ in the graph norm of $H_\lambda$,
is easily computed via the factorization
\eqref{the_fact}, or via the knowledge that the ground state
$I_\lambda$ of the (static form of the) equation \eqref{EL_eqs_red}
is unique. There is a unique eigenfunction, which has zero energy,
and it is given by the formula:
\begin{equation}
    J_\lambda \ = \ r\partial_r I_\lambda \ = \
    k\sin(I_\lambda) \ . \label{J_def}
\end{equation}
In particular, $J_\lambda$ solves the first order ``linearized
Bogomol'nyi equation'':
\begin{equation}
    A_\lambda J_\lambda \ = \ 0 \ . \label{lin_bog}
\end{equation}\ret
The  absolutely continuous spectrum fills the half-line $[0,\infty)$.

Of primary importance for use here will be the conjugate operator formed
by $A_\lambda A^*_\lambda$, the super-symmetric companion $\td{H}_\lambda$
of $H_\lambda$.
\begin{align}
    \td{H}_\lambda \ = \ A_\lambda A^*_\lambda \ &= \
    -\partial_r^2 -\frac{1}{r}\partial_r +
    \frac{k^2+1}{r^2} + \frac{2k}{r^2}
    \cos(I_\lambda) \ , \label{conj_ham} \\
    &= \ -\partial_r^2 -\frac{1}{r}\partial_r +
    V_\lambda(r) \ . \notag
\end{align}
Recall that the Hamiltonians $H_\lambda$ and $\td{H}_\lambda$ are related via an intertwining
relation \eqref{eq:intertwin}. As opposed to $H_\lambda$,
the spectrum of $\td{H}_\lambda$ has only an absolutely continuous component
filling $[0,\infty)$. What is more important is that while in our application
the Hamiltonian $\td{H}_\lambda$ will be \emph{time dependent}, it has
a remarkable structure which allows one to prove strong local energy
decay estimates even if the parameter $\lambda(t)$ grows in an unconstrained
fashion. The needed properties follow from entirely elementary
calculations and are as follows:
\begin{align}
    V_\lambda \ &\geqslant \ \frac{(k-1)^2}{r^2} \ ,
    &\hbox{(Positive)}\ , \label{positivity}\\
    -\, \partial_r V_\lambda \ &= \ \frac{2(k^2+1)}{r^3}
    + \frac{4k}{r^3}\cos(I_\lambda) + \frac{2k^2}{r^3}\sin^2(I_\lambda)
    \ , \notag\\
    &\geqslant \ \frac{2(k-1)^2}{r^3} \ ,
    &\hbox{(Space-Repulsive)}\ , \label{repusivity}\\
    -\, \partial_t V_\lambda \ &= \ \frac{\dlambda}{\lambda}
    \cdot \frac{2k^2}{r^2}
    \sin^2(I_\lambda) \ , &\hbox{(Time-Repulsive)}\ . \label{pos_energy}
\end{align}\ret

We conclude this section by refining the decomposition \eqref{I_u_decomp}.
This will be extremely important for us in
the sequel, and it will ultimately lead us to the precise asymptotic
\eqref{blowup_rate}. What we will need to do is further decompose the radiation
term as:
\begin{equation}
    u \ = \ w_0 + w \ , \notag
\end{equation}
where the ``leading term'' $w_0$ is meant solely to eliminate the
term $\ddot{I}_\lambda$ on the right hand side of \eqref{lin_eq},
while at the same retaining the orthogonality relation:
\begin{equation}
    \langle w_0 , J_\lambda \rangle \ = \ 0 \ . \label{w0_orth}
\end{equation}
The equation we use to generate $w_0$ is (the reason for this will
become apparent in the sequel):
\begin{equation}
 A_\lambda (H_\lambda w_0)= \td{H}_\lambda ( A_\lambda w_0) \ = \ -A_\lambda(\ddot{I}_\lambda)
    \ . \label{w_condition}
\end{equation}
To further the computation, we use the fact that $A_\lambda J_\lambda=0$
to write:
\begin{align*}
    A_\lambda(\ddot{I}_\lambda) \ = \ \frac{\dlambda}{\lambda}
    A_\lambda  \partial_t
    (r\partial_r I_\lambda)
    = \ \frac{\dlambda^2}{\lambda^2} A_\lambda \, r\partial_r
    (J_\lambda)
    = \ \frac{\dlambda^2}{\lambda^2} A_\lambda
    ( J + r\partial_r J)_\lambda
    \ .
\end{align*}
Therefore, by peeling off the common factor of $A_\lambda$ from both
sides of \eqref{w_condition} it suffices to solve:
\begin{equation}
    H_\lambda w_0 \ = \ -\frac{\dlambda^2}{\lambda^2}
    (J + r\partial_r J)_\lambda
    \ . \notag
\end{equation}
We now use the ansatz $w_0 = \dlambda^2\lambda^{-4}\big[ K(\lambda
r)  - \gamma J_\lambda\big]$, where $\gamma$ is a normalization constant
which will be chosen as to satisfy \eqref{w0_orth}.
Eliminating common factors,
and rescaling the spatial variable we find that:
\begin{equation}
    H_1 (K) \ = \ -(J + r\partial_r J) \ , \notag
\end{equation}
from which a direct computation shows that $K(r) = \frac{r^2}{4}
J(r)$ is the desired solution. Therefore we have that:
\begin{equation}
    w_0(t,r) \ = \ \frac{\dlambda^2}{\lambda^4}(t)\cdot\big(aJ_\lambda(r)
    + b (r^2J)_\lambda\big) \ , \label{w0_def}
\end{equation}
where the coefficients are given by the explicit formulas:
\begin{align}
    a \ &= \ -\, \frac{1}{4} \langle J, (r^2 J) \rangle
    \cdot \lp{J}{L^2(rdr)}^{-2}
    \ ,
    &b \ &= \ \frac{1}{4} \ . \label{ab_defs}
\end{align}\ret

Before ending this section, let us translate the abstract function notation
on line \eqref{F_rules} into some specific bounds which will be used
many times in the sequel. These are a consequence of simple explicit
formulas, and the restriction $4\leqslant k$:
\begin{subequations}\label{schematic_not}
\begin{align}
    w_0\ &= \ \dlambda^2\lambda^{-4}( F^{4,4}_\lambda + F^{6,2}_\lambda)\ ,
    &\partial_t(w_0) \ &= \ \ddot{\lambda}\dlambda\lambda^{-4}( F^{4,4}_\lambda
    + F^{6,2}_\lambda)\label{w0_schematic}\\
    & & &\ \ \ \ \ \ \ \  +
    \dlambda^3\lambda^{-5}( F^{4,4}_\lambda + F^{6,2}_\lambda)
    \ ,\notag \\
    |J_\lambda| \ &\lesssim \ F^{4,4}_\lambda \ ,
    &|\partial_t J_\lambda| \ &\lesssim\ \dlambda\lambda^{-1}
    F^{4,4}_\lambda \ , \label{J_schematic}\\
    |\partial_t^2 J_\lambda| \ &\lesssim\ \big(|\ddot{\lambda}|\lambda^{-1}
     +  \dlambda^2\lambda^{-2}\big) F^{4,4}_\lambda \ , \notag\\
     \partial_t(A_\lambda) \ &= \ \dlambda F_\lambda^{7,9} ,
     &|\partial_t^2(A_\lambda)| \ &\lesssim \ \big(|\ddot{\lambda}| +
     \dlambda^2\lambda^{-1}\big) F^{7,9}_\lambda \ , \label{A_schematic}\\
     A_\lambda (F^{m,n}_\lambda) \ &= \ \lambda F^{m-1,n+1}_\lambda
      \ . \notag
\end{align}
We note here that we are assuming $1\leqslant m$ for the last
identity on line \eqref{A_schematic}.
\end{subequations}

\ret\ret

\section{Orbital Stability}\label{orbit_sect}

We now begin with the first step in our proof of Theorem
\ref{mod_thm}. This is to show that one can make a rough
decomposition of the full wave-map $\phi$ into a bulk piece which is
a rescaled soliton, plus a small remainder which we can estimate
in a certain energy space. Of particular importance to us will be
that we can construct this decomposition in such a way as to retain
a certain orthogonality between the bulk piece and the small
``radiation term''. Of course, this is precisely the modulational
approach to orbital stability first pioneered by M. Weinstein in his
study of the non-linear Schr\"odinger equation (see \cite{W_mod}).
What we intend to prove is the following:\\

\begin{lem}[Orbital stability with orthogonal decomposition]\label{orbit_lem}
Suppose that $\phi$ is a solution to the problem \eqref{EL_eqs_red},
and suppose that initially the Cauchy data for $\phi$ decomposes as:
\begin{align}
    \phi(0) \ &= \ I_{\lambda_0} + u_0 \ ,
    &\partial_t\phi\, (0) \ &= \ \dot{\phi}_0 \ , \label{C_data}
\end{align}
Assume that  the energy satisfies $E[\phi] = 4\pi k + \epsilon^2$,
with  $\epsilon$ chosen small enough. Then as long as the solution
$\phi$ exists there is a time dependent parameter
$0<\lambda(t)<\infty$, with $\lambda(0)=\lambda_0$, and the property
that the following conditions hold for all times of existence provided that
they hold initially:
\begin{align}
    E_0[u] \ &= \ \frac{1}{2}
    \int_{\mathbb{R}^+}\ \left[ (\partial_t\phi)^2
    + (\partial_r u)^2 + \frac{k^2}{r^2}u^2\right]\
    rdr \ \lesssim \ \epsilon^2 , \label{energy_size}\\
    0 \ &= \ \langle u(t) , J_{\lambda(t)} \rangle
    \ . \label{orth_cond}
\end{align}
Here $u$ is defined by the relation $u=\phi-I_\lambda$.
In addition we have that:
\begin{equation}
    \big |\frac{\dlambda}{\lambda^2}\big |
    \ \lesssim  \ \epsilon \ . \label{orbit_bound}
\end{equation}
Finally, we
remark that (conversely) the full wave map $\Phi$ is $C^\infty$ up
to any time $T$ as long as $\lambda(t) < \infty$ and
\eqref{energy_size} holds true for any $t\in [0,T]$.
\end{lem}\ret

\begin{proof}[Proof of Lemma \ref{orbit_lem}]
The proof essentially reduces to defining an appropriate
equation for the evolution of  $\lambda(t)$, basic existence and
uniqueness theory of ODEs, followed by the coercive estimate
\eqref{c_app1} proved in Appendix \ref{coercive_app}. We note here
that the last remark of Lemma \ref{orbit_lem} follows from the
local ``small energy implies regularity" statement for symmetric wave-maps, and is
for instance contained in Theorem \ref{C_th}.\\

We will use here  the following observation: Notice that if
$\lambda\to\infty$ or $\lambda\to 0$ on some time interval such that
\eqref{energy_size} holds, then the full wave-map $\Phi$ must break
down also on that time interval because its energy concentrates at
$r=0$ or $r=\infty$ respectively (in fact, by finite speed of
propagation it cannot happen that $\lambda\to 0$ in finite time, so
any ``blow-up'' of this type must occur at $t=\infty$).\\

Now, for a strictly positive real valued function of time
$\lambda(t)$, we define the equation:
\begin{equation}
     \dlambda \Big(2\langle I,J \rangle +
     \langle\phi(\lambda^{-1}r) , r\partial_r J\rangle
     \Big)  \ = \ -\langle\phi_t ,
     J_{\lambda}\rangle\lambda^3 \ . \label{basic_ODE}
\end{equation}
Notice that as long as the wave-map $\phi$ exists and is smooth,
this equation is of the form $\alpha(\lambda,t)\dlambda =
\beta(\lambda,t)$ for two $C^1$ functions $(\lambda,t)$ of the variables $\alpha,\beta$. We
now construct $\lambda(t)$ from \eqref{basic_ODE} via a simple
bootstrapping procedure. Our goal is to provide a strict lower bound
for $|\alpha|$ so that we may simply apply the usual existence
theory to \eqref{basic_ODE} which then produces $\lambda(t)$.\\

First of all, notice that a simple calculation involving the identity:
\begin{equation}
    \langle I, r\pa_r J\rangle \ = \
    -2\langle I, J\rangle - \langle J,J\rangle  \ , \notag
\end{equation}
shows that the bound \eqref{energy_size} initially implies that:
\begin{equation}
    \big|\alpha(\lambda_0,0) + \langle J , J \rangle\big|
    \ \leqslant  \ C\epsilon \ , \notag
\end{equation}
for some, possibly large, constant $C$ which we choose in a moment.
Therefore, by continuity there exists a small time $T^*$ such that
the solution $\phi$ exists and is regular on $[0,T^*]$, and one has
that a solution $\lambda(t)$ to equation \eqref{basic_ODE} exists,
is contained in $(0,\infty)$, and obeys the bounds:
\begin{align}
    \big|\alpha(\lambda(t),t) + \langle J , J \rangle
    \big| \ \leqslant \ 2C\epsilon \ . \label{boot_bound}
\end{align}
The heart of the matter is now the following: We will show that if
$T^*$ is \emph{any} time such that the above holds (i.e. existence
for $\phi$ and $\lambda$ and the bound \eqref{boot_bound}), then we
must also necessarily have the conditions \eqref{orth_cond} and
\eqref{energy_size}, \emph{as well as} the following improvement of
\eqref{boot_bound}:
\begin{align}
    \big|\alpha(\lambda(t),t) + \langle J , J \rangle
    \big| \ \leqslant \ C\epsilon \ . \label{better_bound}
\end{align}
By continuing this process, we will have shown that on any time
interval $[0,T^*]$ such that the solution $\phi$ exists and is
regular, there is a continuous solution of \eqref{basic_ODE}
$\lambda(t)\in(0,\infty)$ such that $\lambda(0)=\lambda_0$ and
\eqref{orth_cond}-\eqref{energy_size} holds (assuming, of course,
that these conditions hold initially).\\

We now show that existence and \eqref{boot_bound}, implies
\eqref{orth_cond} and  \eqref{energy_size}, and that these two
together imply the improved bound \eqref{better_bound}. Everything
rests on the orthogonality condition \eqref{orth_cond}. We define
$u=\phi-I_\lambda$ and use a few simple integration by parts to
compute that:
\begin{align}
    \frac{d}{dt} \langle u , J_\lambda\rangle \ &= \
     \langle u_t , J_\lambda\rangle +
     \frac{\dlambda}{\lambda} \langle u , (r\partial_r J)_\lambda\rangle
     \ , \notag\\
     &= \ \frac{\dlambda}{\lambda^3}
     \Big(2\langle I,J \rangle +
     \langle\phi(\lambda^{-1}r) , r\partial_r J\rangle
     \Big)   + \langle\phi_t ,
     J_{\lambda}\rangle \ , \notag\\
     &= \ 0 \ , \ \notag
\end{align}
where the last line follows from the assumption that $\lambda$
solves \eqref{basic_ODE}. Therefore, since \eqref{orth_cond} holds
initially, it holds on $[0,T^*]$.\\

It remains to show \eqref{energy_size}, and that this bound implies
\eqref{better_bound}. In fact, given the form of $\alpha$ this
latter implication is immediate from \eqref{energy_size} and the
Cauchy-Schwartz inequality, where the sufficiently large constant $C$ is chosen
according to the implicit constant appearing on line
\eqref{energy_size}.\\

Thus, we have reduced things to proving that, assuming that orthogonality
condition \eqref{orth_cond}
holds,  \eqref{energy_size} also
holds with a fixed implicit constant which \emph{does not} depend on
the size of $\lambda$  or $C$. This is where the condition $E[\phi]-E[I_\lambda]
=\epsilon^2$ enters. Notice that this equality follows from the
conservation of energy and our assumption $E[\phi] = 4\pi k +
\epsilon^2$. Computing the energy difference and using line
\eqref{difference_energy} we have that:
\begin{align}
        \epsilon^2 \ &= \ \pi\,
        \int_{\mathbb{R}^+}\ \left[
        (\partial_t\phi)^2 + \big(\partial_r\phi-\frac{k}{r}
        \sin(\phi) \big)^2
        \right]\ rdr \ , \notag\\
        &= \ \pi\,
        \int_{\mathbb{R}^+}\ \left[
        (\partial_t\phi)^2 + \big(\partial_r u-\frac{k}{r}
        \cos(I_\lambda)u \big)^2
        \right]\ rdr \ - \ \pi\,\mathcal{R}_\lambda(u) \ , \notag
\end{align}
where we are defining higher order integrated non-linearity:
\begin{equation}
    \mathcal{R}_\lambda(u) \ = \
    \int_{\mathbb{R}^+}\  \big(\partial_r \phi-\frac{k}{r}\sin(\phi) \big)^2\ rdr
    \ - \ \int_{\mathbb{R}^+}\  \big(\partial_r u-\frac{k}{r}
    \cos(I_\lambda)u \big)^2\ rdr \ . \notag
\end{equation}
Now, using the first coercive bound \eqref{c_app1} of Appendix
\ref{coercive_app}  we have the estimate:
\begin{equation}
     E_0[u]  \ \lesssim \ \epsilon^2 + |\mathcal{R}_\lambda(u)| \ . \label{E0_bound}
\end{equation}
Finally, using the simple algebraic formula for the difference of
squares, the equation \eqref{B_eq}, and writing the first few terms
in the Taylor series for $\sin(I_\lambda + u)$, we easily have the
nonlinear bound:
\begin{align}
    |\mathcal{R}_\lambda(u)| \ &\lesssim \ \int_{\mathbb{R}^+}\
    (|\partial_r u| + \frac{|u|}{r})\cdot(\frac{u^2}{r})\ rdr \ , \notag\\
    &\lesssim \ \big( E_0[u] \big)^\frac{3}{2} \ , \label{nonlin_M_bound}
\end{align}
where the last line follows from Cauchy-Schwartz and the Poincar\'e
type estimate:
\begin{equation}
    |u(r)|^2\ \leqslant \ 2\left(\int_{\mathbb{R}^+}\ (\partial_r u)^2
    \ rdr\right)^\frac{1}{2}
    \cdot\left(\int_{\mathbb{R}^+}\
    \frac{u^2}{r^2}\ rdr\right)^\frac{1}{2} \ . \label{energy_Poincare}
\end{equation}
The bounds \eqref{E0_bound}--\eqref{nonlin_M_bound} taken together
show that we may conclude \eqref{energy_size} for $\epsilon$ small
enough and some universal implicit constant (e.g. this can be shown
through another continuity argument). \\

Having established \eqref{energy_size}  and \eqref{better_bound}
we can easily show from \eqref{basic_ODE} that:
\begin{equation}
    \big |\frac{\dlambda}{\lambda^2} \big |
    \ \lesssim  \ \epsilon \ . \notag
\end{equation}
This concludes our demonstration of Lemma \ref{orbit_lem}.
\end{proof}\ret

\ret


\section{The effective evolution and the main blowup argument}\label{ODE_sect}

We now begin in earnest the proof of the main Theorem
\eqref{mod_thm}. This centers around computing a more effective form of the
modulation ODE \eqref{basic_ODE}. This will be followed by an ODE  analysis giving the
desired blow-up together with its asymptotic profile. The main technical
result of this section is the following:\\

\begin{prop}
[Refined structure for the modulation equation
\eqref{basic_ODE}]\label{main_prop} Consider the scaling parameter
$\lambda(t)$ which is defined through Lemma \ref{orbit_lem} and equation
\eqref{basic_ODE}. Suppose that the initial data \eqref{C_data} is
given according to Theorem \ref{mod_thm}. Then
on a time interval where $\lambda(t)\in (0,\infty)$ and $t\in [0,\epsilon^{-4}]$
the scaling parameter $\lambda(t)$ satisfies a first order ODE:
\begin{equation}
        \big(C_0 - \epsilon_1(t)\big)\dot{\lambda}(t) \ = \
        \epsilon_0\lambda^2(t) \ - \ \lambda^2(t)
        \int_0^t\
        \mathcal{E}(s)\ ds \ .
        \label{refined_main_ODE}
\end{equation}
Here $C_0=\langle J,J\rangle$,
and $\epsilon_0, \epsilon_1(t)$ and $\mathcal {E}(t)$ obey the conditions
($\epsilon_0$ is fixed):
\begin{align}
        |\epsilon_0 - \frac{\epsilon}{\pi}| \ &\lesssim \
        c_0\epsilon\ , \label{eps0_cond}\\
        |\epsilon_1(t)| \ &\lesssim \
        \epsilon\ , \label{eps1_cond}\\
        |\mathcal{E}(t)| \ &\lesssim \ c_0\epsilon^2 +
        c_0^\frac{1}{2}\sup_{0\leqslant s\leqslant t}
        {\dot{\lambda}^4}{\lambda^{-7}}(s) \ . \label{main_bounds}
\end{align}
Here $\epsilon,c_0$ are the small constants from line
\eqref{special_smallness}. In addition, one has the following
``structure bounds'' for the acceleration of $\lambda(t)$:
\begin{equation}
        \big| \ddot{\lambda}(t) - 2 \frac{\dot{\lambda}^2}{\lambda} (t)\big|
        \ \leqslant \ Cc_0^\frac{1}{2}\frac{\dot{\lambda}^2}{\lambda} (t) \ + \
        C\big( c_0\epsilon^2 +
        \sup_{0\leqslant s\leqslant t}
        \frac{\dot{\lambda}^4}{\lambda^7}(s) \big)\lambda^2
        \ . \label{ldd_structure}
\end{equation}
We remark that the constant $C$ is universal, i.e. independent of $c_0$ and $\epsilon$.
\end{prop}\ret

\begin{rem}
The reason the explicit constant $C$ appears in the estimate
\eqref{ldd_structure} is merely a notational convenience. The bound
\eqref{ldd_structure}  will be one of our main bootstrapping assumptions
in the sequel.
\end{rem}\ret

The remainder of this section is divided into two parts. First, we
will show that the identities and estimates
\eqref{refined_main_ODE}--\eqref{ldd_structure} imply that the
parameter $\lambda(t)$ goes to infinity at some time
$T^{**}\in[0,\epsilon^{-4}]$ as long as $\epsilon$ is chosen small
enough. We then establish the bounds \eqref{blowup_rate} on the
asymptotic rate of $\lambda(t)$ as $t\to T^{**}$. Finally, in the
last subsection we state the ``Main Estimate'' (a certain fixed time energy estimate) of our paper,
and we use it to derive all of the assumptions
\eqref{refined_main_ODE}--\eqref{ldd_structure}. This
main technical estimate will be the subject of the final  section of
the paper.\\

\ret

\subsection{Proof of the blowup  and the  universal
bound for $\lambda(t)$}\label{ODE_an_sect}

Using the equation \eqref{refined_main_ODE}, as well as the
assumptions \eqref{eps0_cond}--\eqref{ldd_structure} we now show
that $\lambda$ must blowup in finite time. The basic idea is the
following: without the contribution of the integral on the right
hand side of \eqref{refined_main_ODE} the desired blow-up would
occur in finite time $\sim C_0 \epsilon_0^{-1}$ in a self-similar
Riccati fashion. Therefore, the only problem is that one must
guarantee the integral term (which  adds a negative\footnote{ That
the contribution is overall a negative one follows from Theorem
\ref{C_th}.} contribution, creating a delay effect) does not
interfere to the extent that $\dot{\lambda}$ is driven to zero too
quickly before the blowup can occur. The fact that {\it a priori}
(i.e. again by Theorem \ref{C_th}) at any supposed blow-up time one
\emph{must} have that $\dlambda\lambda^{-2} \to 0$, indicates a very
delicate balancing between the two main terms on right hand side of
\eqref{refined_main_ODE}. This constitutes one of the main technical
difficulties in this paper, and why many of the estimates which
appear in the sequel are so involved. To establish the needed
control, we will first show that the initial self-similar behavior,
approximated by the ODE $C_0\dlambda=\epsilon_0\lambda^2$, forces
$\lambda(t)$ into a different monotonic regime, where in particular
the terms $\dot{\lambda}^4\lambda^{-7}$ dominates the error estimate
for $\mathcal{E}$ (i.e. \eqref{main_bounds}). At that point  blow-up
is assured. Our last task is then to analyze the balance of the
integral and $\epsilon_0$ terms on the  right hand side of
\eqref{refined_main_ODE}, and to derive the precise blow-up bounds
\eqref{blowup_rate} from this. Along the way, we will show that all
of this can be accomplished before the time interval
$[0,\epsilon^{-4}]$ expires, so that we still have access to all of
the structure included
Proposition \ref{main_prop}.\\

We now proceed with the details outlined above. The first
main thing is control the size of the
interval where one cannot guarantee good bounds on $\mathcal{E}(s)$.
The key to this is to show that $\dot{\lambda}^4\lambda^{-7}$
becomes monotonic soon enough and with enough force to cover the
constant $c_0\epsilon^2$ which is lost on line
\eqref{ldd_structure}. Luckily it is not hard to show that these two
things happen at essentially the same time. Computing the time
derivative and then using the identity \eqref{ldd_structure} we
have:
\begin{equation}
        \frac{d}{dt}\left(\frac{\dot{\lambda}^4}{\lambda^7}\right)
        \ = \ \frac{\dot{\lambda}^5}{\lambda^8}
        + O\Big( c_0^\frac{1}{2} \frac{\dlambda^2}{\lambda^3} +
        c_0\epsilon^2 +
        \sup_{0\leqslant s\leqslant t}
        \frac{\dot{\lambda}^4}{\lambda^7}(s) \Big)\frac{\dot{\lambda}^3}
        {\lambda^5} \ . \label{time_der}
\end{equation}
We now let $T_0$ be the first time such that:
\begin{equation}
    \dot{\lambda}^4\lambda^{-7}(T_0)\ = \ C c_0^\frac{1}{2}\epsilon^2
    \ , \label{Tstar_bound}
\end{equation}
for some large constant $C$ which is larger than twice the implicit
constants in \eqref{main_bounds}. First we argue that such a time
must occur.\\

If such a time does not occur, then if $c_0$ is chosen small enough
(we remind the reader that this is done by simply choosing initial
data according to \eqref{special_smallness}, and does not affect the
size of constants in estimates like \eqref{main_bounds}) one sees
immediately from \eqref{main_bounds} and \eqref{refined_main_ODE}
that the following bound holds for all times:
\begin{equation}
    \frac{\dot{\lambda}}{\lambda^2}(t) \ \geqslant \ \frac{1}{2C_0}\epsilon_0
    - \int_0^t \ C c_0\epsilon_0^2 ds \ , \label{eq:star}
\end{equation}
where $C$ is fixed and independent of $c_0,\epsilon$. A simple argument,
which we leave to the reader, shows that if $\lambda(0)=1$ then this
last inequality implies that $\lambda\to\infty$ at some finite time
$T^{**}\leqslant 4C_0\epsilon_0^{-1}$ and that in addition:
\begin{equation}
    \frac{\dlambda}{\lambda^2}(t) \ \geqslant \
     \frac {\epsilon_0}{4C_0} \ , \notag
\end{equation}
for all $t\in [0,T^{**})$. In that case however we would also
have that:
\begin{equation}
    \frac{\dlambda^4}{\lambda^7}(t) \ = \
    \left(\frac{\dlambda}{\lambda^2}(t)\right)^4 \lambda(t)
     \ \geqslant \  \left(\frac{\epsilon_0}{4C_0}\right)^4
     \lambda(t) \ \to  \ \infty\quad\text{as}\quad t\to T^{**} \ , \notag
\end{equation}
which shows that in fact the time $T_0$ defined above must occur.
What's more, a simple analysis of the previous argument shows that
this time must also satisfy the conditions:
 \begin{align}
    T_0 \ &< \ 4 C_0 \epsilon_0^{-1} \ ,
    &\dlambda\lambda^{-2}(T_0) \ &\geqslant \ \frac{1}{4C_0}\epsilon_0
    \ . \label{time_bounds}
\end{align}
Notice  it is also clear that for all times $t\in [0,T_0]$ we have
$\dlambda(t)>0$.\\

Now, by a direct application of the orbital stability
bound \eqref{orbit_bound} we have that for any time where the
solution exists there is the bound:
\begin{equation}
    \frac{\dot{\lambda}^4}{\lambda^7}(t) \ \lesssim \
    \epsilon^2\frac{\dlambda^2}{\lambda^3}(t) \ . \notag
\end{equation}
Applying this on the right hand side of \eqref{time_der} we
see that for all $t\geqslant  T_0$,
and as long as $\dlambda^4\lambda^{-7}(t)\geqslant Cc_0\epsilon^2$, we have
that:
\begin{align}
    \frac{d}{dt}\left(\frac{\dot{\lambda}^4}{\lambda^7}\right)(t)
    \ &\geqslant \ (1-Cc_0^\frac{1}{2}) \frac{\dot{\lambda}^5}{\lambda^8}(t)
    \  > \ 0 \ . \notag
\end{align}
By bootstrapping this argument, we see that $\dlambda^4\lambda^{-7}$ is monotonically
increasing for all times
$t\geqslant T_0$.\\

Before continuing, with the proof of blowup, we pause for a moment to upgrade the bound
\eqref{Tstar_bound}. This will be used in a crucial way in the sequel (see the proof of Proposition
\ref{eng_w_prop} in Section \ref{mor_sect}). We claim that there exists a time $T^*=T_0+O(1)$ such that the
following improvement of  \eqref{Tstar_bound} is valid:
\begin{equation}
    \dot{\lambda}^4\lambda^{-7}(T^*)\ = \  c_0\epsilon \ . \label{better_Tstar_bound}
\end{equation}
Again by contradiction, if this were not the case by equation \eqref{refined_main_ODE}, the
definition \eqref{Tstar_bound} of $T_0$, and the time bounds \eqref{time_bounds},
we would have a bound of the form:
\begin{equation}
    \frac{\dot{\lambda}}{\lambda^2}(t) \ \geqslant \ \frac{1}{4C_0}\epsilon_0
        - \int_{T_0}^t \ C c_0\epsilon_0 ds \ \geqslant \ \frac{1}{8C_0}\epsilon_0 \ , \notag
\end{equation}
for times $T^*=T_0+O(1)$ as long as $c_0\ll 1$. Integrating this, and again applying
\eqref{Tstar_bound} as well as the orbital stability bounds \eqref{orbit_bound}
to the term $\lambda^{-1}(T_0)$, we arrive at the inequality:
\begin{align}
    (t-T_0)\epsilon_0 \ &\lesssim \ \frac{1}{\lambda(T_0)} -
    \frac{1}{\lambda(t)} \ \lesssim \  \ c^{-1}_0 \epsilon^2  \ . \notag
\end{align}
Using the condition that $\epsilon\leqslant c_0^2$ we see that such a bound must expire
in $t-T_0 = O(1)$ time. Finally, notice that by using the time bounds \eqref{time_bounds}, the definition
of $T_0$ \eqref{Tstar_bound}, the definition of $T^*$ \eqref{better_Tstar_bound}, as well
as the relation $T^*=T_0+O(1)$, we may integrate the quantity $\dlambda^4\lambda^7$
over $[0,T^*]$ to obtain:
\begin{equation}
    \int_0^{T^*}\ \frac{\dlambda^4}{\lambda^7}(s)\ ds \ \leqslant \epsilon \ .
    \label{first_int_bound}
\end{equation}\ret

We now return to the main thread of the blowup argument. So far we have achieved the following.
There exists a time $T_0$ such that for $T_0\leqslant t$ we have:\\

\begin{itemize}
    \item $\dot{\lambda}^4\lambda^{-7}(t)$ is a
    monotonically increasing function.\\
    \item $\dot{\lambda}^4\lambda^{-7}(t)\
    \geqslant \ C c_0^\frac{1}{2}\epsilon^2$, where
    $C$ is at least twice the implicit constant in
    the bound \eqref{main_bounds}.\\
    \item The time $T_0$ is associated with the bounds \eqref{time_bounds}.
\end{itemize}\ret

\noindent We are now at the point where blowup with the rate bounds
\eqref{blowup_rate} is assured. From the above conditions,
 we have that for all times $t\geqslant T_0$:
\begin{equation}
     \mathcal{E}\ \lesssim \ c_0^\frac{1}{2}
    \dot{\lambda}^4\lambda^{-7} \ . \label{mon_bounds}
\end{equation}
Writing:
\begin{equation}
    \gamma(t) \ = - \epsilon_0 +
    \int_0^t\
    \mathcal{E}(s) \ ds \ , \notag
\end{equation}
we have that $\dot{\lambda}\lambda^{-2} \ \sim \  -
\gamma$. Differentiating $\gamma$ and using the bound
\eqref{mon_bounds} we see that:
\begin{equation}
    \dot{\gamma} \ \lesssim \ c_0^\frac{1}{2}\frac{\dot{\lambda}^4}{\lambda^7} \ \sim
    \ -\gamma^3 c_0^\frac{1}{2}\frac{\dot\lambda}{\lambda} \ . \notag
\end{equation}
By dividing through by $-\gamma^3$ and integrating both sides of this last inequality
over the interval $[T_0,t_1]$
we arrive at the estimate (we may assume that $0\leqslant -\gamma$ throughout this
argument, as will become apparent on the next line):
\begin{align}
    \frac{1}{\gamma^2(t_1)} \ &\lesssim \
    \frac{1}{\gamma^2(T_0)} + c_0^\frac{1}{2}\ln\left(\lambda(t_1)\right) - c_0^\frac{1}{2}\ln\left(\lambda(T_0)\right)\ , \notag\\
    &\lesssim \ \epsilon^{-2} + c_0^\frac{1}{2}\ln\left(\lambda(t_1)\right) , \notag
\end{align}
where to obtain the last inequality we've used the second bound on line
\eqref{time_bounds} as well as the fact that $1\leqslant \lambda(T_0)$.
We now recast this last expression in terms of
$\dot{\lambda}$:
\begin{equation}
    \frac{1}{\sqrt{\epsilon^{-2} + c_0^\frac{1}{2}\ln\left(\lambda(t_1)\right)}} \ \lesssim \
    \frac{\dlambda}{\lambda^2}(t_1) \ . \label{blowup_ineq}
\end{equation}
Integrating this last line over time intervals past $[0,T_0]$, we see that within
$O(\epsilon^{-1})$ time we must have  $\epsilon^{-2}\leqslant c_0^\frac{1}{2}\ln\left(\lambda(t)\right)$.
Therefore we are assured of the bound:
\begin{align}
    \lambda^{2-\delta}(t) \ &\lesssim \
    \frac{\lambda^2(t)}{\sqrt{c_0^\frac{1}{2}\ln\left(\lambda(t)\right)}} \ \lesssim \
    \dlambda(t) \ ,  &0 < \delta \ll 1 \ , \notag
\end{align}
inside of some interval $[0,C\epsilon^{-1}]$, for a uniform constant $C$. This is
a Ricatti type inequality, which easily implies that $\lambda(t)\to\infty$ in $O(1)$ time
starting with $1\leqslant \lambda$ at the first time where it holds.
Let us call the blowup time $T^{**}$.\\

Finally, we need to recover the rate bounds \eqref{blowup_rate}. By Remark \ref{lower_bound_rem}
we need only establish the upper bound. From \eqref{blowup_ineq} we have the inequality:
\begin{equation}
    c_0^{-\frac{1}{4}} \ \lesssim \
    \sqrt{\ln(\lambda(t))}\frac{\dlambda}{\lambda^2}(t) \ , \notag
\end{equation}
for times $t$ sufficiently close to $T^{**}$.
Making the substitution $\lambda(t)=e^{s^2}$ and
integrating from $t$ to the blowup time we have the estimate:
\begin{align}
    c_0^{-\frac{1}{4}}(T^{**} - t) \ &\lesssim \
    \ 2\int_{\sqrt{\ln(\lambda(t))}}^\infty\
    s^2 e^{-s^2}ds \ , \notag\\
    &= \ \frac{\sqrt{\ln(\lambda(t))}}{\lambda(t)} +
    \int_{\sqrt{\ln(\lambda(t))}}^\infty\  e^{-s^2}ds \ , \notag\\
    &= \ \frac{\sqrt{\ln(\lambda(t))}}{\lambda(t)} +
    \frac{O(1)}{\lambda(t)\sqrt{\ln(\lambda(t))}} \ , \notag\\
    &\lesssim \ \frac{\sqrt{\ln(\lambda(t))}}{\lambda(t)} \ . \notag
\end{align}
where the last line above follows from the well known asymptotics of
the error function. The above identity easily implies that as $t\to
T^{**}$ one has the bound:
\begin{equation}
    \lambda(t) \ \lesssim \ c_0^\frac{1}{4}
    \frac{\sqrt{|\ln(T^{**}-t)|}}{(T^{**}-t)} \ . \notag
\end{equation}
This establishes the upper bound in \eqref{blowup_rate}.\\

\ret

\subsection{Derivation of the Main ODE and its Structure}

We now derive the ODE \eqref{refined_main_ODE}, as well as all of
the accompanying structural assumptions
\eqref{eps0_cond}--\eqref{ldd_structure}. This will require a
certain fixed time energy estimate which is the main technical
estimate of the paper and will be proved in the following section. We start
by recomputing the modulation equation \eqref{basic_ODE} using
the splitting $\phi = I_\lambda + u$.
Differentiating the orthogonality relation \eqref{orth_cond} with
respect to time, we have the simple identity:
\begin{equation}
    \langle \dot{I}_\lambda , J_\lambda \rangle
    \ = \ \langle \phi_t , J_\lambda \rangle + \langle
    u,\dot{J}_\lambda\rangle  \ . \label{baby_mod}
\end{equation}
Differentiating one more time, and rearranging things with a little help
from the equation $\phi_{tt} + H_\lambda u =\mathcal {N}(u)$, we
have that:
 \begin{align}
    \langle \ddot{I}_\lambda , J_\lambda \rangle \ &= \
    2\langle\partial_t u , \dot{J}_\lambda \rangle +
    \langle u,\ddot{J}_\lambda \rangle
    + \langle\mathcal{N}(u) , J_\lambda \rangle \ , \label{two_der_iden1}\\
    &= \
    2\partial_t \langle u , \dot{J}_\lambda \rangle -
    \langle u,\ddot{J}_\lambda \rangle
    + \langle\mathcal{N}(u) , J_\lambda \rangle \ ,
    \label{two_der_iden2}
\end{align}
where $\mathcal{N}(u)$ denotes the nonlinearity:
\begin{align}
    \mathcal{N}(u) \ &= \
    \frac{k^2\sin(2I_\lambda)}{2r^2}
        \cdot(1-\cos(2u)) \ +\  \frac{k^2\cos(2I_\lambda)}{r^2}
        \cdot(u-\frac{1}{2}\sin(2u)) \ , \notag\\
    &= \ \frac{k^2\sin(2I_\lambda)}{r^2}\cdot u^2 + \td{\mathcal{N}}(u)
    \ . \label{N_breakdown}
\end{align}
Notice that:
\begin{equation}
    |\td{\mathcal{N}}(u)| \ \lesssim \ \frac{|u|^3}{r^2} \ . \label{td_N_bound}
\end{equation}
Next, a short computation shows that we have the identity (recall the
definition of $J_\lambda$ from line \eqref{JI_defs}):
\begin{align}
    \langle \ddot{I}_\lambda , J_\lambda \rangle \ = \
    \left(\frac{\ddot{\lambda}}{\lambda} -
    2\frac{\dlambda^2}{\lambda^2}\right)
    \langle J_\lambda , J_\lambda\rangle
    = \ C_0\frac{d}{dt}\left(\frac{\dlambda}{\lambda^2}\right)
    \lambda^{-1} \ \ . \notag
\end{align}
where we have set $\langle J , J \rangle = C_0$. Therefore, from
these last two lines as well as the identities
\eqref{two_der_iden1}--\eqref{two_der_iden2}, we have our two main
structural equations:
\begin{align}
    \ddot{\lambda} - 2\frac{\dlambda^2}{\lambda} \ &= \ C_0^{-1}\left(
    2\langle\partial_t u , \dot{J}_\lambda \rangle +
    \langle u,\ddot{J}_\lambda \rangle
        + \langle\mathcal{N}(u) , J_\lambda \rangle \right)\lambda^3
    \ , \label{accl_iden}\\
        C_0\frac{d}{dt}\left(\frac{\dlambda}{\lambda^2}\right) \ &= \
    2\partial_t \big[\langle u , \dot{J}_\lambda \rangle \lambda\big]
    - 2 \langle u , \dot{J}_\lambda \rangle \dlambda
    - \langle u,\ddot{J}_\lambda \rangle\lambda
        + \langle\mathcal{N}(u) , J_\lambda \rangle\lambda \ . \label{first_iden}
\end{align}
The first equation \eqref{accl_iden} is sufficient for us to prove
the bound \eqref{ldd_structure}. The second equation
\eqref{first_iden} will yield \eqref{refined_main_ODE} upon
integration. Doing this over a time interval $[0,t]$ we have the
identity (recall that $\lambda(0)=1$):
\begin{multline}
    C_0 \frac{\dlambda}{\lambda^2}(t)  \ -  \ 2
    \langle u(t) , \dot{J}_{\lambda(t)} \rangle \lambda(t) \\
    = \ \dlambda(0)\big(C_0 - 2 \langle u(0) , rdr J \rangle\big)
    \ - \ \int_0^t\  \Big(C_*\frac{\dlambda^4}{\lambda^7}(s) + \mathcal{E}(s)
    \Big)\ ds \ , \label{raw_int_iden}
\end{multline}
where we define the constant $C_*$ as follows:
\begin{multline}
    C_* \ = \
    -\ k^2 \Big\langle \frac{\left(a J_\lambda + b(r^2 J)_\lambda\right)^2}
    {r^2}  , \sin(2I_\lambda)\cdot J_\lambda \Big\rangle \\
    + \ \lambda^2 \Big\langle  \left(a J_\lambda +
    b(r^2 J)_\lambda\right) , (r\partial_rJ)_\lambda \Big\rangle
    - \ \lambda^2 \Big\langle r\partial_r \left(a J_\lambda +
    b(r^2 J)_\lambda\right) , (r\partial_rJ)_\lambda \Big\rangle \ . \label{C_exp}
\end{multline}
In Appendix A, it will be shown that $C_*=0$. This indicates that the precise rate
of blowup in the inequality \eqref{blowup_rate} is quite delicate. We'll return to this
in a later work.\\

The error term $\mathcal{E}(s)$ is given by the expression:
\begin{multline}
        \mathcal{E} \ = \  2 \langle w , \dot{J}_\lambda \rangle \dlambda
        \ + \ \langle w,\ddot{J}_\lambda \rangle\lambda
        \ + \ \langle w_0,\big(\ddot{\lambda} -2\frac{\dlambda^2}{\lambda}
        \big)(r\partial_r J)_\lambda \rangle \\
        - \ \langle\frac{w\cdot(2w_0+w)}{r^2} ,\sin(2I_\lambda)\cdot J_\lambda \rangle\lambda
        \ - \ \langle\td{\mathcal{N}}(u) , J_\lambda
        \rangle\lambda \ . \label{cal_E_def}
\end{multline}
Here the terms $w_0,w$ refer to the decomposition on line
\eqref{w0_def} above. We list this here again for the convenience of
the reader:
\begin{align}
    u \ = \ \frac{\dlambda^2}{\lambda^4}\left(a J_\lambda + b(r^2 J)_\lambda\right)
    + w = \ w_0 + w  \ , \label{u_decomp}
\end{align}
where the constants $a,b$ are derived on lines \eqref{ab_defs}. Notice that
the constant $C_*$ defined by \eqref{C_exp} arises from the expression
(and a few integrations by parts):
\begin{multline}
    C_* \dlambda^4\lambda^{-7} \ = \
    2\langle w_0,\dot J_\lambda\rangle\dlambda \ + \ \langle w_0,\ddot J_\lambda\rangle\lambda\\
    - \ \langle\frac{k^2w_0^2}{r^2}, \sin (2I_\lambda) \cdot J_\lambda\rangle\lambda
     \ - \ \langle w_0,\big(\ddot{\lambda} -2\frac{\dlambda^2}{\lambda}
        \big)(r\partial_r J)_\lambda \rangle \ . \notag
\end{multline}\ret

Before commencing with the proof of the estimates
\eqref{main_bounds} and \eqref{ldd_structure}, let us first derive
from \eqref{raw_int_iden} the identity \eqref{refined_main_ODE}, and
also the conditions \eqref{eps0_cond}--\eqref{eps1_cond}. First of
all, notice that from the orbital stability bounds
\eqref{energy_size} we have that:
\begin{align}
    \big| \langle u , \dot{J}_\lambda \rangle \lambda\big|
    \ &\lesssim \ \frac{\dlambda}{\lambda}\
    \lp{r^{-1}u}{L^2(rdr)}\cdot\lp{(r^2\partial_r J)_\lambda}{L^2(rdr)}
    \ , \notag\\
    &\lesssim\ \epsilon\frac{\dlambda}{\lambda^2} \ . \notag
\end{align}
Therefore, we may define $\epsilon_1$ on the left hand side of line
\eqref{refined_main_ODE} as $\epsilon_1 = 2
\lambda^3\dlambda^{-1}\langle u , \dot{J}_\lambda \rangle$ and we
have \eqref{eps1_cond}.\\

Similarly, at the initial time, the identity \eqref{baby_mod} gives
the relation (we are assuming $\lambda(0)=1$):
\begin{equation}
    \dlambda(0) \langle J , J \rangle
        \ = \ \langle \phi_t(0) , J \rangle + \dlambda(0)\langle
        u(0) , r\partial_r J\rangle  \ . \notag
\end{equation}
Substituting into this last relation the initial data
\eqref{special_C_data}, and using the smallness condition
\eqref{special_smallness} we easily have that:
\begin{equation}
    \dlambda(0)\big(C_0 - 2 \langle u(0) , rdr J \rangle\big) \ = \
    \frac{\epsilon}{\pi} + O(c_0\epsilon) \ , \notag
\end{equation}
which gives the condition \eqref{eps0_cond}. Finally, notice that this last line
also implies the initial expansion:
\begin{equation}
    \dot{\lambda}(0) \ = \ \frac{\epsilon}{\pi\langle J , J \rangle} + O(c_0\epsilon) \ . \notag
\end{equation}
Plugging this into the first term on the RHS of formula \eqref{accl_iden},
and using the bounds \eqref{special_smallness} on our chosen initial data \eqref{special_C_data}
to estimate the remaining terms, we see that \eqref{ldd_structure}
holds for the initial time $t=0$.\\
\\ \\ \\


It remains for us to derive the bounds \eqref{main_bounds} and
\eqref{ldd_structure} from the identities \eqref{cal_E_def} and
\eqref{accl_iden} respectively. This will be done through a
bootstrapping process and the use of a special energy estimate for
the function $w$ appearing in those expressions. This brings us to
the main technical estimate of our paper which is the following:\\

\begin{prop}[Main technical estimate]\label{tech_prop}
Let $u=w_0 + w$ be the decomposition from line \eqref{u_decomp}.
Next, let us assume that the estimate \eqref{ldd_structure} holds
with constant $2C$, that is:
\begin{equation}
        \big| \ddot{\lambda}(t) - 2 \frac{\dot{\lambda}^2}
        {\lambda} (t)\big|
        \ \leqslant \ 2Cc_0^\frac{1}{2}\frac{\dot{\lambda}^2}
        {\lambda} (t) \ + \
        2C\big( c_0\epsilon^2 +
        \sup_{0\leqslant s\leqslant t}
        \frac{\dot{\lambda}^4}{\lambda^7}(s) \big)\lambda^2 \ .
        \label{ldd_structure_boot}
\end{equation}
Then as long as the parameter $\lambda(t)$ is monotonically
non-decreasing, one has the following fixed time energy type
estimate for $t\in[0,\epsilon^{-4}]$:
\begin{multline}
        \int_{\mathbb{R}^+}\
        \lambda^{-1}\frac{(\lambda r)^\delta}{1+r^\delta}
        \left[(LA_\lambda w)^2 + \frac{(A_\lambda w)^2}{r^2} \right](t)\ rdr
        \ \lesssim \ c_0^2\epsilon^2 +
        \epsilon\ \sup_{0\leqslant s \leqslant t}\ \frac{\dlambda^4}{\lambda^7}(s)
        \ . \label{special_energy}
\end{multline}
Here the implicit constant depends on $C$ from line
\eqref{ldd_structure_boot} above, but is independent of $\epsilon$
and $c$ from line \eqref{special_smallness}. Also, the operator
$A_\lambda$ is defined on line \eqref{A_ops} above. Lastly,
$L=\partial_t + \partial_r$ is the outgoing null derivative.
\end{prop}\ret

\begin{rem}
Observe that by a direct application of estimate \eqref{c_app2} of
Appendix \ref{coercive_app} that the bound on line
\eqref{special_energy} implies:
\begin{equation}
    \int_{\mathbb{R}^+}\
        \lambda^{-1}\frac{(\lambda r)^\delta}{1+r^\delta}
        \frac{ w^2}{r^4}(t)\ rdr \ \lesssim \ c_0^2\epsilon^2 +
    \epsilon\ \sup_{0\leqslant s \leqslant t}\
    \frac{\dlambda^4}{\lambda^7}(s) \ . \label{special_energy_w}
\end{equation}
This will be used many times in the sequel.
\end{rem}\ret

We now prove the bounds \eqref{main_bounds} and
\eqref{ldd_structure}. This will be done separately and in reverse
order. To prove the second estimate \eqref{ldd_structure}, it will
suffice for us to demonstrate the following set of estimates:\\

\begin{lem}[Estimates for \eqref{ldd_structure}]\label{accl_lem}
Assuming the bootstrapping estimate \eqref{ldd_structure_boot} and
the results of Proposition \ref{tech_prop}, one has the following
estimates where the implicit constant depends on line
\eqref{special_energy}:
\begin{align}
        \big|
        \langle\partial_t u , \dot{J}_\lambda \rangle\lambda^3 \big|
        \ &\lesssim \  c^\frac{3}{4}_0\frac{\dlambda^2}{\lambda} +
        c_0^\frac{1}{4}\big(c_0\epsilon^2 + \sup_{0\leqslant s\leqslant t}
        \frac{\dot{\lambda}^4}{\lambda^7}(s) \big)\lambda^2\ , \label{hardest}\\
        \langle u,\ddot{J}_\lambda \rangle\lambda^3
        \ &= \ \eta_1(t) \ddot{\lambda} + \eta_2(t)\frac{\dlambda^2}{\lambda}
        \ , \label{easiest}\\
        \big| \langle\mathcal{N}(u) , J_\lambda \rangle \lambda^3 \big| \ &\lesssim \
        \epsilon^2\frac{\dlambda^2}{\lambda} + \big(c_0^2\epsilon^2 + \epsilon
        \sup_{0\leqslant s\leqslant t} \frac{\dlambda^4}{\lambda^{7}}(s)\big)
        \lambda^2 \ , \label{medium}
\end{align}
where we also have the bounds\footnote{The implicit constants in \eqref{eta_bounds} depend only on the orbital stability bound
\eqref{energy_size} and thus are independent of the bootstrap constant $C$.}:
\begin{align}
    |\eta_1|\ &\lesssim \ \epsilon \ ,
    &|\eta_2|\ &\lesssim \ \epsilon \ . \label{eta_bounds}
\end{align}
In particular, for $\epsilon^\frac{1}{2}\leqslant c_0$ small enough, we have
that \eqref{ldd_structure} holds.
\end{lem}\ret

\begin{proof}[Proof of estimate \eqref{hardest}]
This is the most involved of the above three estimates. To prove
this, we begin by isolating the explicit piece involving $w_0$
from the expansion \eqref{u_decomp}. Thus, our first task is to
prove that:
\begin{equation}
    \big|
        \langle\partial_t w_0 , \dot{J}_\lambda \rangle
        \lambda^3 \big|
        \ \lesssim \  \epsilon^2 \frac{\dlambda^2}{\lambda}
        \ .  \label{hardest_w0}
\end{equation}
To do this we will employ the bootstrapping assumption
\eqref{ldd_structure_boot}. We will also use the abstract function
notation replacements from lines \eqref{schematic_not}. Doing this,
we see  from lines \eqref{w0_schematic} and \eqref{J_schematic} that
we may write:
\begin{equation}
    \langle\partial_t w_0 , \dot{J}_\lambda \rangle\lambda^3 \ = \
    \dlambda^4\lambda^{-3}\langle F^2_\lambda , F^4_\lambda \rangle +
    \frac{\ddot{\lambda}\dlambda^2}{\lambda^2}
    \langle F^2_\lambda , F^4_\lambda \rangle \ . \label{dtw0_prep}
\end{equation}
Notice that from the bootstrapping assumption
\eqref{ldd_structure_boot} and the estimate \eqref{orbit_bound}, as
well as the monotonicity of $\lambda$ we have that:
\begin{align}
        |\ddot{\lambda}| \  &\lesssim \  \dlambda^2\lambda^{-1} +
        2C\big( c_0\epsilon^2 +
        \sup_{0\leqslant s\leqslant t}
        \dlambda^4\lambda^{-7}(s) \big)\lambda^2 \ , \notag\\
        &\lesssim \
        \dlambda^2\lambda^{-1} + \epsilon^2 \lambda^3 \ \  \lesssim \ \
        \epsilon^2 \lambda^3 \ . \notag
\end{align}
Therefore, plugging this into line \eqref{dtw0_prep} and using that
$|\langle F^2_\lambda , F^4_\lambda \rangle|\lesssim \lambda^{-2}$,
we have the bound (using again \eqref{orbit_bound}):
\begin{align}
    \big|
        \langle\partial_t w_0 , \dot{J}_\lambda \rangle\lambda^3 \big|
     \ &\lesssim  \ \dlambda^4\lambda^{-5} +
    \epsilon^2 \dlambda^2\lambda^{-1} \ , \notag\\
    &\lesssim  \
    \epsilon^2 \dlambda^2\lambda^{-1} \ . \notag
\end{align}
This establishes \eqref{hardest_w0} and therefore \eqref{hardest}
for the $w_0$ portion of $u$.\\

We shall now prove that:
\begin{equation}
        \big| \langle\partial_t w , \dot{J}_\lambda \rangle
        \lambda^3 \big|
        \ \lesssim \  c^\frac{3}{4}_0\frac{\dlambda^2}
        {\lambda} +
        c_0^\frac{1}{4}\big(c_0\epsilon^2 +
        \sup_{0\leqslant s\leqslant t}
        \frac{\dot{\lambda}^4}{\lambda^7}(s) \big)
        \lambda^2\ . \label{hardest_w}
\end{equation}
The complication in this estimate stems from the fact that the energy
estimate \eqref{special_energy} does not provide control of the time
derivative of $w$. In addition, \eqref{hardest_w} is a fixed time estimate
and therefore it is not amenable to the procedure of
integrating out the time derivative as was
done on the line \eqref{first_iden} above. However we will be able to exploit
a structure of the inner product in this expression and convert the $\pa_t w$
into $L A_\lambda w$ derivative, which appears in \eqref{special_energy}.
First of all, we
write this as:
\begin{equation}
    \langle\partial_t w , \dot{J}_\lambda \rangle \ = \
    \dlambda\lambda^{-1}
    \langle\partial_t w , r \partial_r J_\lambda \rangle \ . \notag
\end{equation}
We now employ the following identity:
\begin{equation}
     A_\lambda^*(r J_\lambda) \ = \ 2J_\lambda + 2r\partial_r J_\lambda +
     rA_\lambda J_\lambda \ , \notag
\end{equation}
along with the conditions $A_\lambda J_\lambda=0$ and $\langle w,J_\lambda\rangle=0$
to write:
\begin{align}
        \langle\partial_t w , r \partial_r J_\lambda \rangle \ &= \
        \frac{1}{2} \langle\partial_t w , A_\lambda^*
        (r J_\lambda)\rangle
        - \langle \partial_t w , J_\lambda \rangle \ , \notag\\
        &= \ \frac{1}{2} \langle \partial_t A_\lambda w ,
        (r J)_\lambda\rangle\lambda^{-1}
        + \frac{1}{2} \langle [A_\lambda , \partial_t] w ,
        (r J)_\lambda\rangle\lambda^{-1}
        + \dlambda\lambda^{-1} \langle  w , r\partial_r
    J_\lambda \rangle \ . \notag
\end{align}
Therefore, to show the estimate \eqref{hardest_w} we will establish the following bounds:
\begin{align}
        \big|\langle \partial_t  A_\lambda w, (r J)_\lambda
        \rangle \dlambda\lambda\big|
        \ &\lesssim \ c^\frac{3}{4}_0\frac{\dlambda^2}{\lambda} +
        c_0^\frac{1}{4}\big(c_0\epsilon^2 +
        \sup_{0\leqslant s\leqslant t}
        \frac{\dot{\lambda}^4}{\lambda^7}(s) \big)\lambda^2
        \ , \label{hardest_red1}\\
        \big|\langle [A_\lambda , \partial_t] w , (r
        J)_\lambda\rangle\dlambda\lambda  \big| \ &\lesssim \
        \epsilon^2\frac{\dlambda^2}{\lambda} +
        \big(c_0^2\epsilon^2 + \epsilon \sup_{0\leqslant s\leqslant t}
        \frac{\dot{\lambda}^4}{\lambda^7}(s) \big)\lambda^2
        \ , \label{hardest_red2}\\
        \big| \langle  w , r\partial_r J_\lambda \rangle
        \dlambda^2\lambda \big| \ &\lesssim \
        \epsilon^2 \frac{\dlambda^2}{\lambda} +
        \big(c_0^2\epsilon^2 + \epsilon\sup_{0\leqslant s\leqslant t}
        \frac{\dot{\lambda}^4}{\lambda^7}(s) \big)\lambda^2
        \ . \label{hardest_red3}
\end{align}
We  prove these  three estimates in order.  By the
triangle inequality we have that:
\begin{equation}
        \big|\langle \partial_t A_\lambda w ,
        (r J)_\lambda\rangle \big| \
        \lesssim \  \big|\langle L A_\lambda w ,
        (r J)_\lambda\rangle \big|
        + \big|\langle \partial_r A_\lambda w ,
        (r J)_\lambda\rangle \big|
        \ . \label{the_sum}
\end{equation}
To estimate the left hand side of \eqref{hardest_red1} involving the
first term in the last sum we write:
\begin{align}
        &\big|\langle L A_\lambda w , (r J)_\lambda\rangle\dlambda\lambda
        \big|  \ , \notag\\
        \lesssim \
        &\dlambda\lambda
        \lp{(\lambda r)^\frac{\delta}{2}
        (1+r^\delta)^{-\frac{1}{2}} LA_\lambda w}{L^2(rdr)}\,
        \lp{(1+r^\delta)^\frac{1}{2}
        (r^{1-\frac{\delta}{2}}F^4)_\lambda}{L^2(rdr)} \ ,
        \notag\\
        \lesssim \ &\dlambda\lambda^{-\frac{1}{2}}\
        \lp{\lambda^{-\frac{1}{2}}(\lambda r)^\frac{\delta}{2}
        (1+r^\delta)^{-\frac{1}{2}}L A_\lambda w}{L^2(rdr)}\lambda \ ,
        \notag\\
        \lesssim \ &c_0^\frac{3}{4}\dlambda^2\lambda^{-1} + c_0^{-\frac{3}{4}}
        \lp{\lambda^{-\frac{1}{2}}(\lambda r)^\frac{\delta}{2}
        (1+r^\delta)^{-\frac{1}{2}} LA_\lambda w}{L^2(rdr)}^2\lambda^2 \ ,
        \notag\\
        \lesssim \ &c_0^\frac{3}{4}\dlambda^2\lambda^{-1} + c_0^\frac{1}{4}
        \big(c_0\epsilon^2 + \sup_{0\leqslant s\leqslant t}
        \dlambda^4\lambda^{-7} \big)\lambda^2 \ . \label{h_red_done1}
\end{align}
In the last line above we have used the assumption $\epsilon\leqslant c_0^2$.\\

To conclude the estimate \eqref{hardest_red1}, it remains to bound
the second term on the right hand side of \eqref{the_sum}. To do
this we see that a simple calculation involving the notation on
lines \eqref{F_rules} and \eqref{J_schematic} above, allows us to
write $(\partial_r +r^{-1}) (rJ)_\lambda = \lambda F^4_\lambda$. This
leads us to the estimate:
\begin{align}
        &\ \ \ \ \ \big|\langle \partial_r A_\lambda
        w , (r J)_\lambda\rangle
        \dlambda\lambda \big|
    \ =   \
        \big|\langle A_\lambda w ,F^4_\lambda \rangle
        \dlambda\lambda^2 \big| \ , \notag\\
        &\lesssim  \ \dlambda\lambda\, \lp{(\lambda
        r)^\frac{\delta}{2}(1+r^\delta)^{-\frac{1}{2}}
        r^{-1}A_\lambda w}{L^2(rdr)}\, \lp{(1+r^\delta)^\frac{1}{2}
        (r^{1-\frac{\delta}{2}}F^4)_\lambda }
        {L^2(rdr)} \ , \notag\\
        &\lesssim \ c_0^\frac{3}{4}\dlambda^2\lambda^{-1}
        + c_0^{-\frac{3}{4}}
        \lp{\lambda^{-\frac{1}{2}}(\lambda r)^\frac{\delta}{2}
        (1+r^\delta)^{-\frac{1}{2}} r^{-1}A_\lambda w}{L^2(rdr)}^2\lambda^2 \ ,
        \notag\\
        &\lesssim \ c_0^\frac{3}{4}\dlambda^2\lambda^{-1}
        + c_0^\frac{1}{4}
        \big(c_0\epsilon^2 + \sup_{0\leqslant s\leqslant t}
        \dlambda^4\lambda^{-7} \big)\lambda^2 \ . \label{h_red_done2}
\end{align}
This completes our proof of the estimate \eqref{hardest_red1}.\\

To finish the proof of \eqref{hardest} we need to establish the
estimates \eqref{hardest_red2}--\eqref{hardest_red3} above. As we shall see,
the proof of \eqref{hardest_red3}, with a minor exception, follows almost
verbatim from the estimates we will use for \eqref{hardest_red2}.
Therefore we
now concentrate on  \eqref{hardest_red2}. A simple computation shows that
we may write the commutator as $[
A_\lambda,\partial_t]=-\partial_t(A_\lambda)$, which from line
\eqref{A_schematic} is a multiplication operator given by a function
of the form $\dlambda F^9_\lambda$. Thus, we compute that:
\begin{align}
        &\ \ \ \ \ \big|\langle [A_\lambda , \partial_t] w , (r
        J)_\lambda\rangle\dlambda\lambda  \big|
        \ \ \lesssim  \ \
        \big|\langle  w , F^{12}_\lambda \rangle\dlambda^2
        \lambda  \big| \ , \notag\\
        &\lesssim \ \dlambda^2\lambda^{-\frac{1}{2}}
        \lp{\lambda^{-\frac{1}{2}}(\lambda r)^\frac{\delta}{2}
        (1+r)^{-\frac{\delta}{2}} r^{-2} w}{L^2(rdr)}\
        \lp{(1+r)^\frac{\delta}{2} (r^{2-\frac{\delta}{2}}
        F^{12})_\lambda}{L^2(rdr)}
        \ , \notag\\
        &\lesssim \ \dlambda^2\lambda^{-\frac{3}{2}}\cdot\big(
        c_0^2\epsilon^2 + \epsilon\sup_{0\leqslant s\leqslant t}
        \dot{\lambda}^4\lambda^{-7}(s)
        \big)^\frac{1}{2} \ , \notag\\
        &\lesssim \ \dlambda^4\lambda^{-5} + \big(
        c_0^2\epsilon^2 + \epsilon\sup_{0\leqslant s\leqslant t}
        \dot{\lambda}^4\lambda^{-7}(s)
        \big) \lambda^2\ , \notag\\
    &\lesssim \ \epsilon^2 \dlambda^2\lambda^{-1} + \big(
        c_0^2\epsilon^2 + \epsilon\sup_{0\leqslant s\leqslant t}
        \dot{\lambda}^4\lambda^{-7}(s)
    \big)\lambda^2\ , \notag
\end{align}\ret
where we used that $\dlambda \lambda^{-2}\lesssim \epsilon$
and \eqref{special_energy_w}.

The proof of the bound \eqref{hardest_red3} is
very similar to what was done above. To set things up in terms of
the previous steps, we simply use the notation on line
\eqref{J_schematic} and the Cauchy-Schwartz inequality to write:
\begin{multline}
    \big| \langle  w , r\partial_r J_\lambda \rangle
        \dlambda^2\lambda \big| \\
    \lesssim \
    \dlambda^2\lambda^{-\frac{1}{2}}
    \lp{\lambda^{-\frac{1}{2}}(\lambda r)^\frac{\delta}{2}
    (1+r)^{-\frac{\delta}{2}} r^{-2} w}{L^2(rdr)}\
    \lp{(1+r)^\frac{\delta}{2} (r^{2-\frac{\delta}{2}} F^{4})_\lambda}{L^2(rdr)}
    \ . \notag
\end{multline}
The steps are now identical to what was done in the previous
due to the bound:
\begin{equation}
    \lp{(1+r)^\frac{\delta}{2} (r^{2-\frac{\delta}{2}}
    F^{4})_\lambda}{L^2(rdr)}
     \ \lesssim \ \lambda^{-1} \ . \label{k_restriction}
\end{equation}
This concludes our proof of the estimate \eqref{hardest}.
\end{proof}\ret


\begin{proof}[Proof of the identity \eqref{easiest} and the bounds \eqref{eta_bounds}]
This follows from the orbital stability bound
\eqref{energy_size} and the Cauchy-Schwartz inequality.
Specifically, a short calculation shows that:
\begin{equation}
    \langle u , \ddot{J}_\lambda \rangle\lambda^3
    \ = \ \langle u , (r\partial_r J)_\lambda
    \rangle\ddot{\lambda}\lambda^2
    - \langle u , (r\partial_r J)_\lambda
    \rangle\dlambda^2\lambda +
    \langle u , \big((r\partial_r)^2 J\big)_\lambda
    \rangle\dlambda^2\lambda \ . \notag
\end{equation}
We leave the details of application of the Cauchy-Schwartz and the
estimate \eqref{energy_size} to the reader.
\end{proof}\ret


\begin{proof}[Proof of the estimate \eqref{medium}]
By invoking the explicit formula \eqref{N_breakdown} and
the decomposition \eqref{u_decomp} it
suffices to show the two estimates:
\begin{align}
        \langle\frac{(w_0)^2}{r^2} , J_\lambda \rangle\lambda^3 \ &\lesssim \
        \epsilon^2\frac{\dlambda^2}{\lambda} \ , \label{medium_red1}\\
        \langle\frac{w^2}{r^2} , J_\lambda\rangle\lambda^3 \ &\lesssim \
        \big(c_0^2\epsilon^2 + \epsilon
        \sup_{0\leqslant s\leqslant t} \frac{\dlambda^4}{\lambda^7}(s)\big)
        \lambda^2 \ .  \label{medium_red2}
\end{align}
The proof of \eqref{medium_red1} is a simple and direct calculation
involving the definition \eqref{w0_def} and the estimate
\eqref{orbit_bound}. We leave the details to the reader.\\

The proof of the second estimate \eqref{medium_red2} follows almost immediately from \eqref{special_energy_w}. To
see this, we compute that:
\begin{align}
        \langle\frac{w^2}{r^2} , J_\lambda\rangle\lambda^3 \ &\lesssim \
        \lambda^2\ \lp{\lambda^{-\frac{1}{2}}(\lambda r)^\frac{\delta}{2}
        (1+r)^{-\frac{\delta}{2}} r^{-2} w}{L^2(rdr)}^2\
        \lp{(1+r)^{\frac{\delta}{2}} (r^{2-\frac{\delta}{2}}F^4)_\lambda}{L^\infty}
        \ , \notag\\
        &\lesssim \ \big(c_0^2\epsilon^2 + \epsilon
        \sup_{0\leqslant s\leqslant t} \dlambda^4\lambda^{-7}(s)\big)
        \lambda^2 \ . \notag
\end{align}
This completes our proof of the estimates \eqref{medium}.
\end{proof}\ret

Having now completed our proof of the estimate
\eqref{hardest}--\eqref{medium}, our last task in the section is to
establish the structure estimate \eqref{main_bounds} for the
function $\mathcal{E}$ defined on line \eqref{cal_E_def}. To do this
it clearly suffices to add together the following set of estimates
for the individual terms on the right hand side of
\eqref{cal_E_def}:\\

\begin{lem}[Estimates for \eqref{main_bounds}]\label{cal_E_lem}
Assuming the bootstrapping estimate \eqref{ldd_structure_boot} and
the results of Proposition \ref{tech_prop}, one has the following
estimates where the implicit constant depends on line
\eqref{special_energy}:
\begin{align}
        \big| \langle w , \dot{J}_\lambda \rangle \dlambda\big|
        \ &\lesssim \ c_0\epsilon^2 + c_0^\frac{1}{2}
        \sup_{0\leqslant s\leqslant t}
        \frac{\dot{\lambda}^4}{\lambda^7}(s) \ , \label{cal_E_est1}\\
        \big| \langle w,\ddot{J}_\lambda \rangle\lambda \big|
        \ &\lesssim \ c_0\epsilon^2 + c_0^\frac{1}{2}
        \sup_{0\leqslant s\leqslant t}
        \frac{\dot{\lambda}^4}{\lambda^7}(s) \ , \label{cal_E_est2}\\
        \big| \langle w_0,\big(\ddot{\lambda} -2\frac{\dlambda^2}{\lambda}
        \big)(r\partial_r J)_\lambda \rangle\big|
        \ &\lesssim \
        c_0\epsilon^2 + c_0^\frac{1}{2}
        \sup_{0\leqslant s\leqslant t}
        \frac{\dot{\lambda}^4}{\lambda^7}(s)
        \ , \label{cal_E_est3}\\
        \big|\langle\frac{w\cdot(2w_0+w)}{r^2} ,\sin(2I_\lambda)
        J_\lambda \rangle\lambda\big|
        \ &\lesssim \ c_0\epsilon^2 + c_0^\frac{1}{2}
        \sup_{0\leqslant s\leqslant t}
        \frac{\dot{\lambda}^4}{\lambda^7}(s) \ , \label{cal_E_est4}\\
        \big| \langle\td{\mathcal{N}}(u) , J_\lambda
        \rangle\lambda \big|
        \ &\lesssim \ c_0\epsilon^2 + c_0^\frac{1}{2}
        \sup_{0\leqslant s\leqslant t}
        \frac{\dot{\lambda}^4}{\lambda^7}(s) \ . \label{cal_E_est5}
\end{align}
\end{lem}\ret

The proof of the estimates \eqref{cal_E_est1}--\eqref{cal_E_est5} is
similar to the proof of the estimates in Lemma \ref{accl_lem} above.
We will always follow the three-step strategy: 1) Distribute correct
powers of $r$ and $(1+r)$ inside the inner product. 2) Apply the
Cauchy-Schwartz inequality. 3) Refer to the estimates
\eqref{energy_size}, \eqref{orbit_bound},
\eqref{special_energy}--\eqref{special_energy_w}, and
\eqref{ldd_structure}. We will be a bit more terse here than before,
and leave some of the details to the reader. Each proof will be
written out under an individual heading.\\

\begin{proof}[Proof of the estimates \eqref{cal_E_est1}]
We start with the estimate \eqref{cal_E_est1}. Using our abstract
notation from line \eqref{J_schematic} we have that:
\begin{align}
        &\big| \langle w , \dot{J}_\lambda \rangle
        \dlambda\big| \ , \notag\\
        \lesssim \
        &\dlambda^2\lambda^{-\frac{5}{2}}
        \lp{\lambda^{-\frac{1}{2}}(\lambda r)^\frac{\delta}{2}
        (1+r)^{-\frac{\delta}{2}}
        r^{-2}w}{L^2(rdr)}\
        \lp{(1+r)^\frac{\delta}{2}(r^{2-\frac{\delta}{2}}
        F^4)_\lambda}{L^2(rdr)}
        \ , \notag\\
        \lesssim \ &\dlambda^2\lambda^{-\frac{7}{2}}\cdot\big(
        c_0^2\epsilon^2 + \epsilon\sup_{0\leqslant s\leqslant t}
        \dot{\lambda}^4\lambda^{-7}(s) \big)^\frac{1}{2}
        \ , \notag\\
        \lesssim \ &c_0\epsilon^2 + c_0^\frac{1}{2}
        \sup_{0\leqslant s\leqslant t}
        \dot{\lambda}^4\lambda^{-7}(s) \ . \notag
\end{align}
We remark here that this and similar estimates (below) are the
source of our restriction $4\leqslant k$ on the homotopy class in
Theorem \ref{main_thm}. Notice that one cannot arrive at the desired bound
by simply applying the orbital stability estimate \eqref{energy_size}.
It is crucial that we use \eqref{special_energy_w} here, and this
causes more weights to be placed on $r\partial_r J_\lambda$. It is
likely that one can lower
the value of $k$ in these arguments through a more careful analysis. We will not
pursue this here.
\end{proof}\ret


\begin{proof}[Proof of estimate \eqref{cal_E_est2}]
This is similar to the proof of \eqref{cal_E_est1}. Notice that from the estimate
\eqref{ldd_structure} we have the following bound:
\begin{equation}
    |\ddot{\lambda}| \ \lesssim \ \epsilon^2\lambda^2 +
    \sup_{0\leqslant s\leqslant t} \dlambda^2\lambda^{-1}(s) \ . \notag
\end{equation}
Therefore, using the notation form line \eqref{J_schematic} we have
the chain of inequalities:
\begin{align}
        \big| \langle w,\ddot{J}_\lambda \rangle\lambda \big| \ &\lesssim \
        \epsilon^2\, \big| \langle w,F^4_\lambda \rangle\lambda^2 \big| +
        \sup_{0\leqslant s\leqslant t} \dlambda^2\lambda^{-1}(s)\cdot
        \big| \langle w,F^4_\lambda \rangle \big| \ , \notag\\
        &\lesssim \ \epsilon^3 +
        \sup_{0\leqslant s\leqslant t} \dlambda^2\lambda^{-\frac{7}{2}}(s)\cdot
        \big(
        c_0^2\epsilon^2 + \epsilon\sup_{0\leqslant s\leqslant t}
        \dot{\lambda}^4\lambda^{-7}(s)
        \big)^\frac{1}{2} \ , \notag\\
        &\lesssim \ c_0\epsilon^2 + c_0^\frac{1}{2}
        \sup_{0\leqslant s\leqslant t}
        \dot{\lambda}^4\lambda^{-7}(s) \ . \notag
\end{align}
Notice that we have again used the condition $\epsilon\leqslant
c_0^2$ on this last line.
\end{proof}\ret


\begin{proof}[Proof of estimate \eqref{cal_E_est3}]
This will follow by a direct application of the estimate
\eqref{ldd_structure} and the definition \eqref{w0_def}. Notice that
\eqref{ldd_structure} and \eqref{orbit_bound} taken together imply
that:
\begin{equation}
    \big|\ddot{\lambda} -2\dlambda^2\lambda^{-1}\big|
    \ \lesssim \ \epsilon^2\lambda^2 +
    c_0^\frac{1}{2}\sup_{0\leqslant s\leqslant t} \dlambda^2
    \lambda^{-1}(s) \ . \notag
\end{equation}
Therefore, a simple computation again using \eqref{orbit_bound} shows that:
\begin{align}
    \big| \langle w_0,\big(\ddot{\lambda} -2\dlambda^2\lambda^{-1}
    \big)(r\partial_r J)_\lambda \rangle\big| \ &\lesssim \
    \dlambda^2\lambda^{-6}\cdot\big(\epsilon^2\lambda^2 +
    c_0^\frac{1}{2}\sup_{0\leqslant s\leqslant t} \dlambda^2
    \lambda^{-1}(s)\big) \ , \notag\\
    &\lesssim \ \epsilon^4 + c_0^\frac{1}{2}
    \sup_{0\leqslant s\leqslant t} \dlambda^4
    \lambda^{-7}(s) \ , \notag
\end{align}
which is enough to imply \eqref{cal_E_est3}  since
$\epsilon\leqslant c_0^2$.
\end{proof}\ret


\begin{proof}[Proof of estimate \eqref{cal_E_est4}]
The left hand side of this estimate can be  bounded by the
inequality:
\begin{multline}
    \big|\langle \frac{w^2 + |w_0|\cdot|w|}{r^2} , F^4_\lambda
    \rangle\lambda \big| \ \lesssim \
    \lp{\lambda^{-\frac{1}{2}}(\lambda r)^\frac{\delta}{2}
    (1+r)^{-\frac{\delta}{2}}
    r^{-2}w}{L^2(rdr)}^2 \\
    + \ \lp{\lambda^{-\frac{1}{2}}(\lambda r)^\frac{\delta}{2}
    (1+r)^{-\frac{\delta}{2}}
    r^{-2}w_0}{L^2(rdr)}\cdot
    \lp{\lambda^{-\frac{1}{2}}(\lambda r)^\frac{\delta}{2}
    (1+r)^{-\frac{\delta}{2}}
    r^{-2}w}{L^2(rdr)} \ . \notag
\end{multline}
The estimate \eqref{cal_E_est4} now follows directly from
\eqref{special_energy_w} applied to the terms involving $w$ in this
last line above, and the following bound which is a consequence of
the explicit identity \eqref{w0_def} (or the notation on line
\eqref{w0_schematic}):
\begin{equation}
        \lp{\lambda^{-\frac{1}{2}}(\lambda r)^\frac{\delta}{2}
        (1+r)^{-\frac{\delta}{2}}
        r^{-2}w_0}{L^2(rdr)} \ \lesssim \ \dlambda^2\lambda^{-\frac{7}{2}}
        \ . \label{w0_bound_no_epsl}
\end{equation}
\end{proof}\ret


\begin{proof}[Proof of the inequality \eqref{cal_E_est5}]
To do this, we first note
that by the inequality \eqref{td_N_bound} and the orbital stability estimate
\eqref{energy_size}, used in conjunction with the Poincar\'e type estimate
\eqref{energy_Poincare} as well as the notation from line \eqref{J_schematic},
we have that:
\begin{equation}
    \big| \langle\td{\mathcal{N}}(u) , J_\lambda \rangle\lambda \big|
     \ \lesssim \ \epsilon\ \lp{\lambda^{-\frac{1}{2}}
     (\lambda r)^\frac{\delta}{2}(1+r)^{-\frac{\delta}{2}}
    r^{-2}u}{L^2(rdr)}^2 \ . \notag
\end{equation}
By adding together the estimate \eqref{w0_bound_no_epsl} and the
corresponding bound \eqref{special_energy_w} for $w$ in the
decomposition  $u=w_0+w$, we have the single estimate for $u$:
\begin{equation}
    \lp{\lambda^{-\frac{1}{2}}(\lambda r)^\frac{\delta}{2}
    (1+r)^{-\frac{\delta}{2}}
    r^{-2}u}{L^2(rdr)}^2 \ \lesssim \ c_0^2\epsilon^2
    + \sup_{0\leqslant s \leqslant t }\dlambda^4
    \lambda^{-7}(s) \ . \notag
\end{equation}
Substituting this into the right hand side of the previous line we obtain
the desired bound \eqref{cal_E_est5}. This completes our proof of
Lemma \ref{cal_E_lem}.
\end{proof}

\ret\ret

\section{Space-Time Bounds and the Proof of the Main Estimate}\label{mor_sect}

In this final section of the paper we prove our main technical
estimate \eqref{special_energy}. The crucial role in this will be played
by the remarkable factorization property of the linearized
Hamiltonian \eqref{the_fact}, which allows us to introduce the
``conjugate'' Hamiltonian \eqref{conj_ham}. This new Hamiltonian $\td{H}_\lambda$
possesses
the striking properties \eqref{positivity}--\eqref{pos_energy} which
are ultimately responsible for very strong estimates, proved dynamically and by means of
simple yet quite precise physical space methods,  for the corresponding
Cauchy problem $(\pa_t^2+\td{H}_\lambda)w=F$ .
The key is the physical-space repulsive properties of
the operator \eqref{conj_ham} which lead to the desired estimates
\emph{independent of how violently the scaling parameter $\lambda$ grows}, so long as this
growth is monotonic. This stands in stark contrast to the usual
procedure in asymptotic stability analysis, which attempts to
estimate the linearized operator through \emph{non-dynamic} spectral
analysis (see e.g. \cite{BP_NLS}, \cite{KS_SW}). Such a procedure is not as natural
in the present context, which represents a truly non-linear
situation   not directly amenable to the standard perturbative
techniques. From this point of view, the analysis we present here is close in spirit to the work of
Merle-Raphael on the blow-up for the critical NLS \cite{MR}.\\

We again remind the reader that it is the precise form of
the non-linear equation \eqref{EL_eqs_red},
embodied by the first order Bogomol'nyi
equation \eqref{B_eq}, that is the indispensable structure. \\

The first thing we will need here for the proof of
\eqref{special_energy} is a space-time estimate for general
solutions to the conjugated linearized equation \eqref{conj_ham}.
For us this will take the form of a weighted $L^2$ inequality
involving integration over both space and time variables. These are
commonly referred to as Morawetz estimates, and they have a rich
history in both linear and nonlinear analysis of the dispersive
properties of wave equations. The estimate we use here is the
based on the following energy, defined for sufficiently smooth
and well decaying functions $\psi$ on $\mathbb{R}^+$:
\begin{multline}
    \mathbb{E}_\delta[\psi](t_0,t_1) \ = \
    \sup_{t_0\leqslant s\leqslant t_1}\ \int_{\mathbb{R}^+}\
    \lambda^{-1}\frac{(\lambda r)^\delta}{1+r^\delta}
    \left[(L\psi)^2 + \frac{\psi^2}{r^2} \right](s)\ rdr  \\
    + \ \ \int_{t_0}^{t_1} \int_{\mathbb{R}^+}\ \lambda^{-1}\left[
    \frac{(\lambda r)^\delta}{(1+r^\delta)^2\, r}(L\psi)^2
    \ + \ \frac{(\lambda r)^\delta}{1+r^\delta}
    \frac{\psi^2}{r^3} \right](s)\ rdr\, ds \ , \label{main_mor_eng}
\end{multline}
where  $0 < \delta\ll 1$ is a small fixed constant which will measure
a loss in certain time integrations which appear in the sequel.
Here we have set $L=\partial_t + \partial_r$. The main estimate
we will use is contained in the following:\\

\begin{prop}[Morawetz estimate for $\td{H}_\lambda$]\label{Mor_prop}
Consider the time dependent Hamiltonian $\td{H}_\lambda$. Let
$\psi$ be a smooth function on $[t_0,t_1]\times(0,\infty)$,
satisfying the following uniform bounds:
\begin{align}
        |\psi| \ &\leqslant \ C_{\psi(t)}\, r \ ,
        &|\partial_t\psi| + |\partial_r \psi|
        \ &\leqslant \ C_{\psi(t)} \ ,
        &0\ \leqslant \ r \ \leqslant \ 1 \ .
        \label{vanishing_bounds}
\end{align}
while decaying sufficiently rapidly at $r=\infty$. Furthermore,
suppose that:
\begin{equation}
    \partial_t^2\psi + \td{H}_\lambda\psi \ = \
    \partial_t G + H \ . \label{varphi_structure}
\end{equation}
Then if one has the
pointwise inequalities $0\leqslant \dlambda $ and $\dlambda
\lambda^{-2} \lesssim \epsilon$ for all times $t_0\leqslant s
\leqslant t_1$, one also has the following estimate:
\begin{multline}
    \mathbb{E}_\delta[\psi](t_0,t_1) \ \lesssim \ \delta^{-1}\Big[\
    \int_{t_0}^{t_1} \int_{\mathbb{R}^+}\ \lambda^{-1}(\lambda r)^\delta
    \left[ (\partial_rG)^2 + \epsilon^2 (\lambda G)^2
    + H^2 \right](s)\ r^2dr\, ds \\
    + \ \sup_{t_0\leqslant s\leqslant t_1} \int_{\mathbb{R}^+}\
    \lambda^{-1}\ \frac{(\lambda r)^\delta}{1+r^\delta} G^2(s)\ rdr
    \ + \ \mathbb{E}_\delta[\psi](t_0,t_0) \ \Big]
    \ , \label{main_mor_est}
\end{multline}
which holds with an implicit constant independent of $\lambda$ and $\delta$.
\end{prop}\ret

\begin{rem}
The  constant $0 < \delta\ll 1$ will signify a small loss in time
when we attempt to apply \eqref{main_mor_est} in the proof of
\eqref{special_energy}. This is ultimately why we are restricted to
the time interval $[0,\epsilon^{-4}]$ in the statement of
\eqref{tech_prop} and hence in Proposition \eqref{main_prop}. We
also remark here that this small loss in time can in fact be avoided
through a somewhat more careful analysis involving the precise form
of the equation for $\lambda$ given on line
\eqref{refined_main_ODE}.
\end{rem}\ret

\begin{proof}[Proof of the estimate \eqref{main_mor_est}]
Notice that all of the weights in the inequality are time translation
invariant. Therefore, we may normalize the discussion to $t_0=0$.
We begin by conjugating the equation
\eqref{varphi_structure} by $r^\frac{1}{2}$.
Therefore, we denote the new variable:
\begin{equation}
        \td{\psi} \ = \
        r^\frac{1}{2} \psi \ , \notag
\end{equation}
We note here that the decay $\td{\psi}$  at the
origin ($\sim r^\frac{3}{2}$) will be sufficient to perform the
integration by parts to follow. We also observe that $\partial_t\td{\psi}$ and
$\partial_r\td{\psi}$ may be assumed to be bounded at $r=0$.\\

Next, recall that the original (super-symmetric conjugate)
Hamiltonian has the form:
\begin{equation}
        \td{H}_\lambda \ = \ -\partial_r^2 - \frac{1}{r}\partial_r
        + V_\lambda(r) \ . \notag
\end{equation}
We define the one dimensional Hamiltonian:
\begin{equation}
        \mathcal{H}_\lambda \ = \ -\partial_r^2
        - \frac{1/4}{r^2} + V_\lambda(r) \ . \notag
\end{equation}
Then a quick computation shows that equation
\eqref{varphi_structure} becomes:
\begin{equation}
        \partial_t^2 \td{\psi} + \mathcal{H}_\lambda
        \td{\psi} \ = \
        r^\frac{1}{2} (\partial_t G + H ) \ . \ \label{1d_ham}
\end{equation}
The multiplier we use is the following:
\begin{equation}
        X \ = \ \lambda^{-1+\delta} \frac{r^\delta}{1 + r^\delta}L \ = \
        \lambda^{-1+\delta}\left(1- \frac{1}{1+r^\delta}\right)L \ .
        \notag
\end{equation}
Multiplying the equation \eqref{1d_ham} by the quantity
$X\td{\psi}$ and integrating the resulting expression over the
interval $[0,t]\times(0,\infty)$ we have the identity:
\begin{multline}
        \frac{1}{2}\int_0^t \int_{\mathbb{R}^+}\ \lambda^{-1}
        \frac{(\lambda r)^\delta}{1+r^\delta} \Big[\
        \bL (L\td{\psi})^2 +
        (V_\lambda -\frac{1}{4}r^{-2}) L(\td{\psi})^2\
        \Big] \ dr\, dt  \\
        = \ \int_0^t \int_{\mathbb{R}^+}\ \lambda^{-1}
        \frac{(\lambda r)^\delta}{1+r^\delta}
        \big(\, \bL(G) + \partial_r(G) + H\, \big)
        \cdot L(\td{\psi}) \ r^\frac{1}{2}dr\, dt \ .
        \label{int_iden}
\end{multline}
Here $\bL = \partial_t - \partial_r$ is the incoming null
derivative. We integrate by parts on the left hand side of this last
expression, using the following lower bounds for terms involving the
potential:
\begin{align}
        C\frac{k^2}{r^2} \ \geqslant \ (V_\lambda -\frac{1}{4}r^{-2}) \ &\geqslant \ c\,
        \frac{k^2}{r^2} \ , \notag\\
        -\, L\left[ \lambda^{-1}
        \frac{(\lambda r)^\delta}{1+r^\delta}
        (V_\lambda -\frac{1}{4}r^{-2})\right] \ &\geqslant \ c\,
        \lambda^{-1}\frac{(\lambda
        r)^\delta}{1+r^\delta}\frac{k^2}{r^3} \ , \notag
\end{align}
which follow from \eqref{positivity}-\eqref{pos_energy}, the condition $4\leqslant k$, and
the positivity of $\dlambda$.
Applying  a couple of Cauchy-Schwartz inequalities to the last
two terms on the right hand side of \eqref{int_iden}, and using the positivity
condition $\dlambda \geqslant 0$, we then we arrive at the bound:
\begin{multline}
        \int_{\mathbb{R}^+}\
        \lambda^{-1}\frac{(\lambda r)^\delta}{1+r^\delta}
        \left[(L\td{\psi})^2 + k^2 \frac{\td{\psi}^2}{r^2}
        \right](t)\ dr \
        + \ \delta \int_0^t \int_{\mathbb{R}^+}\ \lambda^{-1}
        \frac{(\lambda r)^\delta}{(1+r^\delta)^2\, r}
        (L\td{\psi})^2
         \ dr\, ds  \\
         + \ k^2 \int_0^t \int_{\mathbb{R}^+}\
         \lambda^{-1} \frac{(\lambda r)^\delta}
         {1+r^\delta}\frac{(\td{\psi})^2}{r^3}
         \ dr\, ds  \ , \\
         \lesssim \
         \left( \int_0^t \int_{\mathbb{R}^+}\ \lambda^{-1}
         (\lambda r)^\delta
         \big((\partial_r G)^2 + H^2\big)\
         r^2dr\, ds\right)^\frac{1}{2}\cdot
         \left( \int_0^t \int_{\mathbb{R}^+}\ \lambda^{-1}
         \frac{(\lambda r)^\delta}{(1+r^\delta)^2\, r}
         (L\td{\psi})^2
         \ dr\, ds \right)^\frac{1}{2}\\
         + \ \int_{\mathbb{R}^+}\
        \lambda^{-1}\frac{(\lambda r)^\delta}{1+r^\delta}
        \left[(L\td{\psi})^2 + k^2 \frac{\td{\psi}^2}{r^2}
        \right](0)\ dr
         \ + \ \int_0^t \int_{\mathbb{R}^+}\ \lambda^{-1}
        \frac{(\lambda r)^\delta}{1+r^\delta}\ \bL(G) \cdot
        L(\td{\psi})\ r^\frac{1}{2}dr\, ds
         \ . \ \label{big_guy}
\end{multline}
It remains to deal with the last integral on the right hand side of
the above expression. To do this, we integrate by parts with respect to
the incoming derivative $\bL$. Employing the pointwise bound:
\begin{align}
        \left| \bL \left(
        \lambda^{-1}
        \frac{(\lambda r)^\delta}{1+r^\delta} r^\frac{1}{2}
        \right)\right| \ &\lesssim \
        |\frac{\dlambda}{\lambda^2}|\cdot
        \frac{(\lambda r)^\delta}{1+r^\delta}
        r^\frac{1}{2} \ + \
        \lambda^{-1}
        \frac{(\lambda r)^\delta}{1+r^\delta} r^{-\frac{1}{2}} \ ,
        \notag\\
        &\lesssim \ \epsilon \frac{(\lambda r)^\delta}{1+r^\delta}
        r^\frac{1}{2} \ + \ \lambda^{-1}
        \frac{(\lambda r)^\delta}{(1+r^\delta) r}
        r^{\frac{1}{2}} \ , \notag
\end{align}
and using the equation \eqref{1d_ham} together with the upper bound
$|V_\lambda|\lesssim r^{-2}$ we have that:
\begin{multline}
        \int_0^t \int_{\mathbb{R}^+}\ \lambda^{-1}
        \frac{(\lambda r)^\delta}{1+r^\delta}\ \bL(G) \cdot
        L(\td{\psi})\ r^\frac{1}{2}dr\, ds  \ , \\
        \lesssim \
        \left(\int_0^t \int_{\mathbb{R}^+}\ \lambda^{-1}
        (\lambda r)^\delta
        \left[ \frac{G^2}{r^2} + \epsilon^2
        (\lambda G)^2 + H^2\right](s)\ r^2dr\, ds
        \right)^\frac{1}{2}\\
        \cdot
        \left(\int_0^t \int_{\mathbb{R}^+}\ \lambda^{-1}\left[
         \frac{(\lambda r)^\delta}{(1+r^\delta)^2\, r}
         (L\td{\psi})^2
        + \frac{(\lambda r)^\delta}
         {1+r^\delta}\frac{(\td{\psi})^2}{r^3}
         \right]\ dr\, ds\right)^\frac{1}{2}\\
         - \ \ \frac{1}{2} \int_0^t \int_{\mathbb{R}^+}\ \lambda^{-1}
        \frac{(\lambda r)^\delta}{1+r^\delta}\
        \partial_t (G^2) \ rdr\, ds-\int_0^t \int_{\mathbb{R}^+}\ \lambda^{-1}
        \frac{(\lambda r)^\delta}{1+r^\delta}\
        G H \ rdr\, ds\\\
        + \ \ \left(\sup_{0\leqslant s \leqslant t}
        \int_{\mathbb{R}^+}\ \lambda^{-1}
        \frac{(\lambda r)^\delta}{1+r^\delta} G^2(s) \ rdr
        \right)^\frac{1}{2}
        \cdot\left(\sup_{0\leqslant s \leqslant t}
        \int_{\mathbb{R}^+}\
        \lambda^{-1}\frac{(\lambda r)^\delta}{1+r^\delta}
        (L\td{\psi})^2(s)\ dr\right)^\frac{1}{2} \ .
        \label{next_bigbound}
\end{multline}
Integrating by parts one more time in the term involving $\partial_t(G)^2$
above, again using the fact that $\dlambda \geqslant 0$, and using also
the following fixed time Poincar\'e type estimate:
\begin{equation}
        \int_{\mathbb{R}^+}\ \lambda^{-1}
    {(\lambda r)^\delta}G^2 \
        dr \ \lesssim \
        \int_{\mathbb{R}^+}\ \lambda^{-1}
        (\lambda r)^\delta (\partial_r G)^2 \
        r^2dr \ , \label{poincare}
\end{equation}
we add together the estimates \eqref{big_guy}--\eqref{poincare} and
take the sup over different times to achieve the bound:
\begin{multline}
        \td{\mathbb{E}}_\delta[\td{\psi}](0,t) \ \lesssim \
        \delta^{-1}\Big[\
        \int_0^t \int_{\mathbb{R}^+}\ \lambda^{-1}(\lambda r)^\delta
        \left[ (\partial_rG)^2 + \epsilon^2 (\lambda G)^2 + H^2 \right](s)\ r^2dr\, ds \\
        + \ \sup_{0\leqslant s\leqslant t} \int_{\mathbb{R}^+}\
        \lambda^{-1}\ \frac{(\lambda r)^\delta}{1+r^\delta} G^2(s)\ rdr
        \ + \ \td{\mathbb{E}}_\delta[\td{\psi}](0,0) \ \Big] \ .
        \label{td_bound}
\end{multline}
where we define the one dimensional energy analogous to
\eqref{main_mor_eng}:
\begin{multline}
    \td{\mathbb{E}}_\delta[\td{\psi}](0,t) \ = \
    \sup_{0\leqslant s\leqslant t}\ \int_{\mathbb{R}^+}\
    \lambda^{-1}\frac{(\lambda r)^\delta}{1+r^\delta}
    \left[(L\td{\psi})^2 + k^2
    \frac{\td{\psi}^2}{r^2} \right](s)\ dr  \\
    + \ \ \int_0^t \int_{\mathbb{R}^+}\ \lambda^{-1}\left[
    \frac{(\lambda r)^\delta}{(1+r^\delta)^2\, r}(L\td{\psi})^2
    \ + \ k^2 \frac{(\lambda r)^\delta}{1+r^\delta}
    \frac{\td{\psi}^2}{r^3} \right](s)\ dr\, ds \ . \notag
\end{multline}
Finally, to complete the proof, we use the expansion:
\begin{equation}
        r^{-\frac{1}{2}}L\td{\psi} \ = \ L\psi +
        \frac{1}{2}r^{-1}\psi \ , \notag
\end{equation}
and the fact that $4\leqslant k$ to bound the energy
$\td{\mathbb{E}}_\delta$ from above and below by  $\mathbb{E}_\delta$. This
completes our proof of the estimate \eqref{main_mor_est}.
\end{proof}\ret\ret

We now turn to the proof of \eqref{special_energy}. The precise
statement of what we need to show is the following:\\

\begin{prop}[Energy estimates for the quantity $w$]\label{eng_w_prop}
Let $u=w_0 + w$ be the decomposition of $u$ given on line
\eqref{u_decomp}, where $u$ itself is part of the decomposition
\eqref{I_u_decomp} of the full field $\phi$. In particular $u$
solves the equation \eqref{lin_eq}. Suppose also that the initial
conditions for $\phi$ are given as to satisfy
\eqref{special_C_data}--\eqref{special_smallness}, and that $u$
obeys the extra decay estimate \eqref{extra_decay} (this will be proved
shortly). Furthermore,
assume that the main assumptions of Proposition \ref{tech_prop}
hold, in particular we have that \eqref{ldd_structure_boot} and
$\dlambda \geqslant 0$. Then the following estimate holds
for $t\in [0,\epsilon^{-4}]$:
\begin{equation}
    \mathbb{E}_\delta[A_\lambda w](0,t) \ \lesssim \
    \delta^{-1} \Big( c_0^2\epsilon^2 +
    \epsilon\sup_{0\leqslant s\leqslant t}\
    \frac{\dlambda^4}{\lambda^7}(s) \Big) \ . \label{conj_w_est}
\end{equation}
The implicit constant depends on $C$ from line
\eqref{ldd_structure_boot} but is independent of the size of $c_0$, $\epsilon$, or
$\delta$.\\

\noindent In particular, for a fixed $0<\delta\ll 1$ we have that the estimate
\eqref{special_energy} holds.
\end{prop}\ret

\begin{proof}[Proof of the estimate \eqref{conj_w_est}]
The first order of business is to reduce the proof to  simpler
bounds. In the sequel, we will only show that:
\begin{equation}
    \mathbb{E}_\delta[A_\lambda w](t_0,t_1) \ \lesssim \
    \delta^{-1} \Big( \mathbb{E}_\delta[A_\lambda w](t_0,t_0) + \epsilon^3 +
    \epsilon\sup_{t_0\leqslant s\leqslant t_1}\
    \frac{\dlambda^4}{\lambda^7}(s) \Big) \ , \label{new_conj_w_est}
\end{equation}
for all time intervals $[t_0,t_1]$ inside the regular interval $[0,T^{**})$, where
again $T^{**}$ is the blowup time,
provided that one \emph{also} has the  inequality:
\begin{equation}
    \int_{t_0}^{t_1} \ \frac{\dlambda^4}{\lambda^7}(s) \ ds \ \leqslant \ \epsilon \ . \label{basic_boot_bound}
\end{equation}
We claim that along with the bootstrapping assumption \eqref{ldd_structure_boot} and the
analysis done in Section \ref{ODE_an_sect}, this is enough to establish
\eqref{conj_w_est}.\\

To verify this  claim, first notice that if we are in the time interval $[0,T^{*}]$
where $T^{*}$ is defined as on line \eqref{better_Tstar_bound}, then we automatically
have \eqref{basic_boot_bound} on account of \eqref{first_int_bound}. Therefore, we
may work inside intervals of the form $[T^*,t]$, and we are only trying to establish:
\begin{equation}
    \mathbb{E}_\delta[A_\lambda w](T^*,t) \ \lesssim \
    \delta^{-1} \epsilon
    \frac{\dlambda^4}{\lambda^7}(t)  \ . \label{new_new_conj_w_est}
\end{equation}
Notice that we are using the monotonicity of $\dlambda^4\lambda^{-7}$ established in Section \ref{ODE_an_sect}.
We now claim that \eqref{new_new_conj_w_est} easily follows from \eqref{new_conj_w_est}
and the bootstrapping assumption \eqref{ldd_structure_boot}. Indeed, let $[t_{k-1},t_k]$
be any interval where \emph{equality} in \eqref{basic_boot_bound} holds. Then we have:
\begin{align}
    \int_{t_{k-1}}^{t_k}\ \frac{d}{ds} \ln(\dlambda^4 \lambda^{-7})\ ds \ &\gtrsim \
    \int_{t_{k-1}}^{t_k}\ \frac{\dlambda}{\lambda}\ ds \ \gtrsim\
    \epsilon^{-3} \int_{t_{k-1}}^{t_k}\ \frac{\dlambda^4}{\lambda^7}\ ds \ = \
    \epsilon^{-2} \ . \label{to_int}
\end{align}
Here we have used the bootstrapping bound \eqref{ldd_structure_boot} in the simple form
$|\ddot{\lambda}-2\dot\lambda^2\lambda^{-1}| \ll \dot\lambda^2\lambda^{-1}$, which holds as long
as we are in the region past $[0,T^*]$ (in particular, one has access to a lower bound
consistent with \eqref{better_Tstar_bound} which allows one to uniformize the RHS of
\eqref{ldd_structure_boot}). Notice that
we have also used the orbital stability bound \eqref{orbit_bound} several times in deriving
the inequalities. Integrating  the inequality \eqref{to_int} we see that:
\begin{equation}
    \frac{\dlambda^4}{\lambda^7}(t_{k-1}) \ \lesssim\ e^{-\epsilon^{-2}}
    \frac{\dlambda^4}{\lambda^7}(t_{k}) \ , \label{crucial_bound}
\end{equation}
on any time interval $[t_{k-1},t_{k}]$ past $[0,T^*]$ where \eqref{basic_boot_bound} also holds.
It is now a simple matter to derive \eqref{new_new_conj_w_est} from \eqref{new_conj_w_est}.
We first decompose the interval $[T^*,t]$ into a finite collection of $N$ subintervals
$[t_{k-1},t_{k}]$ where equality in \eqref{basic_boot_bound} holds.\footnote{On the last interval there may be
a strict inequality in \eqref{basic_boot_bound}, but this single interval may also be estimated
with \eqref{new_conj_w_est}, and the answer may then be directly added into the final bound.} On each of these intervals,
we have the estimate \eqref{new_conj_w_est}. By repeatedly using the bound \eqref{crucial_bound},
each of these estimates may be inductively expanded to yield:
\begin{align}
    \mathbb{E}_\delta[A_\lambda w](t_{k-1},t_k) \ &\lesssim \
    \delta^{-1}C \sum_{i=0}^k \left(\frac{C\delta^{-1}}{C_1}\right)^i
    \cdot \epsilon
    \frac{\dlambda^4}{\lambda^7}(t_k) \  \notag \\
    &\lesssim \  \frac{\delta^{-1}C}{C_1^{N-k}}
    \cdot \epsilon\
    \frac{\dlambda^4}{\lambda^7}(t) \ , \notag
\end{align}
where $C_1\sim e^{\epsilon^{-2}}$ is some incredibly large constant that beats the
(uniform) implicit constant $C\delta^{-1}$ appearing in the estimates \eqref{new_conj_w_est}.
Summing this last line over $0\leqslant k\leqslant N$, we have the bound \eqref{new_new_conj_w_est}.\\

We now prove \eqref{new_conj_w_est} under the additional assumption that \eqref{basic_boot_bound}
also holds.
We start by  providing the general setup, and then reduce the proof to a number of separate estimates to be
dealt with under their own bold-faced headings.
We first record the equation for $w$. Recall that
the purpose of the decomposition \eqref{u_decomp} is to eliminate
the main source term $A_\lambda(\ddot{I}_\lambda)$ on the right hand side of \eqref{lin_eq},
obtained after applying $A_\lambda$.
Therefore, we have that:
\begin{equation}
    A_\lambda\big[ \partial_t^2 w + H_\lambda w\big] \ = \ -A_\lambda(\ddot{w}_0)
    + A_\lambda\mathcal{N}(u) \ . \notag
\end{equation}
To put things in the form where the estimate \eqref{main_mor_est} can be used,
we commute the $A_\lambda$ operator with $\partial_t^2$ on the left hand side
of this last equation, which yields the expression:
\begin{align}
    &\partial_t^2 W + \td{H}_\lambda W \ , \notag\\
     = \  &-A_\lambda(\ddot{w_0}) + [\partial^2_t,A_\lambda]w +
     A_\lambda \mathcal{N}(u) \ , \notag\\
     = \ &\partial_t(A_\lambda)\partial_t(w_0) -
     \partial_t\big( A_\lambda\partial_t(w_0) \big)
     + 2\partial_t\big(\partial_t(A_\lambda)\cdot w\big)
     - \partial_t^2(A_\lambda)\cdot w
     + A_\lambda\mathcal{N}(u) \ , \notag\\
     = \ &M_1 + \partial_t M_2 + \partial_t R_1 + R_2 + R_3\ ,
     \label{detailed_lin}
\end{align}
where we have set $W=A_\lambda w$. The $M$ terms above constitute
the ``main source'' which feeds the quantity $W$ through the
wave-flow of the Hamiltonian $\td{H}_\lambda$. By contrast, the $R$
terms on line \eqref{detailed_lin} are for the most part ``errors'' which
will  be reabsorbed back onto the left hand side of the estimate
\eqref{conj_w_est}. This is where the limits on the time
interval and the decay estimate \eqref{extra_decay} will come in to
play. We now turn to the details of all of this.\\

To estimate the terms $M_1$ and $M_2$ via the general bound
\eqref{main_mor_est}, we will show the following four estimates:
\begin{align}
        \int_{t_0}^{t_1}\int_{\mathbb{R}^+}
        \ \lambda^{-1}(\lambda r )^\delta
        |M_1|^2(s) \ r^2dr\, ds \ &\lesssim \ \epsilon^3 \ + \
        \epsilon \sup_{t_0\leqslant s\leqslant t_1}
        \frac{\dlambda^4}{\lambda^7}(s)
        \ , \label{red_detailed_lin1}\\
        \int_{t_0}^{t_1}\int_{\mathbb{R}^+}\
        \lambda^{-1}(\lambda r )^\delta
        |\partial_r M_2|^2(s) \ r^2dr\, ds
        \ &\lesssim \ \epsilon^3 \ + \
        \epsilon \sup_{t_0\leqslant s\leqslant t_1}
        \frac{\dlambda^4}{\lambda^7}(s)
        \ , \label{red_detailed_lin2}\\
        \int_{t_0}^{t_1}\int_{\mathbb{R}^+}\
        \lambda^{-1}(\lambda r )^\delta
        |\lambda M_2|^2(s) \ r^2dr\, ds
        \ &\lesssim \ \epsilon^3 \ + \
        \epsilon \sup_{t_0\leqslant s\leqslant t_1}
        \frac{\dlambda^4}{\lambda^7}(s)
        \ , \label{red_detailed_lin3}\\
        \sup_{t_0\leqslant s\leqslant t_1}\
        \int_{\mathbb{R}^+}\
        \lambda^{-1}\frac{(\lambda r )^\delta}{1+r^\delta}
        | M_2|^2(s) \ rdr
        \ &\lesssim \ \epsilon^3 \ + \
        \epsilon \sup_{t_0\leqslant s\leqslant t_1} \frac{\dlambda^4}{\lambda^7}(s)
        \ . \label{red_detailed_lin4}
\end{align}
Recall we are assuming that $\epsilon\leqslant c_0^2$, so
these estimates will be enough to generate the right hand side of
\eqref{conj_w_est} for the $M$ terms.\\

We shall estimate the $R$ terms on line \eqref{detailed_lin} in a
nonlinear fashion. Specifically, we will bound them in terms of a
small constant times the energy on the left hand side of
\eqref{conj_w_est}, plus one term involving the nonlinearity $\mathcal{N}$
applied to $w_0$, which fits into the
pattern of the right hand side of
\eqref{red_detailed_lin1}--\eqref{red_detailed_lin4} above. What we
propose to show is the following:
\begin{align}
        \int_{t_0}^{t_1}\int_{\mathbb{R}^+}\
        \lambda^{-1}(\lambda r )^\delta
        |\partial_r R_1|^2(s) \ r^2dr\, ds
        \ &\lesssim \
        \epsilon^2 \mathbb{E}_\delta[W](t_0,t_1)
        \ , \label{R_red_detailed_lin1}\\
        \int_{t_0}^{t_1}\int_{\mathbb{R}^+}\
        \lambda^{-1}(\lambda r )^\delta
        |\lambda R_1|^2(s) \ r^2dr\, ds
        \ &\lesssim \
        \epsilon^2 \mathbb{E}_\delta[W](t_0,t_1)
        \ , \label{R_red_detailed_lin2}\\
        \sup_{t_0\leqslant s\leqslant t_1}\ \int_{\mathbb{R}^+}\
        \lambda^{-1}\frac{(\lambda r )^\delta}{1+r^\delta}
        | R_1|^2(s) \ rdr
        \ &\lesssim \
        \epsilon^2 \mathbb{E}_\delta[W](t_0,t_1)
        \ , \label{R_red_detailed_lin3}\\
        \int_{t_0}^{t_1}\int_{\mathbb{R}^+}\
        \lambda^{-1}(\lambda r )^\delta
        |R_2|^2(s) \ r^2dr\, ds
        \ &\lesssim \
        \epsilon^2 \mathbb{E}_\delta[W](t_0,t_1)
        \ , \label{R_red_detailed_lin4}\\
        \int_{t_0}^{t_1}\int_{\mathbb{R}^+}\
        \lambda^{-1}(\lambda r )^\delta
        |R_3|^2(s) \ r^2dr\, ds
        \ &\lesssim \
        \epsilon\sup_{t_0\leqslant s\leqslant t_1}
        \frac{\dlambda^4}{\lambda^7}(s) +
        t_1^\delta \epsilon^2 \mathbb{E}_\delta[W](t_0,t_1)
        \ . \label{R_red_detailed_lin5}
\end{align}
Using the conditions $\epsilon$ is sufficiently small,
and that $t_1\leqslant  \epsilon^{-4}$, all of the
estimates added together will imply the estimate \eqref{conj_w_est}
for the $R$ terms on line \eqref{detailed_lin} above. \\

We now turn to the details of the proofs of the estimates
\eqref{red_detailed_lin1}--\eqref{red_detailed_lin4} and
\eqref{R_red_detailed_lin1}--\eqref{R_red_detailed_lin5}. We will do
each of these separately and in order.\\

In what follows, we will consistently use the following
``identities'' which are in accordance with our
notation from Section \ref{not_sect}:
\begin{align}
    (\lambda r)^\delta F^m_\lambda\ &=\ F^{m-\delta}_\lambda \ ,
    &(1+r^\delta) F^m_\lambda\ &\leqslant \ F^{m-\delta}_\lambda\ , \notag
\end{align}
where the second ``inequality'' holds provided that $\lambda\geqslant 1$.\\


\subsection*{Proof of estimate \eqref{red_detailed_lin1}}
Here and throughout the sequel we will rely heavily on the abstract
function notation from lines \eqref{schematic_not} above.
Multiplying together estimates from lines \eqref{w0_schematic} and
\eqref{A_schematic}  we have that:
\begin{equation}
    |\partial_t(A_\lambda)\partial_t(w_0)| \ \lesssim \
    \Big(\dlambda^4\lambda^{-5} +
    |\ddot{\lambda}|\dlambda^2\lambda^{-4}\Big)
    F^9_\lambda \ . \notag
\end{equation}
To resolve the second term above which contains the expression
$\ddot{\lambda}$ we do the following. Notice that the orbital stability bound
$\dlambda\lambda^{-2}\lesssim \epsilon$  and the bootstrapping estimate
\eqref{ldd_structure_boot} give   the rough
\emph{pointwise} bound:
\begin{align}
    |\ddot{\lambda}| \ &\lesssim \
    \dlambda^2\lambda^{-1} + (\epsilon^2 + \dlambda^4
    \lambda^{-7})\lambda^2 \ , \notag\\
    &\lesssim \ \dlambda^2\lambda^{-1} + \epsilon^2
    \lambda^2 \ . \label{ddlambda_bound}
\end{align}
Therefore, substituting this estimate back into the previous line we
have that:
\begin{equation}
    |\partial_t(A_\lambda)\partial_t(w_0)| \ \lesssim \
    \Big(\dlambda^4\lambda^{-5} +
    \epsilon^2\dlambda^2\lambda^{-2}\Big)F^9_\lambda \ . \notag
\end{equation}
Using this last line we can now estimate:
\begin{align}
        \hbox{(L.H.S.)}\eqref{red_detailed_lin1} \ &\lesssim \
        \int_{t_0}^{t_1}\int_{\mathbb{R}^+}\ \lambda^{-1}
        \Big(\dlambda^4\lambda^{-5} +
        \epsilon^2\dlambda^2\lambda^{-2}\Big)^2 F^4_\lambda\ r^2dr \, ds \ , \notag\\
        &\lesssim \ \int_{t_0}^{t_1}\ \dlambda^8\lambda^{-14}(s)\ ds
        \ + \ \epsilon^4\int_{t_0}^{t_1}\ \dlambda^4\lambda^{-8}(s)\ ds \ , \notag\\
        &\lesssim \ \sup_{t_0\leqslant s \leqslant t_1}
        \dlambda^4\lambda^{-7}(s)\cdot
        \int_{t_0}^{t_1}\ \dlambda^4\lambda^{-7}(s)\ ds \ + \
        \epsilon^4\sup_{t_0\leqslant s \leqslant t_1}
        \dlambda^3\lambda^{-6}(s)\cdot
        \int_{t_0}^{t_1}\ \dlambda\lambda^{-2}(s)\ ds \ , \notag\\
        &\lesssim \ \epsilon \sup_{t_0\leqslant s \leqslant t_1}
        \dlambda^4\lambda^{-7}(s)
        \ + \ \epsilon^7 \ . \notag
\end{align}
To obtain the last line, we have used both the estimate
\eqref{basic_boot_bound}, the bound \eqref{orbit_bound}, and the assumption
that $\lambda\geqslant 1$.\\


\subsection*{Proof of estimate \eqref{red_detailed_lin2}}
This is very similar to the analysis above, with an addition of a small twist.
First of all, by combining lines \eqref{w0_schematic},
\eqref{A_schematic},  \eqref{F_rules}, and then \eqref{ddlambda_bound}
 we have the abstract notational bound:
\begin{align}
        \big| \partial_r\big( A_\lambda\partial_t(w_0) \big)
        \big|
        \ &\lesssim \ \Big(\dlambda^3\lambda^{-3} +
        |\ddot{\lambda}|\dlambda\lambda^{-2}\Big)F^4_\lambda \ , \notag\\
        &\lesssim \
        \Big( \dlambda^3\lambda^{-3} +
        \epsilon^2\dlambda \Big)F^4_\lambda \ . \notag
\end{align}
Substituting this into the left hand side of \eqref{red_detailed_lin2}
we have the chain of inequalities:
\begin{align}
    \hbox{(L.H.S.)}\eqref{red_detailed_lin2} \ &\lesssim \
    \int_{t_0}^{t_1}\int_{\mathbb{R}^+}\ \lambda^{-1}\Big(
    \dlambda^3\lambda^{-3} + \epsilon^2\dlambda
    \Big)^2 F^4_\lambda\ r^2dr \, ds \ , \notag\\
        &\lesssim \ \int_{t_0}^{t_1}\ \dlambda^6\lambda^{-10}(s)\ ds
    \ + \ \epsilon^4\int_{t_0}^{t_1}\ \dlambda^2\lambda^{-4}(s)\ ds \ , \notag\\
    &\lesssim \ \int_{t_0}^{t_1}\ \dlambda^6\lambda^{-10}(s)\ ds \ + \
    \epsilon^4\sup_{t_0\leqslant s \leqslant t_1}\dlambda\lambda^{-2}(s)\cdot
    \int_{t_0}^{t_1}\ \dlambda\lambda^{-2}(s)\ ds \ , \notag\\
    &\lesssim \ \int_{t_0}^{t_1}\ \dlambda^6\lambda^{-10}(s)\ ds
    \ + \ \epsilon^5 \ . \notag
\end{align}
We now need to do a little work, because the first term on the right
hand side of this last line above is not manifestly of the correct
form. To correct it, we first integrate by parts with respect to
time which yields the identity:
\begin{equation}
    \int_{t_0}^{t_1}\ \dlambda^6\lambda^{-10}\ ds   \ = \
    \dlambda^5\lambda^{-9}(t_1) - \dlambda^5\lambda^{-9}(t_0)
    + 10 \int_{t_0}^{t_1}\ \dlambda^6\lambda^{-10}\ ds -
    5\int_{t_0}^{t_1}\ \ddot{\lambda}\dlambda^4\lambda^{-9}\ ds \ . \notag
\end{equation}
Now, using the bootstrapping assumption \eqref{ldd_structure_boot}
as well as the estimate \eqref{orbit_bound},
this last expression leads to the following
nonlinear bound:
\begin{align}
        &\int_{t_0}^{t_1}\ \dlambda^6\lambda^{-10}\ ds
        \ , \label{6_10_lines}\\
        \lesssim \
        &c_0^\frac{1}{2} \int_{t_0}^{t_1}\
        \dlambda^6\lambda^{-10}(s)\ ds +
        \epsilon\, \sup_{t_0\leqslant s \leqslant t_1}\dlambda^4\lambda^{-7}(s)  +
        \int_{t_0}^{t_1}\ (\epsilon^2 + \sup_{t_0\leqslant \cdot\leqslant s}
        \dlambda^4\lambda^{-7})\cdot \dlambda^4
        \lambda^{-7}(s)\ ds \ , \notag\\
        \lesssim \ &c_0^\frac{1}{2} \int_{t_0}^{t_1}\
        \dlambda^6\lambda^{-10}(s)\ ds + \epsilon
        \sup_{t_0\leqslant s \leqslant t_1}
        \dlambda^4\lambda^{-7}(s)  +
    \epsilon^2 \int_{t_0}^{t_1}\ \dlambda^4\lambda^{-7}(s)\ ds \ , \notag\\
    \lesssim \ &c_0^\frac{1}{2} \int_{t_0}^{t_1}\
        \dlambda^6\lambda^{-10}(s)\ ds + \epsilon
        \sup_{t_0\leqslant s \leqslant t_1}
        \dlambda^4\lambda^{-7}(s) + \epsilon^3 \ . \notag
\end{align}
To go from the first to the second line above, we have used
\eqref{ldd_structure} and the monotonicity established in Section
\ref{ODE_sect} (this works as long as $T^*\leqslant t_0$, whereas in
the other case we may as well assume that $t_0=0$). Notice also that
in the last two lines above we made several uses of the
assumption \eqref{basic_boot_bound}.\\

\subsection*{Proof of estimate \eqref{red_detailed_lin3}}
This is virtually identical to the proof of
\eqref{red_detailed_lin2}. Another simple
calculation using lines \eqref{w0_schematic} and
\eqref{A_schematic}, and then \eqref{ddlambda_bound} gives us:
\begin{align}
        \big| \lambda\big( A_\lambda\partial_t(w_0) \big)
        \big|
        \ &\lesssim \ \Big(\dlambda^3\lambda^{-3} +
        |\ddot{\lambda}|\dlambda\lambda^{-2}\Big)
        F^3_\lambda \ , \notag\\
        &\lesssim \
        \Big( \dlambda^3\lambda^{-3} +
        \epsilon^2\dlambda \Big)F^3_\lambda \ . \notag
\end{align}
Substituting this into the left hand side of \eqref{red_detailed_lin3}
the proof follows verbatim from the calculations done in the
previous paragraph.\\

\subsection*{Proof of estimate \eqref{red_detailed_lin4}}
Once again using lines \eqref{w0_schematic} and \eqref{A_schematic},
and then \eqref{ddlambda_bound} we have that:
\begin{align}
        \big| \big( A_\lambda\partial_t(w_0) \big) \big|
        \ &\lesssim \ \Big(\dlambda^3\lambda^{-4} +
        |\ddot{\lambda}|\dlambda\lambda^{-3}\Big)F^3_\lambda \ , \notag\\
        &\lesssim \
        \Big( \dlambda^3\lambda^{-4} +
        \epsilon^2\dlambda\lambda^{-1} \Big)F^3_\lambda \ . \notag
\end{align}
Plugging this last line into the left hand side of
\eqref{red_detailed_lin4} and simply using the bound
\eqref{orbit_bound} we arrive at the chain of inequalities:
\begin{align}
        \hbox{(L.H.S.)}\eqref{red_detailed_lin4} \ &\lesssim \
        \sup_{t_0\leqslant s\leqslant t_1}\ \int_{\mathbb{R}^+}\
        \lambda^{-1}\Big( \dlambda^3\lambda^{-4} +
        \epsilon^2\dlambda\lambda^{-1} \Big)^2
        F^4_\lambda \ rdr \ , \notag\\
        &\lesssim \ \sup_{t_0\leqslant s\leqslant t_1}
        \dlambda^6\lambda^{-11}
        + \epsilon^4 \sup_{t_0\leqslant s\leqslant t_1}
        \dlambda^2\lambda^{-5} \ , \notag\\
        &\lesssim \ \epsilon^2 \sup_{t_0\leqslant s\leqslant t_1}
        \dlambda^4\lambda^{-7}
        + \epsilon^6  \ . \notag
\end{align}
This concludes our proof of the first list of estimates
\eqref{red_detailed_lin1}--\eqref{red_detailed_lin4} above.\\ \\


\noindent We now turn our attention to the proofs of the estimates
\eqref{R_red_detailed_lin1}--\eqref{R_red_detailed_lin5}.\\

\subsection*{Proof of \eqref{R_red_detailed_lin1}}
We first provide  a pointwise bound for the term $\partial_r R_1$
from line \eqref{detailed_lin} above. This involves a simple
application of the abstract notation \eqref{A_schematic}, the rules
from line \eqref{F_rules}, and the decomposition:
\begin{equation}
    \partial_r w \ = \ -\, W  + \frac{k}{r}\cos(I_\lambda)\cdot w
    \ . \notag
\end{equation}
Together, these give us the following estimate:
\begin{equation}
    \big|\partial_r\big(\partial_t(A_\lambda)\cdot w\big)\big|
    \lesssim \ \big( \big|\dlambda\lambda  w \big|
    + \big|\dlambda  W \big| \big)F^6_\lambda \ . \notag
\end{equation}
Plugging this identity in the left hand side of
\eqref{R_red_detailed_lin1} leads us to the estimates (also using
\eqref{orbit_bound}):
\begin{align}
        &\int_{t_0}^{t_1} \int_{\mathbb{R}^+}\ \lambda^{-1}
        (\lambda r)^\delta
        |\partial_r R_{1}|^2\
        r^2dr\, ds \  , \notag\\
        \lesssim \
        &\int_{t_0}^{t_1} \int_{\mathbb{R}^+} \ \frac{\dlambda^2}{\lambda^5}
        (\lambda r)^\delta\left[ \frac{w^2}{ r^5}
        + \frac{W^2}{ r^3}
        \right]\big((r^4+r^6) F^{12}\big)_\lambda\
        rdr\, dt \ , \notag\\
        \lesssim \ &\epsilon^2\ \int_{t_0}^{t_1} \int_{\mathbb{R}^+} \ \lambda^{-1}
        \frac{(\lambda r)^\delta}{(1+r)^\delta}
        \left[ \frac{w^2}{ r^5}
        + \frac{W^2}{ r^3}
        \right]\ rdr\, dt \ , \notag\\
        \lesssim \ &\epsilon^2\ \int_{t_0}^{t_1}
        \int_{\mathbb{R}^+} \ \lambda^{-1}
        \frac{(\lambda r)^\delta}{(1+r)^\delta}
        \frac{W^2}{ r^3} \ rdr\, dt \ , \notag\\
        \lesssim \  &\epsilon^2\
        \mathbb{E}_\delta[W](t_0,t_1) \ . \notag
\end{align}
To obtain the second to last line above, we have used the comparison
estimate \eqref{c_app3} from Appendix \ref{coercive_app} on the term
involving $w$.\\


\subsection*{Proof of estimate \eqref{R_red_detailed_lin2}}
This is virtually identical to the proof of
\eqref{R_red_detailed_lin1} in the previous paragraph. A simple
calculation using the notation from lines \eqref{schematic_not} and
line \eqref{F_rules} gives us the bound:
\begin{align}
    \big|\lambda\big(\partial_t(A_\lambda)\cdot w\big)\big|
        \lesssim \  \dlambda\lambda  |w|\, F^7_\lambda \ . \notag
\end{align}
The proof now follows line for line from the calculations performed
above.\\


\subsection*{Proof of estimate \eqref{R_red_detailed_lin3}}
We again use the formulas on line
\eqref{A_schematic} which give us:
\begin{equation}
        \big|\big(\partial_t(A_\lambda)\cdot w\big)\big|
        \lesssim \  \dlambda |w|\, F^7_\lambda \ . \notag
\end{equation}
Substituting this last line in the left hand side of
\eqref{R_red_detailed_lin3}, we have the following chain of
inequalities where the second to last line involves the bound
\eqref{c_app3}:
\begin{align}
    \sup_{t_0\leqslant s \leqslant t_1}\ \int_{\mathbb{R}^+} \
    \lambda^{-1}\frac{(\lambda r)^\delta}{1+r^\delta} |R_1|^2\ rdr
     \ &\lesssim \
    \sup_{t_0\leqslant s \leqslant t_1}\ \int_{\mathbb{R}^+}
    \ \frac{\dlambda^2}{\lambda^5}
    \frac{(\lambda r)^\delta}{1+r^\delta}\, \frac{w^2}{ r^4}\,
        \big(r^{4}F^{14}\big)_\lambda\ rdr\, dt \ , \notag\\
    &\lesssim \
    \epsilon^2\ \sup_{t_0\leqslant s \leqslant t_1}\
    \int_{\mathbb{R}^+} \ \lambda^{-1}
    \frac{(\lambda r)^\delta}{(1+r)^\delta}\,
    \frac{w^2}{ r^4}\ rdr\, dt \ , \notag\\
    &\lesssim \
    \epsilon^2\ \sup_{t_0\leqslant s \leqslant t_1}\
    \int_{\mathbb{R}^+} \ \lambda^{-1}
    \frac{(\lambda r)^\delta}{(1+r)^\delta}\,
    \frac{W^2}{ r^2}\ rdr\, dt \ , \notag\\
    &\lesssim \ \epsilon^2\ \mathbb{E}_\delta[W](t_0,t_1) \ . \notag
\end{align}\ret


\subsection*{Proof of estimate \eqref{R_red_detailed_lin4}}
We estimate the $R_2$ term from the line \eqref{detailed_lin}. By
using the second derivative identity from line \eqref{A_schematic},
the estimate \eqref{ddlambda_bound}, and then the bound
\eqref{orbit_bound} we have the following pointwise estimate:
\begin{align}
    \big|\partial_t^2(A_\lambda)\cdot w \big| \ &\lesssim \
    \big(\dlambda^2\lambda^{-1} + \epsilon^2\lambda^2\big)
    F_\lambda^7 \cdot|w| \ , \notag\\
    &\lesssim \ \epsilon^2 \lambda^3 F_\lambda^7 \cdot |w| \ . \notag
\end{align}
We now substitute this estimate into the left hand side of
\eqref{R_red_detailed_lin4} which allows us to estimate:
\begin{align}
    \int_{t_0}^{t_1} \int_{\mathbb{R}^+}\ \lambda^{-1}
    (\lambda r)^\delta |R_2|^2\ r^2dr\, ds
    \ &\lesssim \
    \epsilon^4\ \int_{t_0}^{t_1} \int_{\mathbb{R}^+}\
    \lambda^5(\lambda r)^\delta
    |w|^2\cdot F^{14}_\lambda
    \ r^2 dr\, ds \ , \notag \\
    &\lesssim \
    \epsilon^4\ \int_{t_0}^{t_1} \int_{\mathbb{R}^+}\
     \lambda^{-1}(\lambda r)^\delta
    \frac{w^2}{r^5}\cdot (r^{6}F^{14})_\lambda
    \ rdr\, ds \ , \notag \\
    &\lesssim \
    \epsilon^4\ \int_{t_0}^{t_1} \int_{\mathbb{R}^+}\ \lambda^{-1}
    \frac{(\lambda r)^\delta}{(1+r)^\delta}
    \frac{w^2}{r^5} \ rdr\, ds \ , \notag\\
    &\lesssim \
    \epsilon^4\ \int_{t_0}^{t_1} \int_{\mathbb{R}^+}\ \lambda^{-1}
    \frac{(\lambda r)^\delta}{(1+r)^\delta}
    \frac{W^2}{r^3} \ rdr\, ds \ , \notag\\
    &\lesssim \ \epsilon^4\ \mathbb{E}_\delta[W](t_0,t_1) \ . \notag
\end{align}\ret


\subsection*{Proof of estimate \eqref{R_red_detailed_lin5}}
First of all, using the formula for the nonlinearity
$\mathcal{N}(u)$ given on line \eqref{N_breakdown}, and by making
use of the formula \eqref{A_ops} for the operator $A_\lambda$ as
well as the formula \eqref{w0_def} for $w_0$, we easily have the
pointwise bound:
\begin{equation}
    \big| A_\lambda \mathcal{N}(u)\big|
    \ \lesssim \ \dlambda^4\lambda^{-5}F_\lambda^7
    + \frac{(|u|+|w_0|)\cdot |w|}{r^3} +
    \frac{(|u|+|w_0|)\cdot|W|}{r^2} \ . \label{nonlin_bound}
\end{equation}
We will deal with the first term on the right hand side above by itself.
The other two terms can be handled together.\\

We now substitute the first term on the right hand side of
\eqref{nonlin_bound} for $R_3$ on the left hand side of
\eqref{R_red_detailed_lin5}. Doing this we are left with estimating:
\begin{align}
    \int_{t_0}^{t_1}\int_{\mathbb{R}^+}\ \dlambda^8\lambda^{-11}
    (\lambda r)^\delta F^{12}_\lambda\ r^2dr\, ds
    \ &\lesssim \ \int_{t_0}^{t_1}\ \dlambda^8\lambda^{-14}(s)\ ds \ , \notag\\
    &\lesssim \ \sup_{t_0\leqslant s\leqslant t_1}
    \dlambda^4\lambda^{-7}(s)\cdot\int_{t_0}^{t_1}\
    \dlambda^4\lambda^{-7}(s)\ ds \ , \notag\\
    &\lesssim \ \epsilon \sup_{t_0\leqslant s\leqslant t_1}
    \dlambda^4\lambda^{-7}(s) \ . \notag
\end{align}
This proves the estimate \eqref{R_red_detailed_lin5} for the $w_0$
portion of $R_3$.\\

It remains to deal with \eqref{R_red_detailed_lin5} for the last two
terms on the right hand side of \eqref{nonlin_bound}. Upon
substitution of these into the right hand side of
\eqref{R_red_detailed_lin5} we have that:
\begin{align}
    &\int_ {t_0}^{t_1}\int_{\mathbb{R}^+}\ \lambda^{-1}(\lambda r)^\delta
    \left (|u|^2+|w_0|^2\right)\cdot\Big[ \frac{w^2}{r^6} +
    \frac{W^2}{r^4} \Big]\ r^2dr\, ds \ , \notag\\
    \lesssim \
    &\int_{t_0}^{t_1} \int_{\mathbb{R}^+}\
    \lambda^{-1}\frac{(\lambda r)^\delta}{(1+r)^\delta}
    \Big[\frac{w^2}{r^5} + \frac{W^2}{r^3}\Big]
    \cdot(1+r)^\delta \left (|u|^2+|w_0|^2\right)
    \ rdr\, ds \ , \notag \\
    \lesssim \ &\left (\sup_{\substack{r\\ t_0\leqslant s
     \leqslant t_1}} (1+r)^\delta|u|^2+\epsilon^2\right)\cdot
    \int_{t_0}^{t_1} \int_{\mathbb{R}^+}\ \lambda^{-1}
    \frac{(\lambda r)^\delta}{(1+r)^\delta}
    \Big[\frac{w^2}{r^5} + \frac{W^2}{r^3}\Big]\
    \ rdr\, ds \ , \notag\\
    \lesssim \ &t_1^\delta\epsilon^2\
    \int_{t_0}^{t_1} \int_{\mathbb{R}^+}\ \lambda^{-1}
    \frac{(\lambda r)^\delta}{(1+r)^\delta}
    \frac{W^2}{r^3}\ rdr\, ds \ , \notag\\
    \lesssim \ &t_1^\delta\epsilon^2\
    \mathbb{E}_\delta[W](t_0,t_1) \ . \notag
\end{align}
Notice that in the above estimates we have made crucial use of the
special pointwise estimate \eqref{extra_decay} proved below.
This is the only place in the paper which requires the extra decay of the initial data. This completes
our proof of the estimate \eqref{R_red_detailed_lin5}, and thus our
proof of Proposition \ref{eng_w_prop}.
\end{proof}\ret

\ret

\subsection{A Simple Decay Estimate}
In this subsection we will prove the rough decay estimate:
\begin{equation}
    \sup_r \ (1+r)^\delta|u|^2 \ \lesssim \ t^\delta\epsilon^2
    \ , \label{extra_decay}
\end{equation}
That is, our aim is to show that the reduced field
quantity $u$ enjoys some amount of \emph{pointwise} decay outside of
a sufficiently large cone centered at the space-time origin $t=0$
and $r=0$. \\

\begin{lem}[Decay of $u$ at space-like infinity]
Let $u=\phi-I_\lambda$ be the reduced field quantity as defined in
Lemma \ref{orbit_lem}, which in addition satisfies the initial
conditions of Theorem \ref{mod_thm}. In particular, $u$ is a
solution to the equation \eqref{lin_eq} with initial data
\eqref{special_C_data}--\eqref{special_smallness} and obeys the
estimate \eqref{energy_size} on the time interval $[0,T]$ where $\phi$
exists and remains smooth.  Then $u$ also obeys the following
stronger energy type estimate for any time $t\in [0,T]$, for which in
addition $\lambda\geqslant 1$:
\begin{equation}
    \int_{2t\leqslant r }\ r^2\, \left[ (\partial_t\phi)^2
    + (\partial_r u)^2 + \frac{k^2}{r^2}u^2\right]\
    rdr \ \lesssim \ \epsilon^2 . \label{decay_energy}
\end{equation}
\end{lem}\ret

\begin{rem}
 To transform estimate
\eqref{decay_energy} into an $L^\infty$ bound can be done in
an elementary way by applying the Poincar\'e type estimate
\eqref{energy_Poincare} to the quantity $r\chi_{3t\leqslant r}u$,
and then using the bound \eqref{decay_energy} to estimate the
resulting right hand side. Here $\chi_{3t\leqslant r}$ is a smooth
cutoff onto the region where $3t\leqslant r$ which satisfies the
homogeneity bound $|\chi_{3t\leqslant r}'|\lesssim r^{-1}$.
Therefore, we arrive at the estimate:\ \ $\sup_{3t\leqslant r} r|u|
\ \lesssim \ \epsilon$. By combining this with the pointwise bound $|u|\leqslant \epsilon$
which holds everywhere, we easily have \eqref{extra_decay}
whenever $\delta\leqslant 1$.
\end{rem}\ret

\begin{proof}[Proof of estimate \eqref{decay_energy}]
The proof is an integration by parts argument with a certain
multiplier. We denote by $\alpha(y)$ a smooth
increasing function, supported where $10\leqslant y$, satisfying
$\alpha'\leqslant 3 y$ and the homogeneity bound $y^{-1}\alpha \leqslant \alpha'$. The desired
result will now follow from computing the left hand side of the
identity:
\begin{equation}
    \int_0^t\int_{\mathbb{R}^+}\ \Big[\partial_t^2\phi +
    H_\lambda u -\mathcal{N}(u) \Big]
    \partial_t\phi\cdot \alpha(r-2s)\ rdr\, ds \ = \ 0
    \ , \label{decay_int_iden}
\end{equation}
where $\mathcal{N}(u)=\hbox{R.H.S.}\eqref{lin_eq}$. Also, we will
write the Hamiltonian from line \eqref{the_fact} as
$H_\lambda=-\partial_r^2 - r^{-1}\partial_r +Q_\lambda$. Notice that we have
$Q_\lambda\geqslant ck^2 r^{-2}$ on the support of $\alpha(r-2s)$.\\

Using the factorization \eqref{the_fact} as well as the fact that
$A_\lambda(\dot{I}_\lambda)=0$, thanks to
\eqref{J_def}--\eqref{lin_bog}, we may transform
\eqref{decay_int_iden} into the identity:
\begin{align}
    &- \frac{1}{2}\ \int_{\mathbb{R}^+}\
    \Big[(\partial_t\phi)^2 + (\partial_r u)^2 +
    Q_\lambda\, u^2\Big]\cdot\alpha\  rdr\
    \Bigg|_0^t \ , \label{decay_int_iden_comp}\\
    = \ &\int_0^t\int_{\mathbb{R}^+}\ \Big[(\partial_t\phi)^2 + \partial_t \phi \partial_r u +
    (\partial_r u)^2 + Q_\lambda\, u^2
    \Big]\cdot\alpha'(r-2s)\ rdr\, ds  \notag \\
    &\ \ \ \ \
    - \, k \int_0^t\int_{\mathbb{R}^+}\  \frac{u}{r}\cos(I_\lambda)\,
    \dot{I}_\lambda\,
    \cdot \alpha' \ rdr\, ds \ - \ \frac{1}{2} \int_0^t\int_{\mathbb{R}^+}\
    \dot{Q}_\lambda\, u^2\cdot\alpha \ rdr\, ds \ \notag\\
    &\ \ \ \ \ \ \ \ \ \ \ \ \ \ \ \ -  \
    \int_0^t\int_{\mathbb{R}^+}\ \mathcal{N}(u)\,
    \partial_t\phi\cdot \alpha \ rdr\, ds  \ , \notag\\
    = \ &T_1 + T_2 + T_3 + T_4\ . \notag
\end{align}
The proof will be complete once we show that the terms on the
right hand side of this last expression are either non-negative or
are bounded in absolute value by $C\epsilon^2$. In fact, it is more or less
immediate that we have:
\begin{align}
    |T_2| \ &\lesssim \ \epsilon^2 \ ,
    &T_1+T_3+T_4 \ &\geqslant \ 0 \ . \label{T_ests}
\end{align}
The first estimate above is a consequence of the Cauchy-Schwartz
inequality, the orbital stability bound \eqref{energy_size}, and the
following fixed time estimate valid for $1\leqslant \lambda$:
\begin{equation}
    \lp{\dot{I}_\lambda\cdot\alpha'}{L^2(rdr)} \ \lesssim
    \ \epsilon  \ . \notag
\end{equation}
This last line uses our assumptions that $4\leqslant k$ and
$1\leqslant \lambda$. Specifically,  the ODE bound \eqref{orbit_bound}
$|\dot{\lambda}\lambda^{-1}|\lesssim\epsilon\lambda$ (from
\eqref{orbit_bound} above) and a simple calculation, using the
assumption that $1\leqslant \lambda$ and involving lines
\eqref{J_def} and \eqref{I_form1}, give us the bound
$|\dot{I}_\lambda\cdot\alpha'|\lesssim |r \dot{I}|\lesssim \epsilon (1+r)^{-3}$.\\

The second bound on line \eqref{T_ests} will follow from the
estimate $|T_3+T_4|\lesssim \epsilon T_1$.
The desired result is then a consequence of the homogeneity property
of $\alpha$ and bounds:
\begin{align}
    |\mathcal{N}(u)| \ &\lesssim \ \epsilon\frac{|u|}{r^2} \ ,
    &|\dot{Q}_\lambda| \ \lesssim \ \epsilon\frac{1}{r^3} \ . \notag
\end{align}
The first bound above is a simple consequence of the orbital
stability estimate \eqref{energy_size} together with
\eqref{energy_Poincare}. The second bound follows again from the
estimate $|\dot{\lambda}\lambda^{-1}|\lesssim\epsilon\lambda$ of \eqref{orbit_bound}
and the explicit formulas on lines \eqref{I_form1}--\eqref{I_form2}.\\

The estimate \eqref{decay_energy} now follows from the form of
the left hand side of \eqref{decay_int_iden_comp} and  the smallness
condition \eqref{energy_size}.
\end{proof}

\ret\ret

\ret\ret

\appendix

\section{Computation of the constant $C_*$}\label{C_comp_app}

The purpose of the appendix is to derive an explicit formula $C_*=0$
for the special constant  $C_*$ which appeared on line \eqref{C_exp}.
Here we have written $J=J_1$ according to previous notation. In what
follows we shall also denote $I=I_1=I^k$. Rescaling we have
that:
\begin{equation}
    C_* \ = \ T_1 + T_2 + T_3 \ , \notag
\end{equation}
where:
\begin{align}
    T_1 \ &= \ -\ k^2 \Big\langle \frac{\left(a J + b r^2 J\right)^2}
    {r^2}  \ , \ \sin(2I)\cdot J \Big\rangle  \ , \notag\\
    T_2 \ &= \  \ \Big\langle  a J + br^2 J \ , \
    r\partial_rJ \Big\rangle \ , \notag\\
    T_3 \ &= \ - \  \Big\langle r\partial_r \left(a J +
    br^2 J\right) \ , \ r\partial_r J \Big\rangle \ . \notag
\end{align}
Recall that the constants $a$ and $b$ are given on line
\eqref{ab_defs}. Using now the identity $\sin(2I) = 2\sin(I)\cos(I)$
as well as \eqref{J_def} and \eqref{B_eq}, we have that:
\begin{align}
    T_2 \ &= \ -\, k^2 \int_{\mathbb{R}^+}\  \left(a J + b r^2 J\right)^2
    \partial_r\big( \sin^2(I) \big)\ dr \ , \notag\\
    &= \ -\, \int_{\mathbb{R}^+}\  \left(a  + b r^2 \right)^2 J^2
    \partial_r\big( J^2 \big)\ dr \ , \notag\\
    &= \ 2ab\int_{\mathbb{R}^+}\ J^4 \ rdr \ + \
    2b^2 \int_{\mathbb{R}^+}\ J^4 \ r^3dr \ , \label{T1_comp}
\end{align}
To obtain the last line above, we have used the expansion $J^4 =
k^2J^2\sin^2(I)$, the Pythagorean identity, and the definitions of
$a,b$.\\

We now move on the term $T_2$ above. Here we have directly that:
\begin{align}
        T_2 \ &= \   \int_{\mathbb{R}^+}\ (a J + br^2 J)
        \partial_rJ \ r^2dr \ , \notag\\
        &= \ -\, a\int_{\mathbb{R}^+}\  J^2\ rdr \  -
        2b\int_{\mathbb{R}^+}\ J^2 \ r^3dr \ , \notag\\
        &= \ -\, \frac{1}{4}\int_{\mathbb{R}^+}\
        J^2 \ r^3dr \ . \label{T2_comp}
\end{align}\ret

Finally, we compute that:
\begin{align}
    T_3 \ &= \ - \   a\int_{\mathbb{R}^+}\ (r\partial_r J)^2\ rdr \ - \
    b\int_{\mathbb{R}^+}\  (r\partial_r J)^2 \ r^3dr \ - \
    2b\int_{\mathbb{R}^+}\ J \partial_r J \ r^4dr \ , \notag\\
    &= \ - \   ak^2\int_{\mathbb{R}^+}\ J^2\cos^2(I) \ rdr \ - \
    bk^2 \int_{\mathbb{R}^+}\ J^2\cos^2(I) \ r^3dr \ + \
    4b\int_{\mathbb{R}^+}\ J^2 \ r^3dr \ , \notag\\
    &= \  a\int_{\mathbb{R}^+}\ J^4 \ rdr \ + \
    b \int_{\mathbb{R}^+}\ J^4 \ r^3dr \ + \
    4b\int_{\mathbb{R}^+}\ J^2 \ r^3dr \ .
    \label{T3_comp}
\end{align}\ret

We now add together lines \eqref{T1_comp}--\eqref{T3_comp} into the
single formula:
\begin{equation}
    C_* \ = \   \frac{3}{2}a\int_{\mathbb{R}^+}\ J^4 \ rdr \ + \
    \frac{3}{2}b  \int_{\mathbb{R}^+}\ J^4 \ r^3dr \ + \
    3b\int_{\mathbb{R}^+}\ J^2 \ r^3dr \ . \label{full_C*}
\end{equation}
 It remains to compute the
first two integrals in this last expression.
\begin{align}
    \int_{\mathbb{R}^+}\ J^4 \ rdr \ &= \
    -k \int_{\mathbb{R}^+}\ r\partial_r\big(\cos(I)\big)J^2
    \ rdr \ , \notag\\
    &= \ 2k\int_{\mathbb{R}^+}\ \cos(I)J^2 \ rdr
    + 2k^2\int_{\mathbb{R}^+}\ \cos^2(I)J^2 \ rdr
    \ , \notag\\
    &= \ k^2\int_{\mathbb{R}^+}\ r\partial_r\big(\sin^2(I)\big) \ rdr
    + 2k^2\int_{\mathbb{R}^+}\ J^2 \ rdr
    -2\int_{\mathbb{R}^+}\ J^4 \ rdr \ , \notag\\
    &= \ (2k^2-2)\int_{\mathbb{R}^+}\ J^2 \ rdr -
    2\int_{\mathbb{R}^+}\ J^4 \ rdr \ . \notag
\end{align}
An almost identical calculation also shows that:
\begin{equation}
    \int_{\mathbb{R}^+}\ J^4 \ r^3dr \ = \
    (2k^2-5)\int_{\mathbb{R}^+}\ J^2 \ r^3dr -
    2\int_{\mathbb{R}^+}\ J^4 \ r^3dr \ . \notag
\end{equation}
Therefore, recalling the definition of $a$ and $b$, these last two
calculations together give:
\begin{equation}
    a\int_{\mathbb{R}^+}\ J^4 \ rdr \ + \
    b  \int_{\mathbb{R}^+}\ J^4 \ r^3dr \ = \
    -\frac{1}{3}b\int_{\mathbb{R}^+}\ J^2 \ r^3dr \ . \notag
\end{equation}
Inserting the last line into \eqref{full_C*} and using that
$b=\frac{1}{4}$, we have $C_*=0$ as desired.

\ret\ret

\section{A general functional analysis lemma}\label{coercive_app}

In this appendix, we prove a general form of a coercive estimate we need throughout
the paper. This turns out to be more expedient,
because the required structure is simply a matter of compactness and
weak convergence in various weighted Sobolev spaces. The general
result which we propose to prove here is the following:\\

\begin{lem}[Coercive bounds for first order operators]\label{c_lem}
Let $B_\ell$ be a sequence of first order differential operators
with real smooth (but not necessarily bounded!) coefficients on the
half line $(0,\infty)$, continuously indexed (in the weighted $L^2$
space defined by the LHS of \eqref{sub_coerc} below)  by
$\ell\in[0,1]$ and such that the following sub-coercivity holds for
some continuously indexed (for functions in the norm
\eqref{calH_norm}) function $|f_\ell|\lesssim r^{-2-\gamma +
\sigma}(1+r)^{-2\sigma}$ with $\sigma>0$:
\begin{equation}
    \int_{\mathbb{R}^+}\ \frac{(B_\ell\psi)^2}{r^\gamma}\ rdr\ = \
    \int_{\mathbb{R}^+}\ \big[
    \frac{(\partial_r\psi)^2}{r^\gamma}   +
    h_\ell\frac{\psi^2}{r^{2+\gamma}}
    + f_\ell\psi^2\big]\ rdr \ , \label{sub_coerc}
\end{equation}
for any real valued function $\psi$ with finite norm:
\begin{equation}
    \lp{\psi}{\mathcal{H}^\gamma}^2 \ = \
    \int_{\mathbb{R}^+}\ \big[
    \frac{(\partial_r\psi)^2}{r^\gamma}
      + \frac{\psi^2}{r^{2+\gamma}}\big]
    rdr \ . \label{calH_norm}
\end{equation}
Here $0\leqslant \gamma$ is a fixed parameter, and $0 < C\leqslant
h_\ell$ is some strictly positive function. Then there exists a
universal constant, uniform in $\ell$, such that the following bound holds:
\begin{equation}
    \lp{\psi}{\mathcal{H}^\gamma}^2 \ \lesssim \
    \ \int_{\mathbb{R}^+}\ \frac{(B_\ell\psi)^2}{r^\gamma}\
    rdr \ , \label{coercive_bound}
\end{equation}
for any real $\mathcal{H}^\gamma$ function $\psi$ which also
satisfies:
\begin{equation}
    \int_{\mathbb{R}^+}\ \psi \cdot J^\ell
    \ m_\ell \, rdr \ = \  0 \ ,
    \label{groundstate_orth}
\end{equation}
for some positive weight function $0 < m_\ell$. Here the function
$J^\ell$ is the (nontrivial) ``ground-state'' given by $B_\ell
J^\ell=0$, and we are assuming $m_\ell J^\ell\in
(\mathcal{H}^\gamma)^*$ depends continuously on $\ell$.
\end{lem}\ret

\begin{proof}[Proof of Lemma \ref{c_lem}]
The proof is based on a contradiction argument centered around weak
convergence. Suppose that the estimate \eqref{coercive_bound} was not true.
 Then there would exist a sequence of $\psi_n$
and $\ell_n$ such that:
\begin{equation}
    \int_{\mathbb{R}^+}\ \frac{(B_{\ell_n}\psi_n)^2}{r^\gamma}\
    rdr\
    \leqslant \ c_n\ \lp{\psi_n}{\mathcal{H}^\gamma}^2
    \ , \label{cont_bnds}
\end{equation}
where $c_n\to 0$ is some sequence of constants. We assume that this
sequence is normalized so that
$\lp{\psi_n}{\mathcal{H}^\gamma}=1$. The space
$\mathcal{H}^\gamma$ is a Hilbert space (with an obvious scalar product)
defined as a closure of
$C_0^\infty(\mathbb {R}_+)$ functions in the $\mathcal{H}^\gamma$ norm.
Therefore, we can choose  a subsequence $\psi_{n_k}$
which converges weakly in $\mathcal {H}^\gamma$ to
$\psi_\infty\in\mathcal{H}^\gamma$. Furthermore, we may assume
(by perhaps taking another subsequence) that
$\ell_{n_k}\to\ell_\infty$ for some $\ell_\infty\in[0,1]$. We now
use $\psi_n$ and $\ell_n$ to denote this subsequence.
Also, note that by Cauchy-Schwartz the unit
normalization implies that
$\lp{\psi_\infty}{\mathcal{H}^\gamma}\leqslant 1$.\\

By the continuity of the $J^\ell$, and the uniform boundedness of
the $\psi_n$, we have from the identity \eqref{groundstate_orth}
that the limiting function satisfies:
\begin{equation}
    \int_{\mathbb{R}^+}\ \psi_\infty \cdot
    J^{\ell_\infty} \ m_{\ell_\infty}\,
    rdr \ = \  0 \ .    \label{false_orth}
\end{equation}
Therefore, since the ``ground-state'' $J^{\ell_\infty}$ is unique
(it satisfies a first order ODE) and the measure
$m_{\ell_\infty}rdr$ is strictly positive on $(0,\infty)$, we will
have a contradiction if we can establish that $\psi_\infty$ is
nontrivial. This contradiction would come from again invoking
uniform boundedness and the continuity of the operators $B_\ell$
which implies that:
\begin{equation}
    \int_{\mathbb{R}^+}\ \frac{(B_{\ell_\infty}\psi_n)^2}{r^\gamma}\
    rdr \ \to \ 0 \ , \notag
\end{equation}
so that $\psi_\infty$ is a weak solution of
$B_{\ell_\infty}\psi_\infty=0$, and hence a smooth solution via
ODE regularity, thus violating  uniqueness as
\eqref{false_orth} implies $\psi_\infty \neq \beta J^{\ell_\infty}$
for any constant $\beta\neq 0$.\\

To show that $\psi_\infty$ is nontrivial, we make crucial use of
the sub-coercivity condition \eqref{sub_coerc}. By the unit
normalization, the universal lower bound on $h_\ell$, and the
assumption that \eqref{cont_bnds}, we have that there exists a
universal lower bound to the limit:
\begin{equation}
    0 \ > \ \varlimsup_{n\to\infty}\
    \int_{\mathbb{R}^+}\ f_{\ell_n}\psi^2_n\ rdr \ . \notag
\end{equation}
Therefore we shall have that $\psi_\infty$ is not everywhere zero
if we can show that the sequence $f_{\ell_n}\psi^2_n$ converges
strongly in $L^1(rdr)$. This in turn follows from the universal bounds on
and continuity of $f_\ell$, and fact that $\psi_n$ converges
strongly in the weighted space:
\begin{equation}
     \int_{\mathbb{R}^+}\ r^{-2-\gamma + \sigma}(1+r)^{-2\sigma}
     \psi^2\ rdr
     \ = \ \lp{\psi}{\mathcal{L}^2_{\gamma,\sigma}}^2 \ . \notag
\end{equation}
This latter strong convergence is provided via uniform boundedness
and the compact inclusion $\mathcal{H}^\gamma\Subset
\mathcal{L}^2_{\gamma,\sigma}$ whenever $0 < \sigma$.
\end{proof}\ret

In practice, we will only need two special cases of the Lemma
\ref{c_lem} above. The first case is where $B_\ell\equiv A_1$ and $\gamma=0$, where
$A_1$ is the first order operator from line \eqref{A_ops} above. The
second cases are when we set $\ell=\lambda^{-1}$, with $1\leqslant
\lambda$ and:
\begin{equation}
    B_\ell \ = \ \big(1+(\lambda^{-1}r)\big)^{-\frac{\delta}{2}}
    \, A_1\, \big(1+(\lambda^{-1}r)\big)^{\frac{\delta}{2}} \ , \notag
\end{equation}
where in this case we set
$\psi=(1+(\lambda^{-1}r))^{-\frac{\delta}{2}}u$, as well as
$J^\ell=\big(1+(\lambda^{-1}r)\big)^{-\frac{\delta}{2}}J_1$.
Finally, in this case we set $m_\ell
=\big(1+(\lambda^{-1}r)\big)^\delta$. We apply this to $\gamma=2-\delta$ and $\gamma=3-\delta$.\\

In all of these cases we leave it to the reader to prove that the
condition \eqref{sub_coerc} holds (the continuity is obvious). This
is a simple matter of integration by parts, the fact that
$4\leqslant k$ (notice that this works for our range of $\gamma$, which is the main
thing to check here),
and also that we have chosen $\delta\ll 1$. However, we
do call the readers attention to an important and perhaps subtle point.
In order for the integration by parts to work, it is necessary to show
that a boundary term of the form\  $\lim_{r\to 0}r^{-\gamma}\psi^2$
vanishes for any function $\psi$ in the space $\mathcal{H}^\gamma$.
This follows from the finiteness of that norm, and the fact that the
Poincar\'e type estimate \eqref{energy_Poincare} above applied to
$r^{-\frac{\gamma}{2}}\psi$ implies that this function is
\emph{continuous} on the closed interval $[0,1]$.
This latter fact is perhaps a bit subtle,
and it is crucial for showing
$r^{-\frac{\gamma}{2}}\psi$ vanishes at $r=0$ via the finiteness of
the weighted $L^2$ norm (no derivative) contained in $\mathcal{H}^\gamma$.
Again, we leave the reader to check the details of all this.\\

Now, Applying the above result in these two cases and rescaling by $\lambda$,
we have that:\\

\begin{lem}[Applied version of Lemma \ref{c_lem}]
Suppose that $4\leqslant k$ and $A_\lambda$ is the operator defined
on line \eqref{A_ops}. Then if $u$ is a function which satisfies the
admissibility condition of Lemma \ref{admiss_lem}, and is such that
the energy from line \eqref{energy_size} is finite, and such that
the orthogonality relation holds:
\begin{equation}
    \int_{\mathbb{R}^+}\ u\cdot J_\lambda\ rdr \ = \ 0 \ , \notag
\end{equation}
then one had the following universal bounds whenever
$\delta\ll 1$ is small enough (and $1\leqslant \lambda$ in the last two cases):
\begin{align}
    \int_{\mathbb{R}^+}\ \big[(\partial_r u)^2
    + \frac{u^2}{r^2}\big]\
    rdr \ &\lesssim \ \int_{\mathbb{R}^+}\
    (A_\lambda u)^2\ rdr \ , \label{c_app1}\\
    \int_{\mathbb{R}^+}\ \frac{(\lambda r)^\delta}{(1+r)^\delta}
     \frac{u^2}{r^4}\
    rdr \ &\lesssim \ \int_{\mathbb{R}^+}\
    \frac{(\lambda r)^\delta}{(1+r)^\delta}
    \frac{(A_\lambda u)^2}{r^2}\ rdr \ , \label{c_app2}\\
    \int_{\mathbb{R}^+}\ \frac{(\lambda r)^\delta}{(1+r)^\delta}
     \frac{u^2}{r^5}\
    rdr \ &\lesssim \ \int_{\mathbb{R}^+}\
    \frac{(\lambda r)^\delta}{(1+r)^\delta}
    \frac{(A_\lambda u)^2}{r^3}\ rdr \ . \label{c_app3}
\end{align}
\end{lem}

\ret\ret


\end{document}